%% file: New_Star_critical_Ramsey_numbers_cycles_vs_K5_new.tex
\documentclass[12pt]{article}
\usepackage{pgf,tikz}
\usetikzlibrary{arrows}
\usepackage{amsmath}%
\usepackage{amsfonts}%
\usepackage{amssymb}%
\usepackage{graphicx}

\newtheorem{theorem}{Theorem}

\newtheorem{lemma}{Lemma}

\begin{document}
\date{}
\title{Star-critical Ramsey numbers for cycles versus the complete graph on 5 vertices}
\author{ Chula J. Jayawardene \\
Department of Mathematics\\
University of Colombo \\
Sri Lanka\\
email: c\_jayawardene@yahoo.com
 }

\maketitle
\begin{abstract} Let  $G$, $H$ and $K$ represent three graphs without loops or parallel edges and $n$ represent an integer. Given any red blue coloring of the edges of $G$, we say that $K \rightarrow (G,H)$, if there exists  red copy of $G$  in $K$  or a blue copy of $H$  in $K$.  Let $K_n$ represent a complete graph on $n$ vertices, $C_n$ a cycle on $n$ vertices and $S_n=K_{1,n}$ a star on $n+1$ vertices. The Ramsey number $r(G, H)$ is defined as $\min\{n \mid K_n\rightarrow (G,H)\}$.  Likewise, the star-critical Ramsey number $r_*(H, G)$ is defined $\min\{k \mid  K_{r(G,H)-1} \sqcup K_{1,k} \rightarrow (H, G) \}$. When $n >3$, in this paper we show that $r_*(C_n,K_5)=3n-1$ except $r_*(C_4,K_5)=13$. We also characterize all Ramsey critical $r(C_n,K_5)$  graphs.
\end{abstract}


\section{INTRODUCTION}
\noindent In its classical form, Ramsey's theorem ensures the existence of the Ramsey numbers $r(n,m)$  defined as $r(n,m)=r(K_n,K_m)$. The exact determination of these numbers progresses rapidly in difficulty from  the nearly trivial $r(3,3)=6$ to stubbornly resistant $r(5,5)$ [at present known to be between 43  and 48]. One new branches of classical Ramsey numbers namely \textit{star -critical Ramsey numbers} were introduced by Hook and Isaak in 2010 (see \cite{Ho,HoIs}).  Using the arrow notation,  $r_*(H, G)$ can be also be defined as the smallest positive integer $k$ such that $K_{n}\setminus K_{1,n-k-1} \rightarrow (H, G)$ where $n=r(G,H)$ (as $K_{n-1} \sqcup K_{1,k} =K_{n}\setminus K_{1,n-k-1}$). Obviously, the rational behind finding star critical Ramsey numbers is to get a more in-dept understanding of classical case. Recently many papers have tried to find star critical Ramsey numbers for such as fans vs. complete graphs and cycles versus $K_4$ (see \cite{BaSuBr,HaMaSe,Ho,HoIs,JaNaRa,YaYoRa}).

\vspace{6pt}
\noindent In this  paper we concentrate on finding $r_*(C_n, K_5)$. Most of the work related to this paper, originated with the attempt to solve Erd{\"o}s and Bondy conjecture  which states that $r(C_n,K_m)=(n-1)(m-1)+1$ for all $n\ge m \ge 3$, except when $n=m=3$ (see \cite{BoEr}). However, up to now, the conjecture is unsolved and has been only proved for $n\ge m$ and $n\le 7$ and some other special cases.

\section{NOTATION}

\noindent All graphs $G=(V,E)$ considered in this paper are finite graphs without loops and multiple edges. The order of the graph $G=(V,E)$ is denoted by $|V|$ and  the number of edges in the graph is denoted by $|E|$. For a graph $G=(V,E)$ the \textit{degree} of $v$ is defined as the cardinality of the set  of vertices adjacent to $v$. A set $I$ of $V(G)$, is said to be an \textit{independent set} of a graph $G$ if no pair of vertices of $I$ are connected by an edge in $G$. 

\vspace{6pt}
Suppose that a graph $G$ contains an $n$ cycle $(u_1, u_2, ... , u_{n}, u_1)$ and a vertex (say $y_1$) outside of the cycle such that $y_1$ is adjacent to  exactly two vertices (say $u_i$ and $u_j$) of the $n$  cycle. In such a situation, we say that $y_1$ is adjacent to two vertices of the $\{u_1, u_2, ... , u_{n}\}$  which are length $k$ apart where $k=\min \{ (i-j) \mod n,(j-i) \mod n \}$. Throughout the paper, in any graph which is colored by red and blue, we will denote the red edges by a unbroken line and the blue edges by a broken line. For a red/blue coloring of a graph $G$, and vertices $u,v\in V(G)$ such that $(u,v)\in E(G)$, we say that $u$ is a red (resp. blue) neighbor of $v$ if $(v,u)$ is colored red (resp. blue).  Any 2-coloring of $K_{r(G,H)-1}$ that does not contain a red $G$ or a blue $H$,  is called a \textit{critical coloring}.

\section{PROPERTIES OF $(C_4,K_5)$ RAMSEY CRITICAL GRAPHS}

In order to prove the main result of this paper, namely finding $r_*(C_n, K_5)$, we try to utilize the the critical graphs of $r(C_n,K_4)$ (see \cite{JaNaRa}) and the fact that  $r(C_n,K_5) = 4n - 3$ for $n \geq 4$ (see \cite{BoJaSh,Ra}). In addition, we use the following lemmas to arrive at the main result.

\vspace{6pt}
\noindent The following  four lemmas is a direct consequence of \cite{JaRo,JaRo1,JaNaRa}, written by Jayawardene et al.

\begin{lemma}
\label{l1}
\noindent (\cite{JaRo}, Lemma 2). If $G$ is a graph of order $N$ that contains no $C_m$ and the independent number is less than or equal to $n-1$ then the minimal degree is greater than or equal to $N-r(C_m,K_{n-1})$.
\end{lemma}

\vspace{10pt}
\begin{lemma}
\label{l2}
\noindent (\cite{JaRo1}, Lemma 5). Any $C_5$-free graph of order 11 with no independent set of 4 vertices is isomorphic to one of the graphs $R_{11,1}$, $...$, $R_{11,18}$ (see Figure 1) or  $R_{11,19} \cong 2K_4 \cup K_3$.

\begin{center}
\input{graphl21a.tex}
\input{graphl21b.tex}
\input{graphl22.tex}
\input{graphl23.tex}
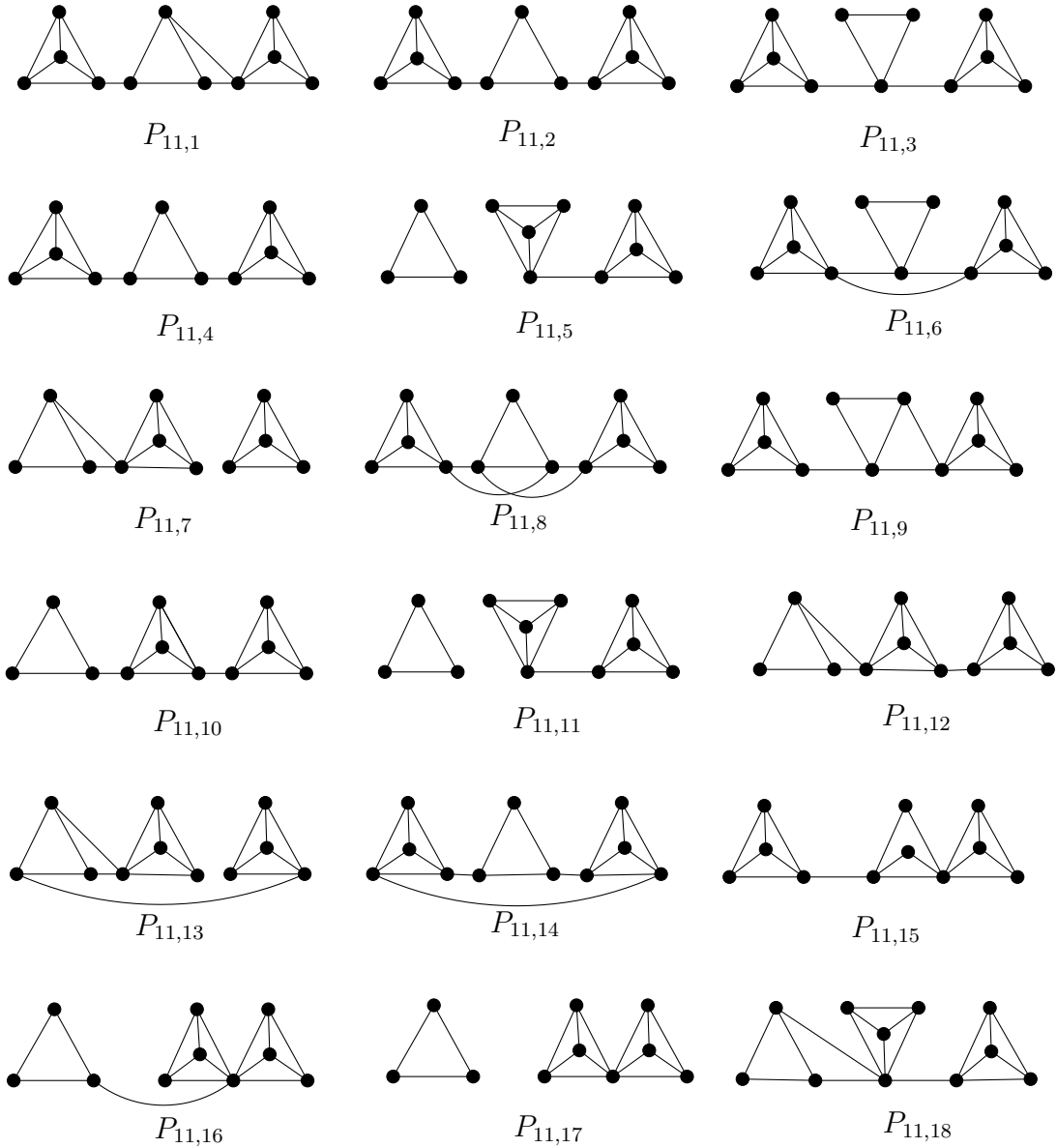
Figure 1. $R_{11,k},\ 1\leq k\leq 18$.
\end{center}
\end{lemma}

\vspace{10pt}
\begin{lemma}
\label{l3}
\noindent (\cite{JaRo1}, Lemma 4). Any $C_5$-free graph of order 12 with no independent set of 4 vertices is isomorphic to one of the graphs $R_{12,1}$, $R_{12,2}$, $R_{12,3}$, $R_{12,4}$, $R_{12,5}$ (see the following figure) or  $R_{12,6} \cong 3K_4$.

\begin{center}
\input{graphl2a.tex}
\input{graphl2b.tex}
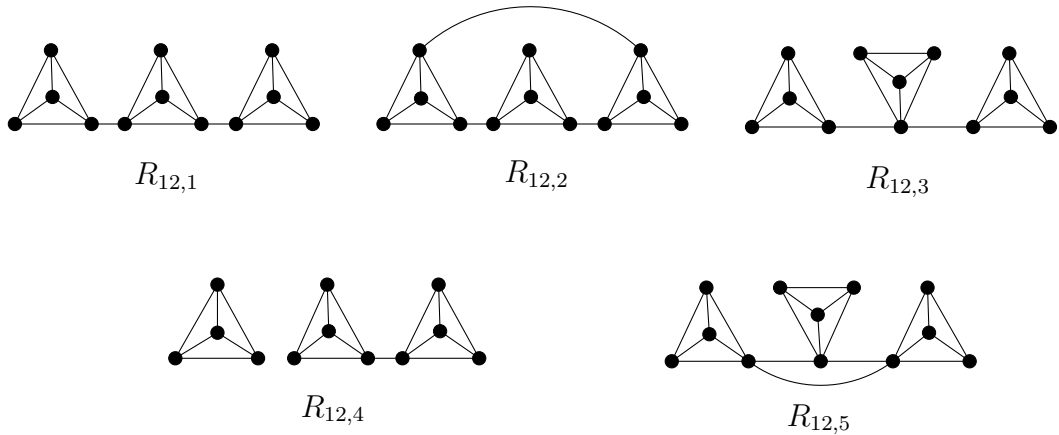
Figure 2. $R_{12,k},\ 1\leq k\leq 5$.
\end{center}
\end{lemma}

\begin{lemma}
\label{l4}
\noindent (\cite{JaNaRa}, Lemma 7). Any $C_6$-free graph of order 15 with no independent set of 4 vertices is isomorphic to one of the graphs five critical graphs corresponding to $K_{15}$, denoted by $R_{15,1}$, $R_{15,2}$, $R_{15,3}$, $R_{15,4}$ or $R_{15,5}$, where $R_{15,4} \cong 3K_{5}+e$ and $R_{15,5} \cong 3K_{5}$. The other three red graphs, namely $R_{15,1}$, $R_{15,2}$, $R_{15,3}$, are illustrated in Figure 3.

\begin{center}
\input{graph71.tex}  
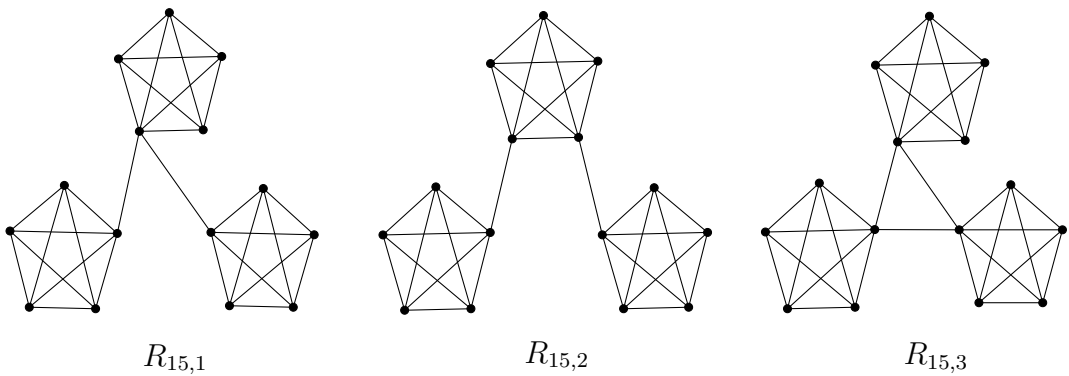
Figure 3. The graphs $R_{15,k},\ 1\leq k\leq 3$.
\label{r3n33}
\end{center}
\end{lemma}

\vspace{10pt}
\noindent The following  lemmas is a direct consequence of \cite{BoJaSh}, written by Bollab{\'a}s et al.

\begin{lemma}
\label{l5}
\noindent Suppose $G$ contains the cycle $(u_1, u_2, ... , u_{n-1}, u_1)$ of length $n - 1$ but
no cycle of length $n$. Let $Y = V(G) \setminus \{u_1, u_2, ... , u_{n-1}\}$. Then,

\vspace{10pt}
\noindent \textbf{(a)} No vertex $x \in Y$ is adjacent to two consecutive vertices on the cycle.

\vspace{6pt}
\noindent \textbf{(b)} If $x \in Y$ is adjacent to $u_i$ and $u_j$ then $u_{i+1}u_{j+1} \notin E(G)$.

\vspace{6pt}
\noindent \textbf{(c)} If $x \in Y$ is adjacent to $u_i$ and $u_j$ then no vertex $x' \in Y$ is adjacent to both
$u_{i+1}$ and $u_{j+2}$.

\vspace{2pt}

\noindent \textbf{(d)} Suppose $\alpha(G) = m - 1$ where $m \leq  \frac {n + 3}{2}$ and $\{x_1, x_2, ... ,x_{m -1} \} \subseteq Y$ is an $(m - 1)$-element independent set. Then no member of this set is adjacent to $m -2$ or more vertices on the cycle.
\end{lemma}

\begin{lemma}
\label{l6}
\noindent A $C_5$ -free graph of order $16$ with no independent set of 5 vertices contains a isomorphic copy of  $4K_{4}$.
\end{lemma}

\noindent {\bf Proof.}

\noindent Let $G$ be a $C_5$-free graph of order 16 with no independence set of 5 vertices. First note that $\delta(G) \geq 3$ by lemma \ref{l3} as $r(C_5, K_4) = 13$.

\vspace{14pt}

\noindent \textbf{Remark 1:} Suppose that $G$ contains a $K_4$, then any vertex outside of $K_4$ can be adjacent to at most one vertex of the $K_4$. Further any two adjacent vertices together can be adjacent to at most one vertex of a $K_4$. 

\vspace{8pt}
\noindent \textbf{Case 1:} The minimum degree is 3. 

\vspace{6pt}
\noindent Say $w$ is a vertex of minimum degree. Then $G[V(G)\setminus \overline{\Gamma(w)}]$ will satisfy the conditions of lemma \ref{l3}. Thus, using all possible cases and using remark 1, any two non adjacent vertices of $G[\Gamma(w)]$ will give rise to an independent set of size five containing them. Therefore, $G[\Gamma(w)] \cong K_3$. This will result in the required $4K_4$. 

\vspace{8pt}
\noindent \textbf{Case 2:} The minimum degree greater than or equal to 4.
 
\vspace{6pt} 
\noindent Say $w$ is a vertex of minimum degree. Then $G[V(G)\setminus V[\overline{\Gamma(w)]}]$ will satisfy the conditions of lemma \ref{l2}. Thus, using all possible cases and using remark 1 one sees that there are no possible extensions. 

\vspace{8pt}
\noindent \textbf{Case 3:} The minimum degree greater than or equal to 5. 

\vspace{6pt}
\noindent Say $w$ is a vertex of minimum degree. Then $G[\Gamma(w)]$ will contain a $C_3 \cup K_2$, $P_3 \cup K_2$,  $K_{1,4}$, $K_{1,3} \cup K_1$, $2K_2$ or at least two isolated vertices.  It is worth noting that,  we will get $C_5$ or else an independent set of size 5 directly in all cases, other than  when $G[\Gamma(w)]=C_3 \cup K_2$. In this case, let $H=G[V(G)\setminus \overline{\Gamma(w)}]$. Then, $H$ has 10 vertices. In order to avoid a $C_5$ each these three vertices of the $C_3$ in $G[\Gamma(w)]$ will each have to be adjacent to two  vertices of $H$ and these neighborhoods will have to be non overlapping. Select $v_1$, $v_2$ and $v_3$ from the three neighborhoods. Moreover, in order to avoid a $C_5$ the two vertices of the $K_2$ in $G[\Gamma(w)]$ will share a common neighbor in $H$, say $v$, as illustrated in the following figure.

\begin{center}
\input{graphl31.tex}
\end{center}
\begin{center}
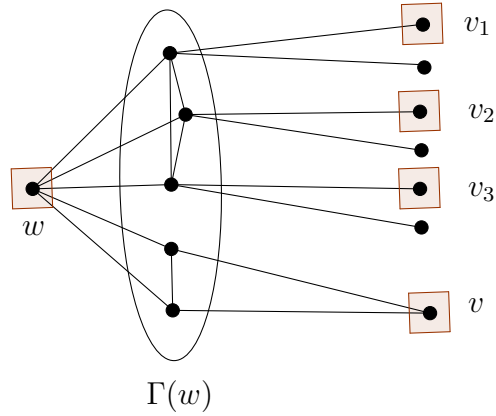
Figure 4. The only option left in case 3
\end{center}

\noindent Then in order to avoid a $C_5$, $\{w,v_1,v_2,v_3,v\}$ will have to be an independent set contrary to the assumption. 

\begin{lemma}
\label{l7}
\noindent A $C_6$ -free graph of order $20$ with no independent set of 5 vertices contains a isomorphic copy of  $4K_{5}$.
\end{lemma}

\noindent {\bf Proof.}

\noindent Let $G$ be a $C_5$-free graph of order 20 with no independence set of 5 vertices. First note that, by lemma \ref{l1}, $\delta(G) \geq 4$  as $r(C_6, K_4) = 16$.

\vspace{14pt}

\noindent \textbf{Remark 2:} Any vertex can be adjacent to one vertex of a disjoint $K_5$. Further any two adjacent vertices together can be  adjacent to at most one vertex of a $K_5$. 

\vspace{8pt}
\noindent \textbf{Case 1:} The minimum degree is 4. 

\vspace{6pt}
\noindent Say $w$ is a vertex of minimum degree. Then $G[V(G)\setminus \overline{\Gamma(w)}]$ will satisfy the conditions of lemma \ref{l4}. Thus, using all possible cases and using remark 2, it follows that any two non adjacent vertices of $G[\Gamma(w)]$ will give rise to an independent set of size five containing them. Therefore, $G[\Gamma(w)] \cong K_4$. This will result in the required $4K_5$. 

\vspace{8pt}
\noindent \textbf{Case 2:} The minimum degree greater than or equal to 5.
 
\vspace{6pt} 
\noindent Suppose that $G$ is a $C_6$ -free graph on $20$  with no independent set of 5 vertices. Then as $r(C_5, K_5)=17$ (see \cite{BoJaSh,Ra}) there exists a cycle $U=(u_l, u_2, ... , u_{5}, u_1)$ of length $5$. Let $X=\{u_l, u_2, ... , u_{5}\}$. Define $H=G[X^c]$  as the induced subgraph of $G$ not containing the vertices of the cycle and $H_1=G[X]$. Then, $|V(H)|=15$ and $|V(H_1)|=5$. 

\vspace{6pt}

\noindent  Suppose that there exists an independent set $Y$ in $H$ of size 4 consisting of the four vertices $y_1$, $y_2$, $y_3$ and $y_4$.  In order to avoid an independent set of size 5, each vertex of $X$ must be adjacent to at least one vertex of $Y$. Clearly, no vertex of $Y$  can be adjacent to three vertices of $X$ as two of these three adjacent vertices will have to be  consecutive vertices of the $C_5$ (which will result in a $C_6$ containing all the vertices of $X$).  Therefore, each vertex of $Y$ is adjacent to at most two vertices of $X$ and if they are adjacent to two vertices of $X$ they must be length 2 apart. As $\vert X \vert =5$, we get that without loss of generality  $y_1$  will have exactly two red neighbors in $X$ length 2 apart (say $u_1$ and $u_3$). Then without loss of generality $u_2$ is  adjacent to $y_2$ and in order to avoid a red $C_6$, $y_2$ cannot be adjacent to any other vertex of $X$. 

\vspace{6pt}
\noindent Next without loss of generality $\{y_3, y_4\}$ will have either two vertices or one vertex or no vertices that are adjacent to two vertices of $X$. Furthermore, any vertex of $\{y_3, y_4\}$ adjacent to 2 vertices of $X$ must be adjacent to either $u_3$ and $u_5$ or $u_1$ and $u_4$ or $u_1$ and $u_3$. Moreover, if $y_3$ is adjacent to $u_3,u_5$ and $y_4$ is adjacent to $u_4,u_1$ then  $(u_1,u_5,y_3,u_3,u_4,y_4,u_1)$ will be a cycle of length $6$, contrary to our assumption. Also if, $y_3$ is adjacent to $u_1,u_3$ we will get that  $y_4$ will be forced to be adjacent to $u_4,u_5$ and then $(u_1,u_2,u_3,u_4,y_4,u_5,u_1)$ will be a cycle of length $6$, contrary to our assumption. Therefore, without loss of generality, by symmetry, we are left with the two possibilities where $y_3$ is adjacent to $u_5$ and $y_4$ is adjacent to $u_4$ or else $y_3$ is adjacent to $u_3,u_5$ and $y_4$ is adjacent to $u_4$. In the first option, in order to avoid an independent set of size 5, $u_1$ must be adjacent to $u_3$. As $y_2$  can not be adjacent to any other vertex of $X$ it must be adjacent to 4 vertices (as minimum degree of $G$ is greater than 5) outside of $X \cup Y$ (say $w_1$,$w_2$,$w_3$ and $w_4$), as illustrated in the following figure.

\begin{center}
\input{graphl41.tex} 
\end{center}
\begin{center} 
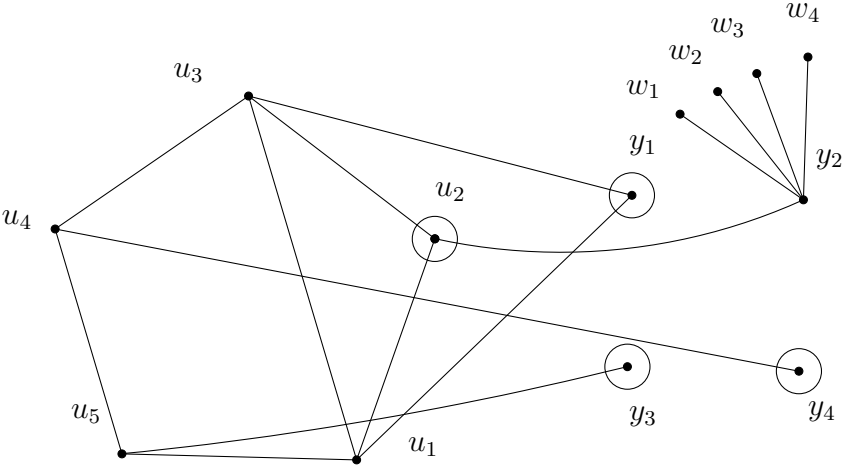
Figure 5. The first possibility
\end{center}

\noindent In order to avoid a independent set of size 5 consisting of $\{w_1,u_2,y_1, y_3, y_4\}$ we get $w_1$ must be adjacent to $u_2$. Similarly, we can conclude that all other vertices of $W=\{w_1,w_2,w_3, w_4\}$ too have to be adjacent to $u_2$. But then $G[W]=K_4$ as otherwise any two non adjacent vertices of $W$ together with $\{y_1,y_3,y_4\}$ will form an independent set of 5 vertices, contrary to the assumption. This will give us a  $C_6$ containing all vertices $W$,  contrary to our assumption.

\begin{center}
\input{graphl42.tex} 
\end{center}
\begin{center} 
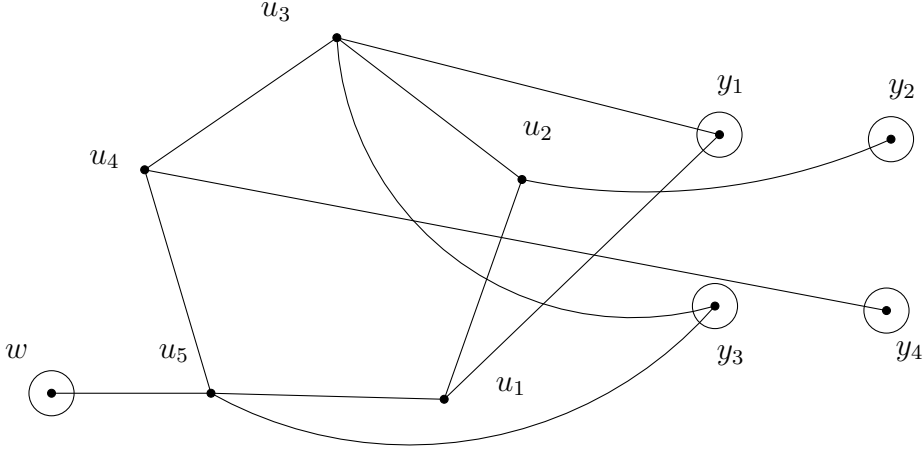
Figure 6: The second possibility
\end{center}

\noindent In the second option, as $u_5$ can not be adjacent to $u_2$ (may be adjacent to $u_3$), it will have to be adjacent to some other vertex (since the  minimum degree of $G$ is greater than 5) outside of $X \cup Y$ (say $w$), as illustrated in the above figure. But then, $w$ can not be adjacent to any vertex of $Y$  as any such occurrence will lead to some $C_6$. Therefore, $Y \cup \{w\}$ will be an independent set of size 5, contrary to our assumption.

\vspace{6pt}
\noindent Thus, we can conclude that the initial assumption is false. That is, there is no independent set of order 4 in $H$. By Lemma \ref{l4}, we can conclude that $H$ is equal to one of the five graphs  $R_{15,1}$, $R_{15,2}$, $R_{15,3}$, $R_{15,4}$ or $R_{15,5}$. Now consider any two vertices of $U$, say $u$ and $v$, and suppose that $\left\{u,v\right\}\not \in E(G)$. Since there is no $C_6$ in $G$, each of the vertices $u$ and $v$ must be adjacent to at most one vertex of each copy of $K_{5}$ in $H$. Therefore,  we can select  vertex $x_1$ in the first $K_{5}$, vertex $x_2$ in the second $K_{5}$ and vertex $x_3$ in the third $K_{5}$ such that $x_1$, $x_2$ and $x_3$ are independent and not adjacent to $u$ or $v$. This gives us that $\{u,v,x_1,x_2,x_3\}$ is an independent set of order 5, a contradiction.   Therefore, $\left\{u,v\right\} \in E(G)$. Since $u,v$ are arbitrary vertices in $U$, we can conclude that $X$ induces a $K_{5}$ as required.\hfill$\square$ 

\begin{lemma}
\label{l8}
\noindent A $C_n$ -free graph (where $n \geq 7$) of order $4(n-1)$ with no independent set of 5 vertices  contains a isomorphic copy of  $4K_{n-1}$.
\end{lemma}

\noindent {\bf Proof.}

\noindent Suppose that $G$ is a $C_n$ -free graph on $4(n-1)$  with no independent set of 5 vertices. Then as $r(C_{n-1}, K_5)=4n-7$ (see \cite{BoJaSh,Ra}) there exists a cycle $(u_l, u_2, ... , u_{n-1}, u_1)$ of length $n - 1$. Let $X=\{u_l, u_2, ... , u_{n-1}\}$. Define $H=G[X^c]$  as the induced subgraph of $G$ not containing the vertices of the cycle and $H_1=G[X]$. Then, $|V(H)|=3(n-1)$ and $|V(H_1)|=n-1$. 

\vspace{6pt}

\noindent Suppose that there exists an independent set $Y$ in $H$ of size 4 consisting of the four vertices $y_1$, $y_2$, $y_3$ and $y_4$. That is, $\alpha(G)=4$. From  lemma \ref{l5}(d) (as $5\leq \frac{n+3}{2}$), it follows that every vertex $y$ is adjacent  to at most two vertex of the cycle $C_{n-1}$.

\vspace{8pt}
\noindent \textbf{Case 1: $n \ge 10$}

\vspace{4pt}

\noindent Then as $n-1>8$, we will get a independent set of size 5, containing $Y$; a contradiction.

\vspace{8pt}
\noindent \textbf{Case 2: $n = 9$}
\vspace{4pt}

\noindent In order to avoid an independent set of size 5, each vertex of $X$ must be adjacent to at least one vertex of $Y$. Therefore as $\vert X \vert =8$, we get that each $y_i$ ($1 \le i \le 4$) will have exactly two disjoint neighbors in $X$ and  without loss of generality then we have the following subcases. In the first subcase, we assume that without loss of generality $Y$ has a vertex which is adjacent to two vertices of $X$ which are length 2 apart in the $C_8$. In the second subcase we assume that without loss of generality $Y$ has a vertex which is adjacent to two vertices of $X$ which are length 3 apart in the $C_8$. In the third subcase  we assume that without loss of generality all vertices of $Y$ are adjacent to two vertices of $X$ which are  length 4 apart in the $C_8$.

\vspace{6pt}
\noindent \textbf{Subcase 2.1:} $Y$ has a vertex which is adjacent to two vertices of $X$ which are length 2 apart in the $C_8$.

\vspace{4pt}
\noindent Without loss of generality assume that $y_1$ is adjacent in to $u_1$ and $u_3$ and $y_2$ is adjacent in to $u_2$. In order to avoid an independent set of size 5 consisting of $\{u_1, u_3,y_2,y_3,y_4\}$, we get that $(u_1,u_3)$ is an edge. Without loss of generality, then $y_2$ is adjacent  to one vertex of $\{u_8, u_7,u_6\}$. In the first possibility, as $(u_2,u_8)$ must be an edge and therefore we will get a cycle  $(u_3, ...,u_8,y_2,u_2,u_1,u_3)$ of length $9$, contrary to the assumption. In the next possibility, as $(u_2,u_7)$ must be an edge and therefore we will get a cycle  $(u_3,u_4,u_5,u_6,u_7,y_2,u_2,u_1,y_1,u_3)$ of length $9$, contrary to the assumption. In the last possibility without loss of generality gives rise to two scenarios as illustrated in the following two figures.

\begin{center}
\input{graph1.tex} 
\end{center}
\begin{center} 
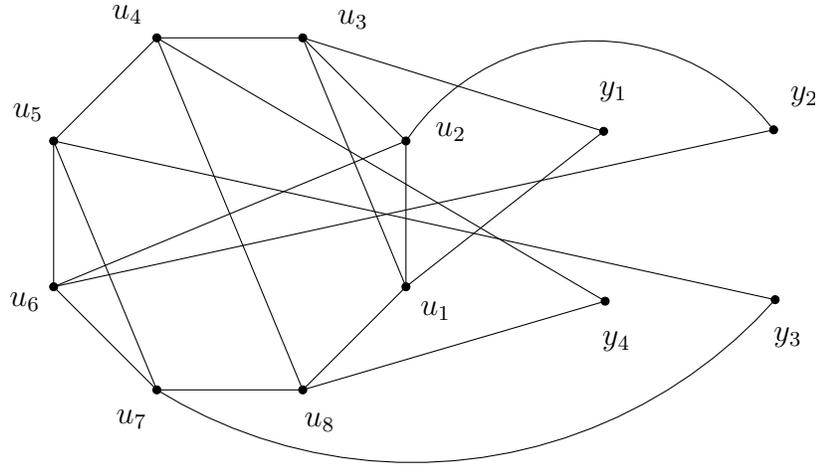
Figure 7. The first scenario
\end{center}

\vspace{12pt}
\noindent  In the first scenario, without loss of generality, $y_3$ is adjacent   to $u_5$ and $u_7$ and $y_4$ is adjacent  to $u_4$ and $u_8$. Moreover, $(u_2,u_6)$, $(u_4,u_8)$ and $(u_5,u_7)$  must be edges. Thus, we get a cycle  $(u_2,u_1,y_1,u_3,u_4,u_8,u_7,u_5,u_6,u_2)$ of length $9$, contrary to the assumption.
\vspace{4pt}

\vspace{8pt}
\noindent  In the second scenario, without loss of generality, $y_3$ is adjacent  to $u_5$ and $u_8$ and $y_4$ is adjacent to $u_4$ and $u_7$. Moreover, $(u_2,u_6)$, $(u_5,u_8)$ and $(u_4,u_7)$  must be edges. Thus, we get a  cycle $(u_1,y_1,u_3,u_4,u_7,u_8,u_5,u_6,u_2,u_1)$ of length $9$, contrary to the assumption.

\begin{center}
\input{graph2.tex} 
\end{center}
\begin{center} 
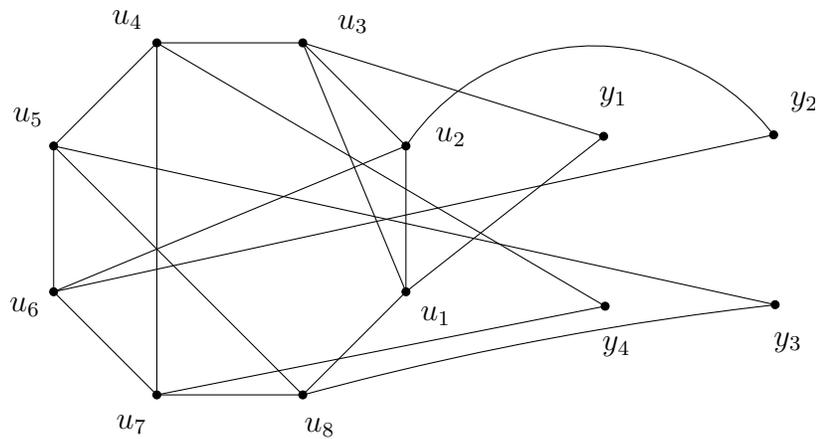
Figure 8. The second scenario
\end{center}

\vspace{6pt}
\noindent \textbf{Subcase 2.2:  } $Y$ has a vertex which is adjacent to two vertices of $X$ which are length 3 apart in the $C_8$

\vspace{6pt}
\noindent Without loss of generality assume that $y_1$ is adjacent  to $u_1$ and $u_4$ and $y_2$ is adjacent to $u_3$. In order to avoid an independent set of size 5 consisting of $\{u_1,u_4,y_2,y_3,y_4\}$, we get that $(u_1,u_4)$ is an edge. Without loss of generality, then $y_2$ must be adjacent to one vertex of 
$\{u_2,u_7,u_8,u_6\}$.

\vspace{6pt}
\noindent The first possibility, is impossible as we get a cycle  $(u_1,u_2, y_2, u_3,...,u_8,u_1)$ of length $9$, contrary to the assumption. In the second possibility we get a cycle  $(u_7,y_2, u_3,u_2,u_1,$ $y_1,u_4,u_5,u_6,u_7)$ of length $9$, contrary to the assumption. In the third possibility we get a  cycle  $(u_5,u_4, u_1,u_2,u_3,y_2,u_8,u_7,u_6,u_5)$ of length $9$, contrary to the assumption. In the fourth possibility is without loss of generality $y_3$ is adjacent to $u_2$. But then the  only options for $y_3$ to be adjacent   to is to two vertices of $\{u_5,u_8,u_7\}$. However,  if $y_3$ is adjacent to $u_2$ and $u_5$ it leads to  $(u_1,u_4,u_3,u_2,y_3,u_5,u_6,u_7,u_8,u_1)$, a  cycle of length 9 and if $y_3$ is adjacent in red to $u_2$ and $u_8$ it leads to  $(u_8,y_3,u_2,u_1,u_4,u_3,y_2,u_6,u_7,u_8)$ a cycle of length 9, contrary to the assumption.  Therefore, we are left with the scenario when $y_3$ is adjacent to $u_2$ and $u_7$. But then in order to avoid an independent set of size 5, we get that $(u_2,u_7)$ and $(u_3,u_6)$ are edges. This is illustrated in the following figure.  

\begin{center}
\input{graph3.tex} 
\end{center}
\begin{center} 
Figure 9. The only remaining scenario
\end{center}

\vspace{8pt}
\noindent  However, even in this scenario too we get a  cycle $(u_1,u_8,u_7,u_2,u_3,u_6,u_5,u_4,y_1,u_1)$ of length $9$, contrary to the assumption.

\vspace{14pt}

\noindent \textbf{Subcase 2.3:} All vertices of $Y$ are adjacent to two vertices of $X$ which are  length 4 apart in the $C_8$

\vspace{8pt}
\noindent Without loss of generality assume that $y_1$ is adjacent in to $u_1$ and $u_5$ and $y_2$ is adjacent to $u_2$ and $u_6$ and $y_3$ is adjacent  to $u_4$ and $u_8$

\vspace{8pt}

\noindent In order to avoid an independent set of size 5, we get that $(u_1,u_5)$, $(u_2,u_6)$ and $(u_4,u_8)$ are edges. But then we get cycle $(u_3,u_4,u_8,u_7,u_6,u_5,y_1,u_1,u_2,u_3)$ of length $9$, contrary to the assumption.

\vspace{12pt}

\noindent \textbf{Case 3: $n = 8$}
\vspace{8pt}

\noindent In order to avoid an independent set of size 5, each vertex of $X$ must be adjacent to at least one vertex of $Y$. Therefore as $\vert X \vert =7$, we get that without loss of generality each $y_i$ ($1 \le i \le 2$) will have exactly two disjoint  neighbors  in $X$  such that $\Gamma_R(x) \cap \Gamma_R(x) \cap X=\emptyset$  for all $x \in \{y_1,y_2\}$ and $y \in \{y_3,y_4\}$. This will give rise to the following subcases. In the first subcase, we assume that without loss of generality $\{y_1,y_2\}$ has a vertex which is adjacent to two vertices of $X$ which are length 2 apart in the $C_7$. In the second subcase we assume that without loss of generality both $\{y_1,y_2\}$ are adjacent to two vertices of $X$  which are length 3 apart in the $C_7$.  

\vspace{14pt}
\noindent \textbf{Subcase 3.1:} $\{y_1,y_2\}$ has a vertex which is adjacent to two vertices of $X$ which are length 2 apart in the $C_7$.

\vspace{4pt}
\noindent Without loss of generality assume that $y_1$ is adjacent to $u_1$ and $u_3$. Moreover, $y_2$ is adjacent to $u_2$ and $u_4$ or $y_2$ is adjacent  to $u_2$ and $u_5$ or else $y_2$ is adjacent to $u_4$ and $u_6$. In order to avoid an independent set of size 5 consisting of $\{u_1, u_3,y_2,y_3,y_4\}$, we get that $(u_1,u_3)$ is a  edge.  In the first possibility we will get a cycle $(u_1,u_3,u_2,y_2,u_4,u_5,u_6,u_7,$ $u_1)$ of length $8$, contrary to the assumption. In the second possibility, we will get a cycle $(u_1,y_1,u_3,u_2,y_2,u_5,u_6,u_7,u_1)$ of length $8$, contrary to the assumption. Therefore, we may assume the only remaining option that  $y_2$ is adjacent  to $u_4$ and $u_6$.  In order to avoid an independent set of size 5 consisting of $\{u_1, u_3,y_2,y_3,y_4\}$, we get that $(u_4,u_6)$ is a  edge. Since without loss of generality, $y_3$ is adjacent to two vertices distinct from  $\{u_1, u_3,u_4,u_6\}$, we get three possibilities given by $y_3$ is adjacent  to $u_2$ and $u_5$ or $y_3$ is adjacent  to $u_2$ and $u_7$ or else $y_3$ is adjacent  to $u_5$ and $u_7$. This will give a cycles $(u_1,u_3,u_2,y_3,u_5,u_4,u_6,u_7,u_1)$, $(u_4,u_3,u_1,u_2,y_3,u_7,u_6,u_5,u_4)$, $(u_1,u_2,u_3,u_4,u_6,u_5,y_3,u_7,u_1)$ respectively in the three cases, contrary to the assumption.

\vspace{14pt}
\noindent \textbf{Subcase 3.2:} Both $\{y_1,y_2\}$ are adjacent to two vertices of $X$ which are length 3 apart in the $C_7$.

\vspace{4pt}
\noindent Without loss of generality assume that $y_1$ is adjacent to $u_1$ and $u_4$. By symmetry, we may assume that without loss of generality $y_2$ is adjacent   to $u_2$ and $u_5$ or $u_2$ and $u_6$. In order to avoid an independent set of size 5 consisting of $\{u_1, u_4,y_2,y_3,y_4\}$, we get that $(u_1,u_4)$ is an edge.  In the first possibility, $y_2$ is adjacent to $u_2$ and $u_5$. However, in order to avoid an independent set of size 5, $(u_2,u_5)$ is an edge and thus we  will get a  cycle $(u_1,u_4,u_3,u_2,y_2,u_5,u_6,u_7,u_1)$ of length $8$, contrary to the assumption. In the second possibility, $y_2$ is adjacent  to $u_2$ and $u_6$. In order to avoid an independent set of size 5 consisting of $\{u_2, u_6,y_1,y_3,y_4\}$, we get that $(u_2,u_6)$ is a  edge.  This scenario is illustrated in the following figure 10.

\begin{center}
\input{graph4.tex} 
\end{center}
\begin{center} 
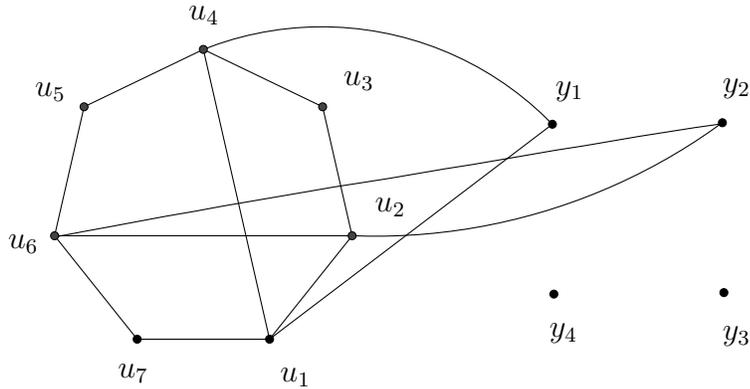
Figure 10. The only remaining scenario of Subcase 3.2
\end{center}

\vspace{8pt}
\noindent  Thus we  will get a cycle $(u_1,y_1,u_4,u_3,u_2,y_2,u_6,u_7,u_1)$ of length $8$, contrary to the assumption.

\vspace{10pt}

\noindent \textbf{Case 4: $n = 7$}
\vspace{8pt}

\noindent In order to avoid an independent set of size 5, each vertex of $X$ must be adjacent to at least one vertex of $Y$. Note that in order to avoid a $C_7$, no vertex of $Y$ can be adjacent to two vertices of $X$ which are length 1 apart in the $C_6$. Thus, we assume that without loss of generality there are two main subcases generated by when $Y$ has a vertex which is adjacent to two vertices of $X$ which are length 2 or 3 apart in the  $C_6$ and no other vertex of $Y$ is adjacent  to these two neighbors (subcase 4.1)  and when such a situation doesn't exist (subcase 4.2). 

\vspace{14pt}
\noindent \textbf{Subcase 4.1:} $Y$ has a vertex which is adjacent to two vertices of $X$ which are length 2 or 3 apart in the $C_6$ and no other vertex of $Y$ is adjacent inv to these two neighbors

\vspace{14pt}
\noindent \textbf{Subcase 4.1.1:} $Y$ has a vertex which is adjacent to two vertices of $X$ which are length 2 apart in the $C_6$ and no other vertex of $Y$ is adjacent to these two neighbors

\vspace{6pt}
\noindent As $Y$ has a vertex which is adjacent to two vertices of $X$ which are length 2 apart in the  $C_6$,  without loss of generality, assume that $y_1$ is adjacent  to $u_1$ and $u_3$ and no other vertex of $Y$ is adjacent  to  $u_1$ and $u_3$. Then, either $y_2$ is adjacent  to $u_2$ and $u_4$ or $y_2$ is adjacent  to $u_2$ and $u_5$ or else $y_2$ is adjacent  to $u_4$ and $u_6$ . In order to avoid an independent set of size 5 consisting of $\{u_1, u_3,y_2,y_3,y_4\}$, we get that $(u_1,u_3)$ is an edge.  In the first possibility, we will get a cycle $(u_1,u_3,u_2,y_2,u_4,u_5,u_6,u_1)$ of length $7$, contrary to the assumption. In the second possibility we get a cycle  $(u_1,y_1,u_3,u_2,y_2,u_5,u_6,u_1)$  of length $7$, contrary to the assumption. In the last option, we get that $y_2$ is adjacent  to $u_4$ and $u_6$.  Next note that, $y_3$  can not be adjacent to both $u_2$ and $u_5$  as then we will get  a  cycle $(u_1,y_1,u_3,u_2,y_2,u_5,u_6,u_1)$ of length $7$, contrary to the assumption. Similarly, $y_4$  can not be adjacent to both $u_2$ and $u_5$. Therefore, without loss of generality (excluding the posibilities we have already discused in case 4.1.1), we get that  $y_3$ is adjacent  to $u_5$ and $y_4$ is adjacent to $u_2$. Moreover, one can observe that in order to avoid a  $C_7$ and the cases already discussed in case 4.1.1 both $y_3$ and $y_4$ can not be adjacent to no more vertices of $X$. In order to avoid  an independent set of size 5, we get that $(u_4,u_6)$ is an edge. Also in the induced subgraph of $X$, in order to avoid a $C_7$, the only other possible  edge is either $(u_1,u_4)$ or $(u_3,u_6)$, but not both. Therefore, in the graph of $G[X \cup Y]$ illustrated in the following figure 11. Without loss of generality $(u_3,u_6)$ is marked as a broken edge to indicate that it is a possible edge (see figure 5) and its worth noting that $\{y_1,y_2,u_2,u_5\}$ forms an independent set of size 4. 

\begin{center}
\input{graph5.tex} 
\end{center}
\begin{center} 
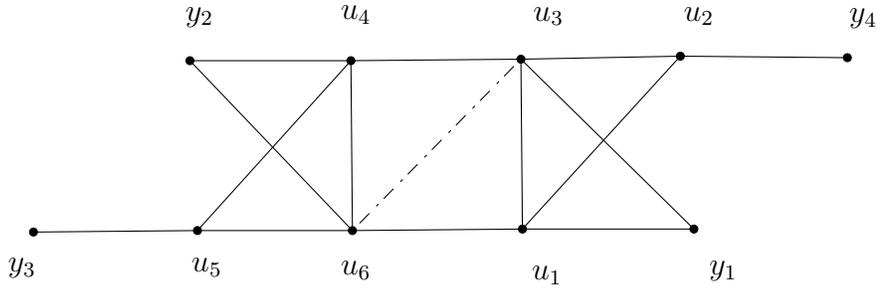
Figure 11. The structure of the graph induced by $G[X \cup Y]$ 
\end{center}

\noindent The graph $G$ is an order 24 graph with no $C_7$ and no independent set of size 5. Therefore,  by applying lemma \ref{l1}, we get that the degree of each of the vertices of $G$ is greater than or equal to 5. Therefore, $u_4$ has to be adjacent to some vertex outside $X \cup Y$ (say $w$). In order to avoid an independent set of size 5 consisting of $\{w,y_1,y_2,u_2,u_5\}$, we get that $w$ must be adjacent  to one of the vertices of $\{y_1,y_2,u_2,u_5\}$. This will give a  $C_7$ containing $w$, contrary to the assumption.

\vspace{14pt}
\noindent \textbf{Subcase 4.1.2:} $Y$ has a vertex which is adjacent to two vertices of $X$ which are length 3 apart in the  $C_6$ and no other vertex of $Y$ is adjacent   to these two neighbors

\vspace{6pt}
\noindent Without loss of generality, assume that $y_1$ is adjacent  to $u_1$ and $u_4$ and no other vertex of $Y$ is adjacent $u_1$ and $u_4$. In order to avoid an independent set of size 5, we get that $(u_1,u_4)$ is an edge. Hence, one of the vertices in $\{y_2,y_3,y_4\}$ (say $y_2$), must be adjacent to either $u_2$ and $u_5$ or  $u_2$ and $u_6$ or $u_3$ and $u_5$ or else $u_3$ and $u_6$. However, if  $y_2$ is  adjacent to $u_2$ and $u_5$ (see figure 12), we will get a cycle $(u_1,u_4,u_3,u_2,y_2,u_5,u_6,u_1)$ of length $7$, contrary to the assumption. Similarly, if  $y_2$ is  adjacent to $u_2$ and $u_6$ , we will get a cycle $(u_4,y_1,u_1,u_6,y_2,u_2,u_3,u_4)$ and if $y_2$ is  adjacent to  $u_3$ and $u_5$, we will get a cycle $(u_1,y_1,u_4,u_3,y_2,u_5,u_6,u_1)$  and if $y_2$ is  adjacent to  $u_3$ and $u_6$, we will get a cycle $(u_1,u_4,u_5,u_6,y_2,u_3,u_2,u_1)$, contrary to the assumption.

\begin{center}
\input{graph6.tex} 
\end{center}
\begin{center} 
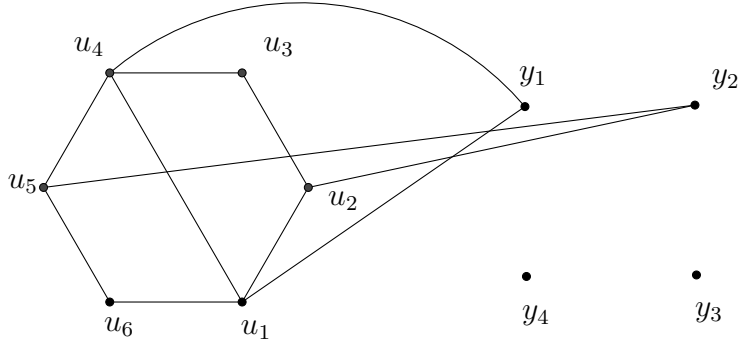
Figure 12. In subcase 4.2 when  $y_2$ is  adjacent to $u_2$ and $u_5$ 
\end{center}

\vspace{10pt}
\noindent \textbf{Subcase 4.2:} $Y$ has a vertex which is adjacent to two vertices of $X$ which are length 2 or 3 apart in the $C_6$ and some other vertex of $Y$ is adjacent in to at least one of these two neighbors

\vspace{6pt}
\noindent Without loss of generality, assume that $y_1$ is adjacent  to $u_1$ and $w$ length 2 or 3 apart. If one other vertex in $Y$ is adjacent to only one vertex $\{u_1,w\}$ or both vertices of $\{u_1,w\}$ and not adjacent to any other vertex in $X$ outside $\{u_1,w\}$ it would lead to one of the previous two subcases. Therefore, without loss of generality we may assume that, 
 $y_1$ is adjacent to $u_1$ and $u_3$ and $y_2$ is adjacent to $u_1$ and $u_4$ or else $y_1$ is adjacent to $u_1$ and $u_3$ and $y_2$ is adjacent to $u_1$ and $u_5$. In the first possibility, in order to avoid the previous subcases, without loss of generality the only option left is for $y_3$ to be adjacent to $u_2$ and $u_5$ or $y_3$ is adjacent to $u_2$ and $u_6$. If $y_3$ to be adjacent to $u_2$ and $u_5$, we get a cycle $(u_1,u_6,u_5,y_3,u_2,u_3,y_1,u_1)$ of length $7$ and if $y_3$ to be adjacent to $u_2$ and $u_6$, we get a cycle $(u_4,y_2,u_1,u_6,y_3,u_2,u_3,u_4)$ of length $7$, contrary to the assumption. In the second possibility, in order to avoid a $C_7$, we get $\{u_2,u_4,u_6\}$ will induce an independent set of size 3. But then  $\{y_1,y_2,u_2,u_4,u_6\}$ will induce an independent set of size 5 contrary to the assumption. 
\vspace{4pt}

\noindent Now continuing with the proof of the lemma for all four cases we may assume that $H$ is a $C_n$ -free graph of order $3(n-1)$ with no independent set of size 4. By \cite{JaNaRa},  we can deduce that $H$ is equal to $3K_{3n-3,1}$, $3K_{3n-3,2}$, ...,$2K_{3n-3,6}$, as $n \geq 5$. In any case $H$ contains a copy of a $3K_{n-1}$. Now consider any two vertices of $V(C_{n-1})$ say $u$ and $v$. Since there is no $C_n$ in $G$, each of the vertices $u$ and $v$ must be adjacent to at most one vertex each of the three copies of $K_4$ in $H$. Also there can be at most two  edges connecting a given $K_{n-1}$ in the copy of the $3K_{n-1}$ to the remaining  $2K_{n-1}$ subgraphs.  Therefore, as $n>5$, we can select three vertices $x_1$, $x_2$ and $x_3$ in first $K_{n-1}$, second $K_{n-1}$ and third $K_{n-1}$ respectively, such that $x_1$, $x_2$ and $x_3$ and are not adjacent to $u$ or $v$. As $\{u,v,x_1,x_2,x_3\}$ cannot be an independent set of size 5,  this will force $(u,v) \in E(G)$. Therefore, $V(C_n)$ will induce a $K_{n-1}$, since $u,v$ are arbitrary elements of $V(C_n)$. Hence the lemma.

\vspace{26pt}

\section{MAIN RESULT}

\begin{theorem}
\label{t1}

If $n \ge 4$, then
\[
r_*(C_n, K_5) =
\begin{cases} 
$ $ 13 & \text{if } n=4, \\
\hspace{20pt} & \\
$ $ 3n-1 & \text{if }   n \ge 5.  \\
\end{cases}
\]
\end{theorem}

\noindent {\bf Proof.} We break up the proof into two cases.

\vspace{6pt}
\noindent \textbf{Case $n=4$}

\vspace{6pt}
\noindent Let $R_{13}$ (see Figure 13)  represent the unique Ramsey critical $(C_4,K_5)$ graph (cf.\cite{JaRo}). 

\begin{center}
\input{grapht71.tex} 
\end{center}
\begin{center} 
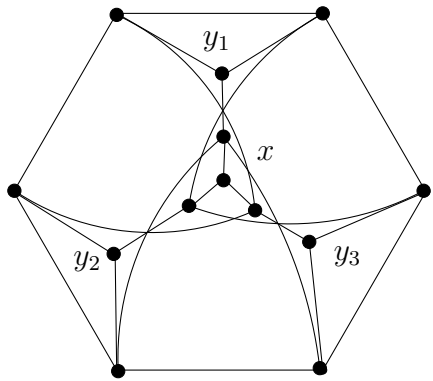
Figure 13. The unique $r(C_4,K_5)$ - critical  graph $R_{13}$ 
\end{center}

\noindent Let $R_{13}^*$ represent the graph of order 14, obtained from $R_{13}$ by adding a vertex $y$ and connecting it to the vertices $y_1,y_2,y_3$. Color the edges of $K_{13} \sqcup K_{1,12}\cong K_{14}-e$, using red and blue such that the red graph is isomorphic to $R_{13}^*$. Then the blue graph will be is isomorphic to the complement of the graph $R_{13}^*$ in $K_{14}-e$ where $e=(x,y)$.

\begin{center}
\input{grapht72.tex} 
\end{center}
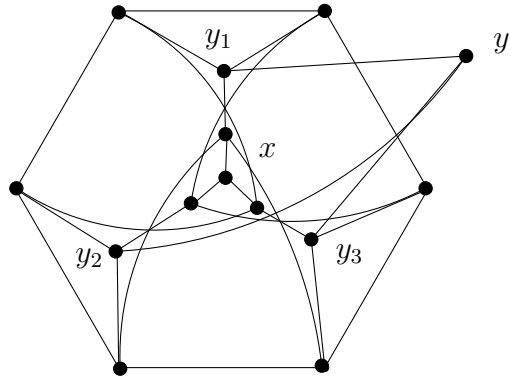
Figure 14. A graph isomorphic to the red coloring of $K_{14}-e$ which contains no  $C_4$ and which contains no $K_5$ in the complement with respect to $K_{14}-e$. Notice the missing edge between the nodes labeled $y$ and $x$ is represented by $e$.

\vspace{6pt}
\noindent \textbf{Case $n\geq 5$}

\vspace{6pt}
\noindent  Color the graph $K_{4(n-1)+1}\setminus K_{n-2}$ using red and blue colors, such that the red graph consists of a $3K_{n-1} \cup ( K_{n-1} \sqcup K_{1,1} )$ as illustrated in the following figure.  

\begin{center}
\input{graph7.tex} 
\end{center}
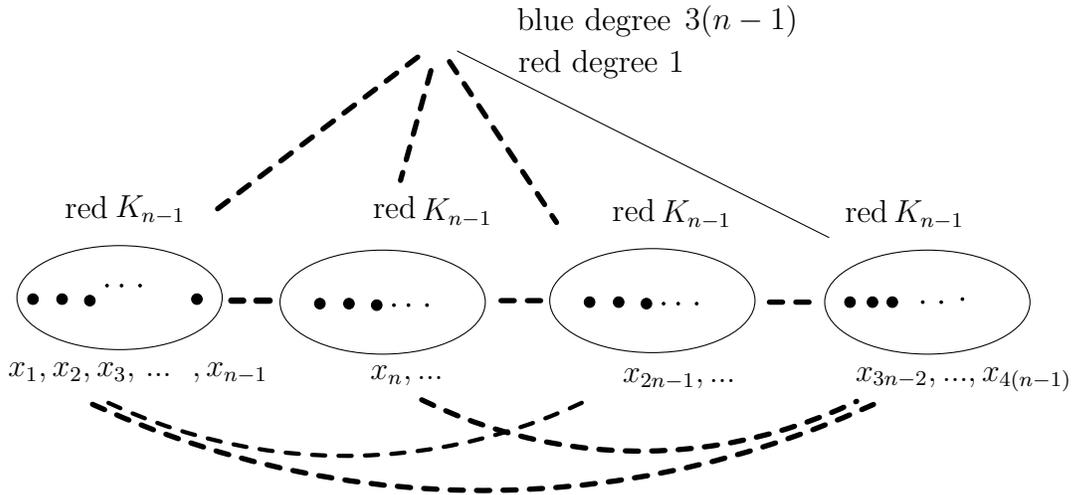
Figure 15. A coloring of $K_{4(n-1)+1}-K_{n-2}$ which  contains no red $C_n$ and no blue $K_5$

\vspace{6pt}
\noindent Therefore, $r_*(C_n, K_5)\geq 3n-1$. Next to show that, $r_*(C_n, K_5)\leq 3n-1$, consider any red/blue coloring of a graph $G=K_{4(n-1)+1}\setminus K_{n-2}$ contains no red $C_n$ and no blue $K_5$. Let $H$ be the graph obtained by deleting the vertex of degree $3n-1$ (say $v$) from $G$.

\vspace{6pt}

\noindent Then $H$ is a graph on $4(n-1)$ vertices such that it contains no red $C_n$ or a blue $K_4$. Therefore, by lemma \ref{l6}, lemma \ref{l7} and lemma \ref{l8} we get that $H$ with a red $4K_{n-1}$.  Let us denote the  sets of vertices of the four connected components by $V_1$, $V_2$, $V_3$ and $V_4$ respectively.  Since there is no red $C_n$ in the coloring, $v$ has at most one red neighbor in each of the four sets $V_1$, $V_2$, $V_3$ and $V_4$.

\vspace{16pt}

\noindent \textbf{Case 1:} $v$ is adjacent to at exactly three vertices of  some $V_i$ ($1\leq i\leq4$).  

\vspace{4pt}
\noindent Without loss of generality, we may assume that $v$ is adjacent to all the vertices of $V_1$, $V_2$ and $V_3$. That is, $v$ is adjacent to at least 5 vertices of each $V_i$ ($1\leq i\leq3$).  Select  $v_1 \in V_1$,  $v_2 \in V_2$ and $v_3 \in V_3$ such that  $v_1$ has no red neighbors in $G[V_2 \cup V_3 \cup V_4 \cup \{v \}]$,  $v_2$ has no red neighbors in $G[V_1 \cup V_3 \cup V_4 \cup \{v \}]$ and $v_3$ has no red neighbors in $G[V_1 \cup V_2 \cup V_4 \cup \{v \}]$ (this is possible because each $V_i$ can have at most 4 vertices with red neighbors outside $V_i$). Because there are no
red $C_n$ 's in the coloring, we can find a $v_4 \in V_4$ such that $\{v, v_4 \}$ is colored blue. Then, $\{v, v_1, v_2, v_3, v_4 \}$ will induce a blue $K_5$, a contradiction.

\vspace{16pt}

\noindent \textbf{Case 2:} $v$ is adjacent to  exactly  two vertices of some $V_i$ (say $V_4$). 

\vspace{4pt}

\noindent Without loss of generality, we may assume that $v$ is adjacent to all the vertices of $V_1$, $V_2$ and all but one vertex of $V_3$. That is, $v$ is adjacent to at least 5 vertices of each $V_i$ ($1\leq i\leq2$) and at least 4 vertices of $V_3$. Since there are at least two vertices of $V_4$ adjacent to $v$ in blue, select $v_4 \in V_4$  such that $v_3$ has no red neighbors in $G[ V_3 \cup \{v \}]$. Next select any vertex $v_3 \in V_3$ such that $(v_3,v)$ is colored blue. Then, $\{v, v_3, v_4 \}$ will induce a blue $K_3$. Finally select  $v_1 \in V_1$ and  $v_2 \in V_2$ such that  $v_1$ has no red neighbors in $G[V_2 \cup V_3 \cup V_4 \cup \{v \}]$ and  $v_2$ has no red neighbors in $G[V_1 \cup V_3 \cup V_4 \cup \{v \}]$ (this is possible because each $V_i$ can have at most 4 vertices with red neighbors outside $V_i$). Then, $\{v, v_1, v_2, v_3, v_4 \}$ will induce a blue $K_5$, a contradiction.

\vspace{6pt}

\noindent \textbf{Case 3:} Both case 1 and 2, do not hold.  

\vspace{4pt}
\noindent Without loss of generality, we may assume that $v$ is adjacent to at least 5 vertices of each $V_1$  and at least 4 vertices of ($2\leq i\leq4$).  Since there are at least three vertices of $V_4$ adjacent to $v$ in blue, select $v_4\in V_4$  such that $v_4$ has no red neighbors in $G[ V_2 \cup V_3 \cup \{v \}]$. Next select $v_2 \in V_2$ and $v_3 \in V_3$ such that $(v_2,v)$, $(v_3,v)$ $(v_2,v_3)$ are colored blue. Then, $\{v, v_2, v_3, v_4 \}$ will induce a blue $K_4$. Finally select  $v_1 \in V_1$  such that  $v_1$ has no red neighbors in $G[V_2 \cup V_3 \cup V_4 \cup \{v \}]$ (this is possible because each $V_i$ can have at most 4 vertices with red neighbors outside $V_i$). Then, $\{v, v_1, v_2, v_3, v_4 \}$ will induce a blue $K_5$, a contradiction.

\vspace{6pt}
\noindent Hence  $r_*(C_n, K_4)\leq 3n-1$. Therefore, $r_*(C_n, K_4)=3n-1$.

\section{ALL RAMSEY $(C_n,K_5)$ CRITICAL GRAPHS WHEN $n\geq 4$}

\noindent The thirty  graphs $R_{4n-4,1}$, $...$  ,$R_{4n-4,30}$ are illustrated in the following figure.  
   
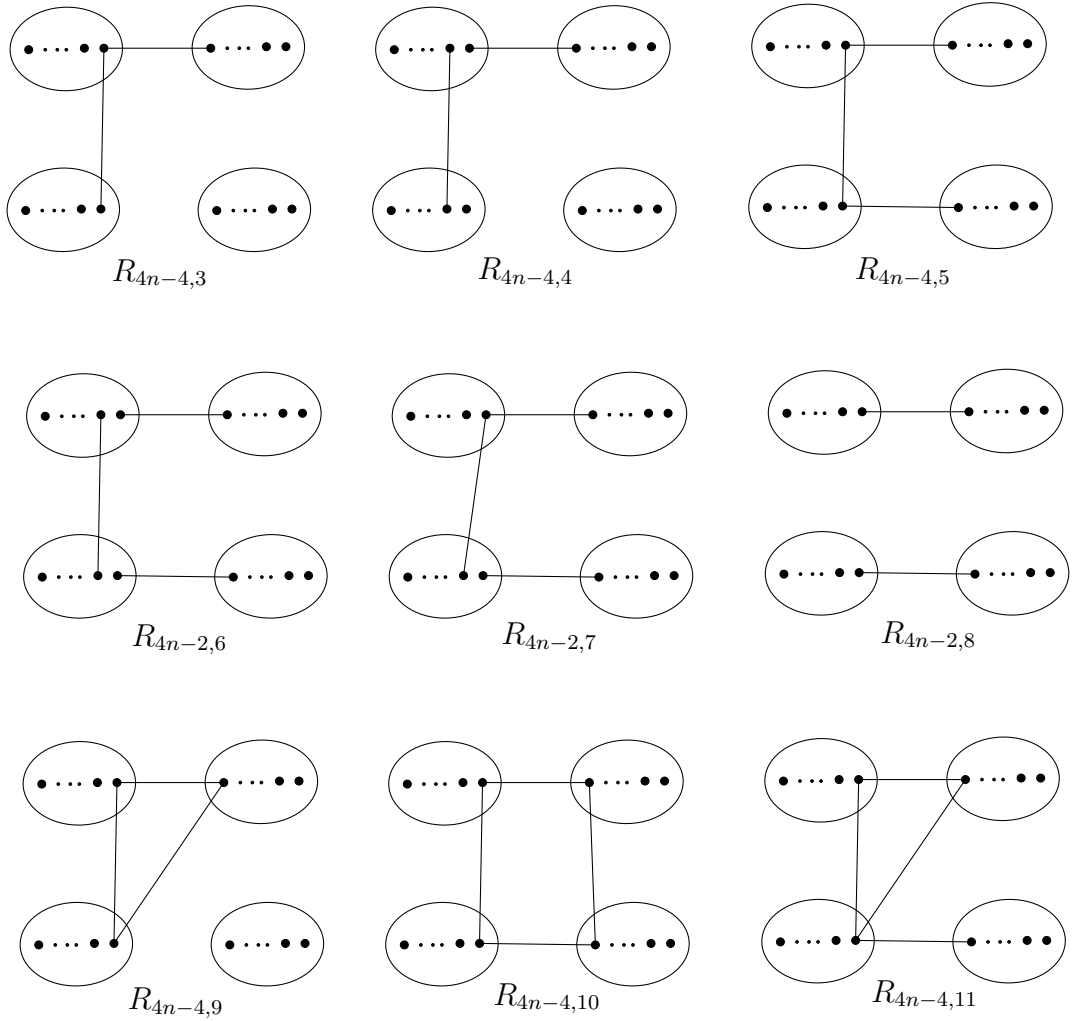
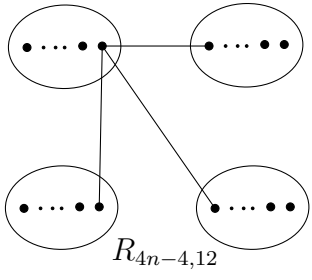
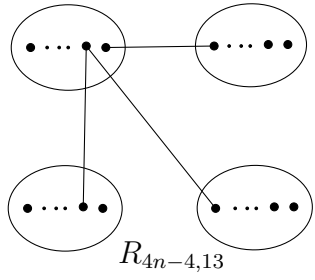
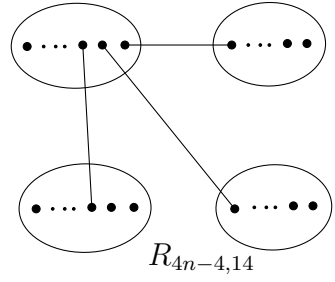
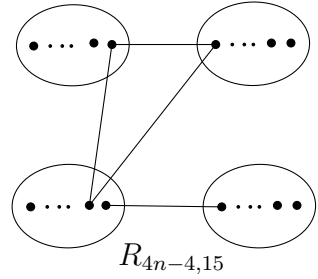
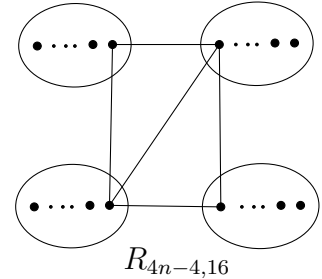
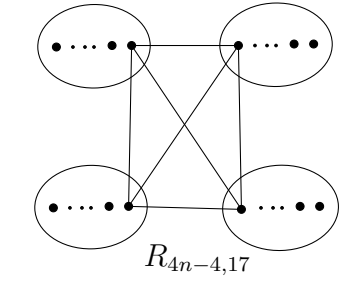
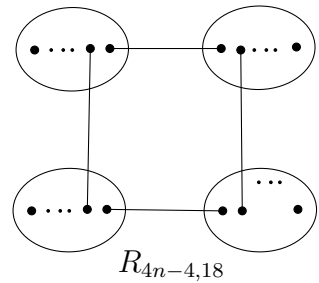
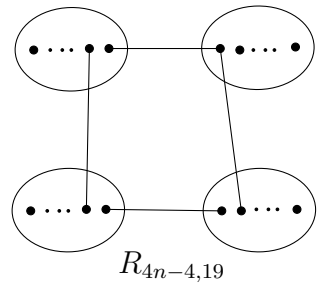
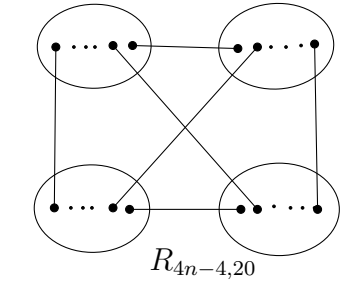
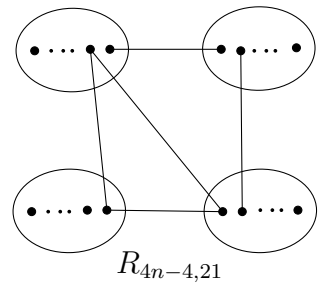
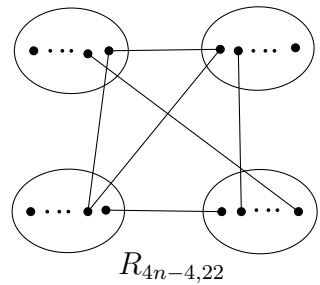
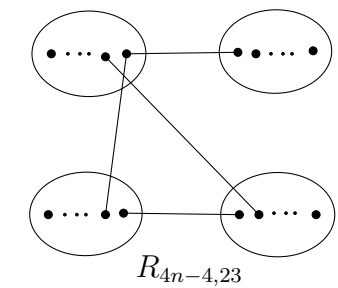
\begin{center}
\input{graph81.tex} 
\vspace{4pt}
\input{graph82.tex} 
\vspace{4pt}
\input{graph83.tex} 
\vspace{4pt}
\input{graph84.tex} 
\input{graph85.tex} 
\input{graph86.tex} 
\input{graph87.tex} 
\input{graph88.tex} 
\input{graph89.tex} 
\input{graph810.tex} 
\end{center}
\begin{center}
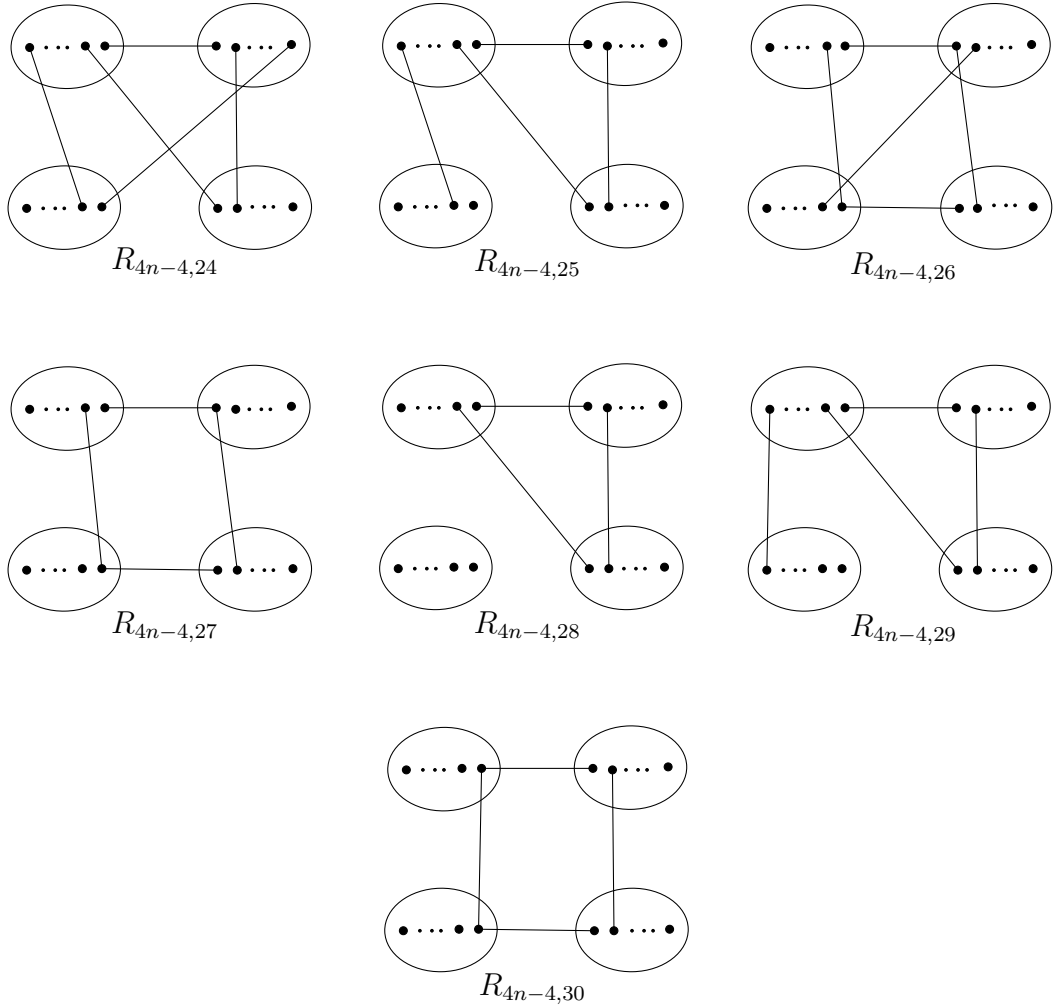
Figure 16. The red graphs $R_{4n-4,1}$, $R_{4n-4,2}$, $...$  ,$R_{4n-4,30}$ when $n \ge6$
\end{center}

\vspace{6pt}

\begin{lemma}
\label{l9}

\vspace{10pt}
When $n>3$, all $r(C_n,K_5)$ critical graphs will be generated by one of the following.

\vspace{10pt}
\noindent \textbf{(a)}  One critical generated when $n=4$, with the red graph, of the red/blue coloring, corresponding to $R_{13}$ (see Figure 13). 

\vspace{10pt}
\noindent \textbf{(b)} Thirty critical graph generated when $n=5$, with the red graph, of the red/blue coloring, corresponding to $R_{16,1}$, $R_{16,2}$, $...$ ,$R_{16,30}$ (Figure 15).

\vspace{10pt}
\noindent \textbf{(c)}  Nineteen critical graph generated when $n=6$, with the red graph, of the red/blue coloring, corresponding to $R_{20,1}$, $R_{20,2}$, $...$ ,$R_{20,19}$ (Figure 15).

\vspace{10pt}
\noindent \textbf{(d)} Eighteen critical graph generated when $n=7$, with the red graph, of the red/blue coloring, corresponding to $R_{24,1}$, $R_{24,2}$, $...$ ,$R_{24,18}$  (Figure 15).

\vspace{10pt}
\noindent \textbf{(e)}  Seventeen critical graph generated when  $n \geq 8$, with the red graph, of the red/blue coloring, corresponding to $R_{4n-4,1}$, $R_{4n-4,2}$, $...$ ,$n_{4n-4,16}$ or $R_{4n-4,17}$  (Figure 15). 

\end{lemma}

\noindent {\bf Proof.} \noindent \textbf{(a)} When $n=4$ the result follows from  (\cite{JaRo}, Lemma 11) and the fact $r(C_4,K_5)=14$ (cf. \cite{JaRo}).

\noindent \textbf{(b)-(e)} When $n\geq8$, for critical $r(C_n,K_5)$ the red graph of the red/blue coloring corresponding to $K_{4(n-1)}$,  must contain a  $4K_{n-1}$ (by Lemma \ref{l8}).  Note that in order to avoid a  $C_n$ there can be at most one red edge between any two of the  red $K_{n-1}$ graphs. We see that these graphs are generated by the 6 connected subgraphs on 4 vertices ($R_{4n-4,k}$ for $k \in \{5,10,11,12,16,17\}$), 5 disconnected subgraphs on 5 vertices ($R_{4n-4,k}$ for $k \in \{1,2,3,8,9\}$), 5 possible spliting of vertices of trees ($R_{4n-4,k}$ for $k \in \{4,6,7,13,14\}$), and one possible splitting of vertices of a $K_{1,3} +e$ ($R_{4n-4,15}$). Thus,  we see that there are only 17 distinct colorings and the corresponding red graphs  are given by $R_{4n-4,k}$ for $1\leq k\leq 17 $. Next, note that  $R_{24,18}$ is the only possible $C_7$-free graph which contains a $C_8$ and has independence number less than 5. It is generated by splitting of 4 distinct vertices of the $C_4$ ($R_{24,18}$). Thus, for $n=7$ there are only 18 distinct colorings and the corresponding red graphs  are given by $R_{24,k}$ for $1\leq k\leq 18 $. Next, note that $R_{20,19}$ is the only  possible $C_6$-free graph which contains a $C_7$ and has independence number less than 5. It  is generated by splitting of 3 distinct vertices of the $C_4$ ($R_{20,19}$). Thus, we see that for $n=6$ there are only 19 distinct colorings and the corresponding red graphs  are given by $R_{20,k}$ for $1\leq k\leq 19 $.  Finally, note that $R_{16,k}$ for $20\leq k\leq 30 $ are the only possible $C_5$-free graphs which contains a $C_6$ and has independence number less than 5 (each one of them is generated by appropriate vertex splittings of subgraphs containing a 3 or 4 cycle). Thus, we see that for $n=5$ there are 30 distinct colorings and the corresponding red graphs  are given by $R_{16,k}$ for $1\leq k\leq 30 $.

\end{document}

%% file: graphl21a.tex
\begin{tikzpicture}[line cap=round,line join=round,>=triangle 45,x=1.0cm,y=1.0cm]
\clip(-5.576705385965406,4.460708819318563) rectangle (9.019075258598198,6.969358617602924);
\draw (-5.08,5.62)-- (-4.06,5.62);
\draw (-4.6,6.6)-- (-4.06,5.62);
\draw (-4.6,6.6)-- (-5.08,5.62);
\draw (-3.62,5.62)-- (-2.6,5.62);
\draw (-3.14,6.6)-- (-2.6,5.62);
\draw (-3.14,6.6)-- (-3.62,5.62);
\draw (-2.14,5.62)-- (-1.12,5.62);
\draw (-1.66,6.6)-- (-1.12,5.62);
\draw (-1.66,6.6)-- (-2.14,5.62);
\draw (-4.06,5.62)-- (-3.62,5.62);
\draw (-2.6,5.62)-- (-2.14,5.62);
\draw (-0.18,5.62)-- (0.84,5.62);
\draw (0.3,6.6)-- (0.84,5.62);
\draw (0.3,6.6)-- (-0.18,5.62);
\draw (1.28,5.62)-- (2.3,5.62);
\draw (1.76,6.6)-- (2.3,5.62);
\draw (1.76,6.6)-- (1.28,5.62);
\draw (2.76,5.62)-- (3.78,5.62);
\draw (3.24,6.6)-- (3.78,5.62);
\draw (3.24,6.6)-- (2.76,5.62);
\draw (0.84,5.62)-- (1.28,5.62);
\draw (2.3,5.62)-- (2.76,5.62);
\draw (4.72,5.58)-- (5.74,5.58);
\draw (5.2,6.56)-- (5.74,5.58);
\draw (5.2,6.56)-- (4.72,5.58);
\draw (6.16,6.56)-- (6.7,5.58);
\draw (7.66,5.58)-- (8.68,5.58);
\draw (8.14,6.56)-- (8.68,5.58);
\draw (8.14,6.56)-- (7.66,5.58);
\draw (6.7,5.58)-- (7.66,5.58);
\draw (-3.6001934236807513,5.220905727889582) node[anchor=north west] {$P_{11,1}$};
\draw (1.2840717138880586,5.2589155733181325) node[anchor=north west] {$P_{11,2}$};
\draw (6.24435654231397,5.182895882461031) node[anchor=north west] {$P_{11,3}$};
\draw (5.74,5.58)-- (6.7,5.58);
\draw (7.14,6.56)-- (6.7,5.58);
\draw (7.14,6.56)-- (6.16,6.56);
\draw (-4.58,5.98)-- (-5.08,5.62);
\draw (-4.58,5.98)-- (-4.6,6.6);
\draw (-4.58,5.98)-- (-4.06,5.62);
\draw (-1.64,5.98)-- (-2.14,5.62);
\draw (-1.64,5.98)-- (-1.66,6.6);
\draw (-1.64,5.98)-- (-1.12,5.62);
\draw (-0.18,5.62)-- (0.32,5.96);
\draw (0.32,5.96)-- (0.3,6.6);
\draw (0.32,5.96)-- (0.84,5.62);
\draw (2.76,5.62)-- (3.28,5.98);
\draw (3.28,5.98)-- (3.78,5.62);
\draw (3.28,5.98)-- (3.24,6.6);
\draw (5.22,5.96)-- (5.2,6.56);
\draw (5.22,5.96)-- (4.72,5.58);
\draw (5.22,5.96)-- (5.74,5.58);
\draw (8.16,5.98)-- (8.14,6.56);
\draw (8.16,5.98)-- (7.66,5.58);
\draw (8.16,5.98)-- (8.68,5.58);
\draw (-3.14,6.6)-- (-2.14,5.62);
\begin{scriptsize}
\draw [fill=black] (-4.6,6.6) circle (2.5pt);
\draw [fill=black] (-4.06,5.62) circle (2.5pt);
\draw [fill=black] (-5.08,5.62) circle (2.5pt);
\draw [fill=black] (-3.14,6.6) circle (2.5pt);
\draw [fill=black] (-2.6,5.62) circle (2.5pt);
\draw [fill=black] (-3.62,5.62) circle (2.5pt);
\draw [fill=black] (-1.66,6.6) circle (2.5pt);
\draw [fill=black] (-1.12,5.62) circle (2.5pt);
\draw [fill=black] (-2.14,5.62) circle (2.5pt);
\draw [fill=black] (0.3,6.6) circle (2.5pt);
\draw [fill=black] (0.84,5.62) circle (2.5pt);
\draw [fill=black] (-0.18,5.62) circle (2.5pt);
\draw [fill=black] (1.76,6.6) circle (2.5pt);
\draw [fill=black] (2.3,5.62) circle (2.5pt);
\draw [fill=black] (1.28,5.62) circle (2.5pt);
\draw [fill=black] (3.24,6.6) circle (2.5pt);
\draw [fill=black] (3.78,5.62) circle (2.5pt);
\draw [fill=black] (2.76,5.62) circle (2.5pt);
\draw [fill=black] (5.2,6.56) circle (2.5pt);
\draw [fill=black] (5.74,5.58) circle (2.5pt);
\draw [fill=black] (4.72,5.58) circle (2.5pt);
\draw [fill=black] (6.16,6.56) circle (2.5pt);
\draw [fill=black] (6.7,5.58) circle (2.5pt);
\draw [fill=black] (8.14,6.56) circle (2.5pt);
\draw [fill=black] (8.68,5.58) circle (2.5pt);
\draw [fill=black] (7.66,5.58) circle (2.5pt);
\draw [fill=black] (7.14,6.56) circle (2.5pt);
\draw [fill=black] (-4.58,5.98) circle (2.5pt);
\draw [fill=black] (-1.64,5.98) circle (2.5pt);
\draw [fill=black] (0.32,5.96) circle (2.5pt);
\draw [fill=black] (3.28,5.98) circle (2.5pt);
\draw [fill=black] (5.22,5.96) circle (2.5pt);
\draw [fill=black] (8.16,5.98) circle (2.5pt);
\end{scriptsize}
\end{tikzpicture}

%% file: graphl21b.tex
\begin{tikzpicture}[line cap=round,line join=round,>=triangle 45,x=1.0cm,y=1.0cm]
\clip(-5.614715231393957,1.5909654894629648) rectangle (9.228129408455228,4.270659592175806);
\draw (-5.12,2.78)-- (-4.02,2.78);
\draw (-4.56,3.76)-- (-4.02,2.78);
\draw (-4.56,3.76)-- (-5.12,2.78);
\draw (-3.1,3.76)-- (-2.56,2.78);
\draw (-2.1,2.78)-- (-1.08,2.78);
\draw (-1.62,3.76)-- (-1.08,2.78);
\draw (-1.62,3.76)-- (-2.1,2.78);
\draw (-2.56,2.78)-- (-2.1,2.78);
\draw (0.0,2.8)-- (1.0009950772857246,2.8);
\draw (0.46,3.78)-- (1.0009950772857246,2.8);
\draw (0.46,3.78)-- (0.0,2.8);
\draw (1.44,3.78)-- (1.9609950772857245,2.8);
\draw (2.94,2.8)-- (3.96,2.8);
\draw (3.4,3.78)-- (3.96,2.8);
\draw (3.4,3.78)-- (2.94,2.8);
\draw (1.9609950772857245,2.8)-- (2.94,2.8);
\draw (-3.334124505680894,2.4271820888910858) node[anchor=north west] {$P_{11,4}$};
\draw (1.6071554000307424,2.4651919343196367) node[anchor=north west] {$P_{11,5}$};
\draw (2.42,3.78)-- (1.9609950772857245,2.8);
\draw (2.42,3.78)-- (1.44,3.78);
\draw (-3.1,3.76)-- (-3.54,2.78);
\draw (-3.54,2.78)-- (-2.56,2.78);
\draw (-4.56,3.12)-- (-5.12,2.78);
\draw (-4.56,3.12)-- (-4.56,3.76);
\draw (-4.56,3.12)-- (-4.02,2.78);
\draw (-2.1,2.78)-- (-1.6,3.14);
\draw (-1.6,3.14)-- (-1.62,3.76);
\draw (-1.6,3.14)-- (-1.08,2.78);
\draw (1.9609950772857245,2.8)-- (1.94,3.42);
\draw (1.94,3.42)-- (2.42,3.78);
\draw (1.94,3.42)-- (1.44,3.78);
\draw (3.42,3.18)-- (2.94,2.8);
\draw (3.42,3.18)-- (3.4,3.78);
\draw (3.42,3.18)-- (3.96,2.8);
\draw (5.081174401885775,2.8535463697474692)-- (6.101174401885774,2.8535463697474692);
\draw (5.541174401885776,3.833546369747469)-- (6.101174401885774,2.8535463697474692);
\draw (5.541174401885776,3.833546369747469)-- (5.081174401885775,2.8535463697474692);
\draw (6.521174401885776,3.833546369747469)-- (7.061174401885776,2.8535463697474692);
\draw (8.021174401885776,2.8535463697474692)-- (9.041174401885776,2.8535463697474692);
\draw (8.481174401885777,3.833546369747469)-- (9.041174401885776,2.8535463697474692);
\draw (8.481174401885777,3.833546369747469)-- (8.021174401885776,2.8535463697474692);
\draw (7.061174401885776,2.8535463697474692)-- (8.021174401885776,2.8535463697474692);
\draw (6.681469764742308,2.5032017797481876) node[anchor=north west] {$P_{11,6}$};
\draw (6.101174401885774,2.8535463697474692)-- (7.061174401885776,2.8535463697474692);
\draw (7.501174401885778,3.833546369747469)-- (7.061174401885776,2.8535463697474692);
\draw (7.501174401885778,3.833546369747469)-- (6.521174401885776,3.833546369747469);
\draw (5.081174401885775,2.8535463697474692)-- (5.581174401885775,3.2135463697474695);
\draw (5.581174401885775,3.2135463697474695)-- (6.101174401885774,2.8535463697474692);
\draw (5.581174401885775,3.2135463697474695)-- (5.541174401885776,3.833546369747469);
\draw (8.501174401885777,3.2335463697474696)-- (8.021174401885776,2.8535463697474692);
\draw (8.501174401885777,3.2335463697474696)-- (8.481174401885777,3.833546369747469);
\draw (8.501174401885777,3.2335463697474696)-- (9.041174401885776,2.8535463697474692);
\draw (-4.02,2.78)-- (-3.54,2.78);
\draw [shift={(7.061174401885774,4.27780909085363)}] plot[domain=4.119303967112796:5.305473993656586,variable=\t]({1.0*1.717592588110676*cos(\t r)+-0.0*1.717592588110676*sin(\t r)},{0.0*1.717592588110676*cos(\t r)+1.0*1.717592588110676*sin(\t r)});
\begin{scriptsize}
\draw [fill=black] (-4.56,3.76) circle (2.5pt);
\draw [fill=black] (-4.02,2.78) circle (2.5pt);
\draw [fill=black] (-5.12,2.78) circle (2.5pt);
\draw [fill=black] (-3.1,3.76) circle (2.5pt);
\draw [fill=black] (-2.56,2.78) circle (2.5pt);
\draw [fill=black] (-1.62,3.76) circle (2.5pt);
\draw [fill=black] (-1.08,2.78) circle (2.5pt);
\draw [fill=black] (-2.1,2.78) circle (2.5pt);
\draw [fill=black] (0.46,3.78) circle (2.5pt);
\draw [fill=black] (1.0009950772857246,2.8) circle (2.5pt);
\draw [fill=black] (0.0,2.8) circle (2.5pt);
\draw [fill=black] (1.44,3.78) circle (2.5pt);
\draw [fill=black] (1.9609950772857245,2.8) circle (2.5pt);
\draw [fill=black] (3.4,3.78) circle (2.5pt);
\draw [fill=black] (3.96,2.8) circle (2.5pt);
\draw [fill=black] (2.94,2.8) circle (2.5pt);
\draw [fill=black] (2.42,3.78) circle (2.5pt);
\draw [fill=black] (-3.54,2.78) circle (2.5pt);
\draw [fill=black] (-4.56,3.12) circle (2.5pt);
\draw [fill=black] (-1.6,3.14) circle (2.5pt);
\draw [fill=black] (1.94,3.42) circle (2.5pt);
\draw [fill=black] (3.42,3.18) circle (2.5pt);
\draw [fill=black] (5.541174401885776,3.833546369747469) circle (2.5pt);
\draw [fill=black] (6.101174401885774,2.8535463697474692) circle (2.5pt);
\draw [fill=black] (5.081174401885775,2.8535463697474692) circle (2.5pt);
\draw [fill=black] (6.521174401885776,3.833546369747469) circle (2.5pt);
\draw [fill=black] (7.061174401885776,2.8535463697474692) circle (2.5pt);
\draw [fill=black] (8.481174401885777,3.833546369747469) circle (2.5pt);
\draw [fill=black] (9.041174401885776,2.8535463697474692) circle (2.5pt);
\draw [fill=black] (8.021174401885776,2.8535463697474692) circle (2.5pt);
\draw [fill=black] (7.501174401885778,3.833546369747469) circle (2.5pt);
\draw [fill=black] (5.581174401885775,3.2135463697474695) circle (2.5pt);
\draw [fill=black] (8.501174401885777,3.2335463697474696) circle (2.5pt);
\end{scriptsize}
\end{tikzpicture}

%% file: graphl22.tex
\begin{tikzpicture}[line cap=round,line join=round,>=triangle 45,x=1.0cm,y=1.0cm]
\clip(-5.576705385965406,1.7810147166057222) rectangle (9.266139253883779,6.9883635403171995);
\draw (-5.08,5.62)-- (-4.06,5.62);
\draw (-4.6,6.6)-- (-4.06,5.62);
\draw (-4.6,6.6)-- (-5.08,5.62);
\draw (-3.14,6.6)-- (-3.62,5.62);
\draw (-1.66,6.6)-- (-1.12,5.62);
\draw (-1.66,6.6)-- (-2.14,5.62);
\draw (-4.06,5.62)-- (-3.62,5.62);
\draw (-0.18,5.62)-- (0.84,5.62);
\draw (0.3,6.6)-- (0.84,5.62);
\draw (0.3,6.6)-- (-0.18,5.62);
\draw (1.28,5.62)-- (2.3,5.62);
\draw (1.76,6.6)-- (2.3,5.62);
\draw (1.76,6.6)-- (1.28,5.62);
\draw (2.76,5.62)-- (3.78,5.62);
\draw (3.24,6.6)-- (3.78,5.62);
\draw (3.24,6.6)-- (2.76,5.62);
\draw (0.84,5.62)-- (1.28,5.62);
\draw (2.3,5.62)-- (2.76,5.62);
\draw (4.72,5.58)-- (5.74,5.58);
\draw (5.2,6.56)-- (5.74,5.58);
\draw (5.2,6.56)-- (4.72,5.58);
\draw (6.16,6.56)-- (6.7,5.58);
\draw (7.66,5.58)-- (8.68,5.58);
\draw (8.14,6.56)-- (8.68,5.58);
\draw (8.14,6.56)-- (7.66,5.58);
\draw (6.7,5.58)-- (7.66,5.58);
\draw (-3.6001934236807513,5.220905727889582) node[anchor=north west] {$P_{11,7}$};
\draw (1.2840717138880586,5.2589155733181325) node[anchor=north west] {$P_{11,8}$};
\draw (6.24435654231397,5.182895882461031) node[anchor=north west] {$P_{11,9}$};
\draw (5.74,5.58)-- (6.7,5.58);
\draw (7.14,6.56)-- (6.7,5.58);
\draw (7.14,6.56)-- (6.16,6.56);
\draw (-5.12,2.78)-- (-4.02,2.78);
\draw (-4.56,3.76)-- (-4.02,2.78);
\draw (-4.56,3.76)-- (-5.12,2.78);
\draw (-3.1,3.76)-- (-2.56,2.78);
\draw (-2.1,2.78)-- (-1.08,2.78);
\draw (-1.62,3.76)-- (-1.08,2.78);
\draw (-1.62,3.76)-- (-2.1,2.78);
\draw (-2.56,2.78)-- (-2.1,2.78);
\draw (0.0,2.8)-- (1.0009950772857246,2.8);
\draw (0.46,3.78)-- (1.0009950772857246,2.8);
\draw (0.46,3.78)-- (0.0,2.8);
\draw (1.44,3.78)-- (1.9609950772857245,2.8);
\draw (2.94,2.8)-- (3.96,2.8);
\draw (3.4,3.78)-- (3.96,2.8);
\draw (3.4,3.78)-- (2.94,2.8);
\draw (1.9609950772857245,2.8)-- (2.94,2.8);
\draw (-3.334124505680894,2.427182088891088) node[anchor=north west] {$P_{11,10}$};
\draw (1.6071554000307424,2.465191934319639) node[anchor=north west] {$P_{11,11}$};
\draw (2.42,3.78)-- (1.9609950772857245,2.8);
\draw (2.42,3.78)-- (1.44,3.78);
\draw (-3.1,3.76)-- (-3.54,2.78);
\draw (-3.54,2.78)-- (-2.56,2.78);
\draw (-1.64,5.98)-- (-2.14,5.62);
\draw (-1.64,5.98)-- (-1.66,6.6);
\draw (-1.64,5.98)-- (-1.12,5.62);
\draw (-0.18,5.62)-- (0.32,5.96);
\draw (0.32,5.96)-- (0.3,6.6);
\draw (0.32,5.96)-- (0.84,5.62);
\draw (2.76,5.62)-- (3.28,5.98);
\draw (3.28,5.98)-- (3.78,5.62);
\draw (3.28,5.98)-- (3.24,6.6);
\draw (5.22,5.96)-- (5.2,6.56);
\draw (5.22,5.96)-- (4.72,5.58);
\draw (5.22,5.96)-- (5.74,5.58);
\draw (8.16,5.98)-- (8.14,6.56);
\draw (8.16,5.98)-- (7.66,5.58);
\draw (8.16,5.98)-- (8.68,5.58);
\draw (-2.1,2.78)-- (-1.6,3.14);
\draw (-1.6,3.14)-- (-1.62,3.76);
\draw (-1.6,3.14)-- (-1.08,2.78);
\draw (1.9609950772857245,2.8)-- (1.94,3.42);
\draw (1.94,3.42)-- (2.42,3.78);
\draw (1.94,3.42)-- (1.44,3.78);
\draw (3.42,3.18)-- (2.94,2.8);
\draw (3.42,3.18)-- (3.4,3.78);
\draw (3.42,3.18)-- (3.96,2.8);
\draw (6.681469764742308,2.50320177974819) node[anchor=north west] {$P_{11,12}$};
\draw (-4.02,2.78)-- (-3.54,2.78);
\draw (-1.12,5.62)-- (-2.14,5.62);
\draw (-3.1,5.98)-- (-3.14,6.6);
\draw (-3.1,5.98)-- (-3.62,5.62);
\draw (-4.6,6.6)-- (-3.62,5.62);
\draw [shift={(1.5699999999999998,6.123200700385838)}] plot[domain=3.745112119360146:5.679665841409234,variable=\t]({1.0*0.8866289781350469*cos(\t r)+-0.0*0.8866289781350469*sin(\t r)},{0.0*0.8866289781350469*cos(\t r)+1.0*0.8866289781350469*sin(\t r)});
\draw [shift={(2.02,6.065221745285946)}] plot[domain=3.6832252455070806:5.741552715262299,variable=\t]({1.0*0.863610098641432*cos(\t r)+-0.0*0.863610098641432*sin(\t r)},{0.0*0.863610098641432*cos(\t r)+1.0*0.863610098641432*sin(\t r)});
\draw (7.14,6.56)-- (7.66,5.58);
\draw (-3.1,3.76)-- (-2.56,2.78);
\draw (-3.06,3.14)-- (-2.56,2.78);
\draw (-3.06,3.14)-- (-3.1,3.76);
\draw (-3.06,3.14)-- (-3.54,2.78);
\draw (-2.5929325198241484,5.601004182175089)-- (-3.62,5.62);
\draw (-2.5929325198241484,5.601004182175089)-- (-3.14,6.6);
\draw (-2.5929325198241484,5.601004182175089)-- (-3.1,5.98);
\draw (5.156199015457151,2.834541447033195)-- (6.176199015457152,2.834541447033195);
\draw (5.636199015457152,3.8145414470331946)-- (6.176199015457152,2.834541447033195);
\draw (5.636199015457152,3.8145414470331946)-- (5.156199015457151,2.834541447033195);
\draw (7.0961990154571515,3.8145414470331946)-- (6.616199015457152,2.834541447033195);
\draw (8.57619901545715,3.8145414470331946)-- (9.116199015457148,2.834541447033195);
\draw (8.57619901545715,3.8145414470331946)-- (8.096199015457152,2.834541447033195);
\draw (6.176199015457152,2.834541447033195)-- (6.616199015457152,2.834541447033195);
\draw (8.596199015457149,3.1945414470331954)-- (8.096199015457152,2.834541447033195);
\draw (8.596199015457149,3.1945414470331954)-- (8.57619901545715,3.8145414470331946);
\draw (8.596199015457149,3.1945414470331954)-- (9.116199015457148,2.834541447033195);
\draw (9.116199015457148,2.834541447033195)-- (8.096199015457152,2.834541447033195);
\draw (7.136199015457152,3.1945414470331954)-- (7.0961990154571515,3.8145414470331946);
\draw (7.136199015457152,3.1945414470331954)-- (6.616199015457152,2.834541447033195);
\draw (5.636199015457152,3.8145414470331946)-- (6.616199015457152,2.834541447033195);
\draw (7.643266495633003,2.815545629208284)-- (6.616199015457152,2.834541447033195);
\draw (7.643266495633003,2.815545629208284)-- (7.0961990154571515,3.8145414470331946);
\draw (7.643266495633003,2.815545629208284)-- (7.136199015457152,3.1945414470331954);
\draw (7.643266495633003,2.815545629208284)-- (8.096199015457152,2.834541447033195);
\begin{scriptsize}
\draw [fill=black] (-4.6,6.6) circle (2.5pt);
\draw [fill=black] (-4.06,5.62) circle (2.5pt);
\draw [fill=black] (-5.08,5.62) circle (2.5pt);
\draw [fill=black] (-3.14,6.6) circle (2.5pt);
\draw [fill=black] (-3.62,5.62) circle (2.5pt);
\draw [fill=black] (-1.66,6.6) circle (2.5pt);
\draw [fill=black] (-1.12,5.62) circle (2.5pt);
\draw [fill=black] (-2.14,5.62) circle (2.5pt);
\draw [fill=black] (0.3,6.6) circle (2.5pt);
\draw [fill=black] (0.84,5.62) circle (2.5pt);
\draw [fill=black] (-0.18,5.62) circle (2.5pt);
\draw [fill=black] (1.76,6.6) circle (2.5pt);
\draw [fill=black] (2.3,5.62) circle (2.5pt);
\draw [fill=black] (1.28,5.62) circle (2.5pt);
\draw [fill=black] (3.24,6.6) circle (2.5pt);
\draw [fill=black] (3.78,5.62) circle (2.5pt);
\draw [fill=black] (2.76,5.62) circle (2.5pt);
\draw [fill=black] (5.2,6.56) circle (2.5pt);
\draw [fill=black] (5.74,5.58) circle (2.5pt);
\draw [fill=black] (4.72,5.58) circle (2.5pt);
\draw [fill=black] (6.16,6.56) circle (2.5pt);
\draw [fill=black] (6.7,5.58) circle (2.5pt);
\draw [fill=black] (8.14,6.56) circle (2.5pt);
\draw [fill=black] (8.68,5.58) circle (2.5pt);
\draw [fill=black] (7.66,5.58) circle (2.5pt);
\draw [fill=black] (7.14,6.56) circle (2.5pt);
\draw [fill=black] (-4.56,3.76) circle (2.5pt);
\draw [fill=black] (-4.02,2.78) circle (2.5pt);
\draw [fill=black] (-5.12,2.78) circle (2.5pt);
\draw [fill=black] (-3.1,3.76) circle (2.5pt);
\draw [fill=black] (-2.56,2.78) circle (2.5pt);
\draw [fill=black] (-1.62,3.76) circle (2.5pt);
\draw [fill=black] (-1.08,2.78) circle (2.5pt);
\draw [fill=black] (-2.1,2.78) circle (2.5pt);
\draw [fill=black] (0.46,3.78) circle (2.5pt);
\draw [fill=black] (1.0009950772857246,2.8) circle (2.5pt);
\draw [fill=black] (0.0,2.8) circle (2.5pt);
\draw [fill=black] (1.44,3.78) circle (2.5pt);
\draw [fill=black] (1.9609950772857245,2.8) circle (2.5pt);
\draw [fill=black] (3.4,3.78) circle (2.5pt);
\draw [fill=black] (3.96,2.8) circle (2.5pt);
\draw [fill=black] (2.94,2.8) circle (2.5pt);
\draw [fill=black] (2.42,3.78) circle (2.5pt);
\draw [fill=black] (-3.54,2.78) circle (2.5pt);
\draw [fill=black] (-1.64,5.98) circle (2.5pt);
\draw [fill=black] (0.32,5.96) circle (2.5pt);
\draw [fill=black] (3.28,5.98) circle (2.5pt);
\draw [fill=black] (5.22,5.96) circle (2.5pt);
\draw [fill=black] (8.16,5.98) circle (2.5pt);
\draw [fill=black] (-1.6,3.14) circle (2.5pt);
\draw [fill=black] (1.94,3.42) circle (2.5pt);
\draw [fill=black] (3.42,3.18) circle (2.5pt);
\draw [fill=black] (-3.1,5.98) circle (2.5pt);
\draw [fill=black] (-3.06,3.14) circle (2.5pt);
\draw [fill=black] (-2.5929325198241484,5.601004182175089) circle (2.5pt);
\draw [fill=black] (5.636199015457152,3.8145414470331946) circle (2.5pt);
\draw [fill=black] (6.176199015457152,2.834541447033195) circle (2.5pt);
\draw [fill=black] (5.156199015457151,2.834541447033195) circle (2.5pt);
\draw [fill=black] (7.0961990154571515,3.8145414470331946) circle (2.5pt);
\draw [fill=black] (6.616199015457152,2.834541447033195) circle (2.5pt);
\draw [fill=black] (8.57619901545715,3.8145414470331946) circle (2.5pt);
\draw [fill=black] (9.116199015457148,2.834541447033195) circle (2.5pt);
\draw [fill=black] (8.096199015457152,2.834541447033195) circle (2.5pt);
\draw [fill=black] (8.596199015457149,3.1945414470331954) circle (2.5pt);
\draw [fill=black] (7.136199015457152,3.1945414470331954) circle (2.5pt);
\draw [fill=black] (7.643266495633003,2.815545629208284) circle (2.5pt);
\end{scriptsize}
\end{tikzpicture}

%% file: graphl23.tex
\begin{tikzpicture}[line cap=round,line join=round,>=triangle 45,x=1.0cm,y=1.0cm]
\clip(-5.5196906178225795,1.5529556440344143) rectangle (9.323154022026605,7.349457071888431);
\draw (-5.08,5.62)-- (-4.06,5.62);
\draw (-4.6,6.6)-- (-4.06,5.62);
\draw (-4.6,6.6)-- (-5.08,5.62);
\draw (-3.14,6.6)-- (-3.62,5.62);
\draw (-1.66,6.6)-- (-1.12,5.62);
\draw (-1.66,6.6)-- (-2.14,5.62);
\draw (-4.06,5.62)-- (-3.62,5.62);
\draw (-0.18,5.62)-- (0.84,5.62);
\draw (0.3,6.6)-- (0.84,5.62);
\draw (0.3,6.6)-- (-0.18,5.62);
\draw (1.28,5.600995077285725)-- (2.3,5.62);
\draw (1.76,6.6)-- (2.3,5.62);
\draw (1.76,6.6)-- (1.28,5.600995077285725);
\draw (2.76,5.600995077285725)-- (3.78,5.62);
\draw (3.24,6.6)-- (3.78,5.62);
\draw (3.24,6.6)-- (2.76,5.600995077285725);
\draw (0.84,5.62)-- (1.28,5.600995077285725);
\draw (2.3,5.62)-- (2.76,5.600995077285725);
\draw (4.72,5.58)-- (5.74,5.58);
\draw (5.2,6.56)-- (5.74,5.58);
\draw (5.2,6.56)-- (4.72,5.58);
\draw (7.66,5.58)-- (8.68,5.58);
\draw (8.14,6.56)-- (8.68,5.58);
\draw (8.14,6.56)-- (7.66,5.58);
\draw (6.7,5.58)-- (7.66,5.58);
\draw (-3.6001934236807513,5.220905727889579) node[anchor=north west] {$P_{11,13}$};
\draw (1.2840717138880586,5.25891557331813) node[anchor=north west] {$P_{11,14}$};
\draw (6.24435654231397,5.182895882461028) node[anchor=north west] {$P_{11,15}$};
\draw (5.74,5.58)-- (6.7,5.58);
\draw (7.14,6.56)-- (6.7,5.58);
\draw (-5.12,2.78)-- (-4.02,2.78);
\draw (-4.56,3.76)-- (-4.02,2.78);
\draw (-4.56,3.76)-- (-5.12,2.78);
\draw (-2.1,2.78)-- (-1.08,2.78);
\draw (-1.62,3.76)-- (-1.08,2.78);
\draw (-1.62,3.76)-- (-2.1,2.78);
\draw (4.9,2.8)-- (5.90099507728572,2.7809950772857244);
\draw (5.36,3.78)-- (5.90099507728572,2.7809950772857244);
\draw (5.36,3.78)-- (4.9,2.8);
\draw (6.34,3.78)-- (6.8609950772857164,2.7809950772857244);
\draw (7.84,2.7809950772857244)-- (8.86,2.8);
\draw (8.3,3.78)-- (8.86,2.8);
\draw (8.3,3.78)-- (7.84,2.7809950772857244);
\draw (6.8609950772857164,2.7809950772857244)-- (7.84,2.7809950772857244);
\draw (-3.334124505680894,2.4271820888910858) node[anchor=north west] {$P_{11,16}$};
\draw (1.6071554000307424,2.4461870116053612) node[anchor=north west] {$P_{11,17}$};
\draw (7.32,3.78)-- (6.8609950772857164,2.7809950772857244);
\draw (7.32,3.78)-- (6.34,3.78);
\draw (-2.6,3.76)-- (-3.04,2.78);
\draw (-1.64,5.98)-- (-2.14,5.62);
\draw (-1.64,5.98)-- (-1.66,6.6);
\draw (-1.64,5.98)-- (-1.12,5.62);
\draw (-0.18,5.62)-- (0.32,5.96);
\draw (0.32,5.96)-- (0.3,6.6);
\draw (0.32,5.96)-- (0.84,5.62);
\draw (2.76,5.600995077285725)-- (3.28,5.98);
\draw (3.28,5.98)-- (3.78,5.62);
\draw (3.28,5.98)-- (3.24,6.6);
\draw (5.22,5.96)-- (5.2,6.56);
\draw (5.22,5.96)-- (4.72,5.58);
\draw (5.22,5.96)-- (5.74,5.58);
\draw (8.16,5.98)-- (8.14,6.56);
\draw (8.16,5.98)-- (7.66,5.58);
\draw (8.16,5.98)-- (8.68,5.58);
\draw (-2.1,2.78)-- (-1.6,3.14);
\draw (-1.6,3.14)-- (-1.62,3.76);
\draw (-1.6,3.14)-- (-1.08,2.78);
\draw (6.8609950772857164,2.7809950772857244)-- (6.84,3.42);
\draw (6.84,3.42)-- (7.32,3.78);
\draw (6.84,3.42)-- (6.34,3.78);
\draw (8.32,3.18)-- (7.84,2.7809950772857244);
\draw (8.32,3.18)-- (8.3,3.78);
\draw (8.32,3.18)-- (8.86,2.8);
\draw (6.681469764742308,2.5032017797481876) node[anchor=north west] {$P_{11,18}$};
\draw (-1.12,5.62)-- (-2.14,5.62);
\draw (-3.1,5.98)-- (-3.14,6.6);
\draw (-3.1,5.98)-- (-3.62,5.62);
\draw (-4.6,6.6)-- (-3.62,5.62);
\draw (7.14,6.56)-- (7.66,5.58);
\draw (-2.56,3.14)-- (-2.6,3.76);
\draw (-2.56,3.14)-- (-3.04,2.78);
\draw (-2.5929325198241484,5.601004182175089)-- (-3.62,5.62);
\draw (-2.5929325198241484,5.601004182175089)-- (-3.14,6.6);
\draw (-2.5929325198241484,5.601004182175089)-- (-3.1,5.98);
\draw [shift={(-3.1000000000000005,10.09125364266483)}] plot[domain=4.295514584349766:5.129263376419614,variable=\t]({1.0*4.890041833874585*cos(\t r)+-0.0*4.890041833874585*sin(\t r)},{0.0*4.890041833874585*cos(\t r)+1.0*4.890041833874585*sin(\t r)});
\draw [shift={(1.8,9.8818710142161)}] plot[domain=4.277472841011921:5.1473051197574575,variable=\t]({1.0*4.69935575816679*cos(\t r)+-0.0*4.69935575816679*sin(\t r)},{0.0*4.69935575816679*cos(\t r)+1.0*4.69935575816679*sin(\t r)});
\draw (7.175597755313471,5.924087868317772)-- (7.66,5.58);
\draw (7.175597755313471,5.924087868317772)-- (6.7,5.58);
\draw (-2.6,3.76)-- (-2.1,2.78);
\draw (-2.1,2.78)-- (-3.04,2.78);
\draw (-2.1,2.78)-- (-2.56,3.14);
\draw [shift={(-3.06,3.963690441925441)}] plot[domain=4.030962832452406:5.393815128316974,variable=\t]({1.0*1.5240482480242041*cos(\t r)+-0.0*1.5240482480242041*sin(\t r)},{0.0*1.5240482480242041*cos(\t r)+1.0*1.5240482480242041*sin(\t r)});
\draw (6.8609950772857164,2.7809950772857244)-- (5.36,3.78);
\draw (0.09690926431696645,2.8345414470331947)-- (1.1969092643169668,2.8345414470331947);
\draw (0.6569092643169663,3.8145414470331946)-- (1.1969092643169668,2.8345414470331947);
\draw (0.6569092643169663,3.8145414470331946)-- (0.09690926431696645,2.8345414470331947);
\draw (3.116909264316967,2.8345414470331947)-- (4.1369092643169685,2.8345414470331947);
\draw (3.5969092643169676,3.8145414470331946)-- (4.1369092643169685,2.8345414470331947);
\draw (3.5969092643169676,3.8145414470331946)-- (3.116909264316967,2.8345414470331947);
\draw (2.616909264316968,3.8145414470331946)-- (2.1769092643169676,2.8345414470331947);
\draw (3.116909264316967,2.8345414470331947)-- (3.616909264316967,3.1945414470331945);
\draw (3.616909264316967,3.1945414470331945)-- (3.5969092643169676,3.8145414470331946);
\draw (3.616909264316967,3.1945414470331945)-- (4.1369092643169685,2.8345414470331947);
\draw (2.656909264316968,3.1945414470331945)-- (2.616909264316968,3.8145414470331946);
\draw (2.656909264316968,3.1945414470331945)-- (2.1769092643169676,2.8345414470331947);
\draw (2.616909264316968,3.8145414470331946)-- (3.116909264316967,2.8345414470331947);
\draw (3.116909264316967,2.8345414470331947)-- (2.1769092643169676,2.8345414470331947);
\draw (3.116909264316967,2.8345414470331947)-- (2.656909264316968,3.1945414470331945);
\draw (6.8609950772857164,2.7809950772857244)-- (5.90099507728572,2.7809950772857244);
\begin{scriptsize}
\draw [fill=black] (-4.6,6.6) circle (2.5pt);
\draw [fill=black] (-4.06,5.62) circle (2.5pt);
\draw [fill=black] (-5.08,5.62) circle (2.5pt);
\draw [fill=black] (-3.14,6.6) circle (2.5pt);
\draw [fill=black] (-3.62,5.62) circle (2.5pt);
\draw [fill=black] (-1.66,6.6) circle (2.5pt);
\draw [fill=black] (-1.12,5.62) circle (2.5pt);
\draw [fill=black] (-2.14,5.62) circle (2.5pt);
\draw [fill=black] (0.3,6.6) circle (2.5pt);
\draw [fill=black] (0.84,5.62) circle (2.5pt);
\draw [fill=black] (-0.18,5.62) circle (2.5pt);
\draw [fill=black] (1.76,6.6) circle (2.5pt);
\draw [fill=black] (2.3,5.62) circle (2.5pt);
\draw [fill=black] (1.28,5.600995077285725) circle (2.5pt);
\draw [fill=black] (3.24,6.6) circle (2.5pt);
\draw [fill=black] (3.78,5.62) circle (2.5pt);
\draw [fill=black] (2.76,5.600995077285725) circle (2.5pt);
\draw [fill=black] (5.2,6.56) circle (2.5pt);
\draw [fill=black] (5.74,5.58) circle (2.5pt);
\draw [fill=black] (4.72,5.58) circle (2.5pt);
\draw [fill=black] (6.7,5.58) circle (2.5pt);
\draw [fill=black] (8.14,6.56) circle (2.5pt);
\draw [fill=black] (8.68,5.58) circle (2.5pt);
\draw [fill=black] (7.66,5.58) circle (2.5pt);
\draw [fill=black] (7.14,6.56) circle (2.5pt);
\draw [fill=black] (-4.56,3.76) circle (2.5pt);
\draw [fill=black] (-4.02,2.78) circle (2.5pt);
\draw [fill=black] (-5.12,2.78) circle (2.5pt);
\draw [fill=black] (-2.6,3.76) circle (2.5pt);
\draw [fill=black] (-1.62,3.76) circle (2.5pt);
\draw [fill=black] (-1.08,2.78) circle (2.5pt);
\draw [fill=black] (-2.1,2.78) circle (2.5pt);
\draw [fill=black] (5.36,3.78) circle (2.5pt);
\draw [fill=black] (5.90099507728572,2.7809950772857244) circle (2.5pt);
\draw [fill=black] (4.9,2.8) circle (2.5pt);
\draw [fill=black] (6.34,3.78) circle (2.5pt);
\draw [fill=black] (6.8609950772857164,2.7809950772857244) circle (2.5pt);
\draw [fill=black] (8.3,3.78) circle (2.5pt);
\draw [fill=black] (8.86,2.8) circle (2.5pt);
\draw [fill=black] (7.84,2.7809950772857244) circle (2.5pt);
\draw [fill=black] (7.32,3.78) circle (2.5pt);
\draw [fill=black] (-3.04,2.78) circle (2.5pt);
\draw [fill=black] (-1.64,5.98) circle (2.5pt);
\draw [fill=black] (0.32,5.96) circle (2.5pt);
\draw [fill=black] (3.28,5.98) circle (2.5pt);
\draw [fill=black] (5.22,5.96) circle (2.5pt);
\draw [fill=black] (8.16,5.98) circle (2.5pt);
\draw [fill=black] (-1.6,3.14) circle (2.5pt);
\draw [fill=black] (6.84,3.42) circle (2.5pt);
\draw [fill=black] (8.32,3.18) circle (2.5pt);
\draw [fill=black] (-3.1,5.98) circle (2.5pt);
\draw [fill=black] (-2.56,3.14) circle (2.5pt);
\draw [fill=black] (-2.5929325198241484,5.601004182175089) circle (2.5pt);
\draw [fill=black] (7.175597755313471,5.924087868317772) circle (2.5pt);
\draw [fill=black] (0.6569092643169663,3.8145414470331946) circle (2.5pt);
\draw [fill=black] (1.1969092643169668,2.8345414470331947) circle (2.5pt);
\draw [fill=black] (0.09690926431696645,2.8345414470331947) circle (2.5pt);
\draw [fill=black] (2.616909264316968,3.8145414470331946) circle (2.5pt);
\draw [fill=black] (3.5969092643169676,3.8145414470331946) circle (2.5pt);
\draw [fill=black] (4.1369092643169685,2.8345414470331947) circle (2.5pt);
\draw [fill=black] (3.116909264316967,2.8345414470331947) circle (2.5pt);
\draw [fill=black] (2.1769092643169676,2.8345414470331947) circle (2.5pt);
\draw [fill=black] (3.616909264316967,3.1945414470331945) circle (2.5pt);
\draw [fill=black] (2.656909264316968,3.1945414470331945) circle (2.5pt);
\end{scriptsize}
\end{tikzpicture}

%% file: graphl2a.tex
\begin{tikzpicture}[line cap=round,line join=round,>=triangle 45,x=1.0cm,y=1.0cm]
\clip(-5.560000000000005,4.139999999999996) rectangle (9.460000000000012,7.319999999999992);
\draw (-5.08,5.62)-- (-4.06,5.62);
\draw (-4.6,6.6)-- (-4.06,5.62);
\draw (-4.6,6.6)-- (-5.08,5.62);
\draw (-3.62,5.62)-- (-2.6,5.62);
\draw (-3.14,6.6)-- (-2.6,5.62);
\draw (-3.14,6.6)-- (-3.62,5.62);
\draw (-2.14,5.62)-- (-1.12,5.62);
\draw (-1.66,6.6)-- (-1.12,5.62);
\draw (-1.66,6.6)-- (-2.14,5.62);
\draw (-4.06,5.62)-- (-3.62,5.62);
\draw (-2.6,5.62)-- (-2.14,5.62);
\draw (-0.18,5.62)-- (0.84,5.62);
\draw (0.3,6.6)-- (0.84,5.62);
\draw (0.3,6.6)-- (-0.18,5.62);
\draw (1.28,5.62)-- (2.3,5.62);
\draw (1.76,6.6)-- (2.3,5.62);
\draw (1.76,6.6)-- (1.28,5.62);
\draw (2.76,5.62)-- (3.78,5.62);
\draw (3.24,6.6)-- (3.78,5.62);
\draw (3.24,6.6)-- (2.76,5.62);
\draw (0.84,5.62)-- (1.28,5.62);
\draw (2.3,5.62)-- (2.76,5.62);
\draw (4.72,5.58)-- (5.74,5.58);
\draw (5.2,6.56)-- (5.74,5.58);
\draw (5.2,6.56)-- (4.72,5.58);
\draw (6.16,6.56)-- (6.7,5.58);
\draw (7.66,5.58)-- (8.68,5.58);
\draw (8.14,6.56)-- (8.68,5.58);
\draw (8.14,6.56)-- (7.66,5.58);
\draw (6.7,5.58)-- (7.66,5.58);
\draw (-3.640000000000003,5.2599999999999945) node[anchor=north west] {$R_{12,1}$};
\draw (1.2800000000000025,5.279999999999995) node[anchor=north west] {$R_{12,2}$};
\draw (6.080000000000007,5.199999999999995) node[anchor=north west] {$R_{12,3}$};
\draw [shift={(1.7700000000000002,5.031379310344827)}] plot[domain=0.8178425167397448:2.323750136850048,variable=\t]({1.0*2.1497606536575806*cos(\t r)+-0.0*2.1497606536575806*sin(\t r)},{0.0*2.1497606536575806*cos(\t r)+1.0*2.1497606536575806*sin(\t r)});
\draw (5.74,5.58)-- (6.7,5.58);
\draw (7.14,6.56)-- (6.7,5.58);
\draw (7.14,6.56)-- (6.16,6.56);
\draw (5.020000000000007,2.2799999999999985) node[anchor=north west] {$R_{12,5}$};
\draw (-4.58,5.98)-- (-5.08,5.62);
\draw (-4.58,5.98)-- (-4.6,6.6);
\draw (-4.58,5.98)-- (-4.06,5.62);
\draw (-3.12,5.98)-- (-3.62,5.62);
\draw (-3.12,5.98)-- (-3.14,6.6);
\draw (-3.12,5.98)-- (-2.6,5.62);
\draw (-1.64,5.98)-- (-2.14,5.62);
\draw (-1.64,5.98)-- (-1.66,6.6);
\draw (-1.64,5.98)-- (-1.12,5.62);
\draw (-0.18,5.62)-- (0.32,5.96);
\draw (0.32,5.96)-- (0.3,6.6);
\draw (0.32,5.96)-- (0.84,5.62);
\draw (1.28,5.62)-- (1.78,5.98);
\draw (1.78,5.98)-- (1.76,6.6);
\draw (1.78,5.98)-- (2.3,5.62);
\draw (2.76,5.62)-- (3.28,5.98);
\draw (3.28,5.98)-- (3.78,5.62);
\draw (3.28,5.98)-- (3.24,6.6);
\draw (5.22,5.96)-- (5.2,6.56);
\draw (5.22,5.96)-- (4.72,5.58);
\draw (5.22,5.96)-- (5.74,5.58);
\draw (6.68,6.18)-- (6.7,5.58);
\draw (6.68,6.18)-- (7.14,6.56);
\draw (6.68,6.18)-- (6.16,6.56);
\draw (8.16,5.98)-- (8.14,6.56);
\draw (8.16,5.98)-- (7.66,5.58);
\draw (8.16,5.98)-- (8.68,5.58);
\begin{scriptsize}
\draw [fill=black] (-4.6,6.6) circle (2.5pt);
\draw [fill=black] (-4.06,5.62) circle (2.5pt);
\draw [fill=black] (-5.08,5.62) circle (2.5pt);
\draw [fill=black] (-3.14,6.6) circle (2.5pt);
\draw [fill=black] (-2.6,5.62) circle (2.5pt);
\draw [fill=black] (-3.62,5.62) circle (2.5pt);
\draw [fill=black] (-1.66,6.6) circle (2.5pt);
\draw [fill=black] (-1.12,5.62) circle (2.5pt);
\draw [fill=black] (-2.14,5.62) circle (2.5pt);
\draw [fill=black] (0.3,6.6) circle (2.5pt);
\draw [fill=black] (0.84,5.62) circle (2.5pt);
\draw [fill=black] (-0.18,5.62) circle (2.5pt);
\draw [fill=black] (1.76,6.6) circle (2.5pt);
\draw [fill=black] (2.3,5.62) circle (2.5pt);
\draw [fill=black] (1.28,5.62) circle (2.5pt);
\draw [fill=black] (3.24,6.6) circle (2.5pt);
\draw [fill=black] (3.78,5.62) circle (2.5pt);
\draw [fill=black] (2.76,5.62) circle (2.5pt);
\draw [fill=black] (5.2,6.56) circle (2.5pt);
\draw [fill=black] (5.74,5.58) circle (2.5pt);
\draw [fill=black] (4.72,5.58) circle (2.5pt);
\draw [fill=black] (6.16,6.56) circle (2.5pt);
\draw [fill=black] (6.7,5.58) circle (2.5pt);
\draw [fill=black] (8.14,6.56) circle (2.5pt);
\draw [fill=black] (8.68,5.58) circle (2.5pt);
\draw [fill=black] (7.66,5.58) circle (2.5pt);
\draw [fill=black] (7.14,6.56) circle (2.5pt);
\draw [fill=black] (-4.58,5.98) circle (2.5pt);
\draw [fill=black] (-3.12,5.98) circle (2.5pt);
\draw [fill=black] (-1.64,5.98) circle (2.5pt);
\draw [fill=black] (0.32,5.96) circle (2.5pt);
\draw [fill=black] (1.78,5.98) circle (2.5pt);
\draw [fill=black] (3.28,5.98) circle (2.5pt);
\draw [fill=black] (5.22,5.96) circle (2.5pt);
\draw [fill=black] (6.68,6.18) circle (2.5pt);
\draw [fill=black] (8.16,5.98) circle (2.5pt);
\end{scriptsize}
\end{tikzpicture}

%% file: graphl2b.tex
\begin{tikzpicture}[line cap=round,line join=round,>=triangle 45,x=1.0cm,y=1.0cm]
\clip(-5.5000000000000036,1.2999999999999985) rectangle (9.520000000000012,4.339999999999996);
\draw (-2.96,2.74)-- (-1.86,2.74);
\draw (-2.4,3.72)-- (-1.86,2.74);
\draw (-2.4,3.72)-- (-2.96,2.74);
\draw (-0.94,3.72)-- (-0.4,2.74);
\draw (0.06,2.74)-- (1.08,2.74);
\draw (0.54,3.72)-- (1.08,2.74);
\draw (0.54,3.72)-- (0.06,2.74);
\draw (-0.4,2.74)-- (0.06,2.74);
\draw (3.64,2.7)-- (4.66,2.7);
\draw (4.1,3.68)-- (4.66,2.7);
\draw (4.1,3.68)-- (3.64,2.7);
\draw (5.08,3.68)-- (5.62,2.7);
\draw (6.58,2.7)-- (7.6,2.7);
\draw (7.04,3.68)-- (7.6,2.7);
\draw (7.04,3.68)-- (6.58,2.7);
\draw (5.62,2.7)-- (6.58,2.7);
\draw (-1.4399999999999993,2.3999999999999977) node[anchor=north west] {$R_{12,4}$};
\draw (5.020000000000008,2.279999999999998) node[anchor=north west] {$R_{12,5}$};
\draw (4.66,2.7)-- (5.62,2.7);
\draw (6.06,3.68)-- (5.62,2.7);
\draw (6.06,3.68)-- (5.08,3.68);
\draw (-0.94,3.72)-- (-1.38,2.74);
\draw (-1.38,2.74)-- (-0.4,2.74);
\draw (-2.4,3.08)-- (-2.96,2.74);
\draw (-2.4,3.08)-- (-2.4,3.72);
\draw (-2.4,3.08)-- (-1.86,2.74);
\draw (-1.38,2.74)-- (-0.92,3.1);
\draw (-0.92,3.1)-- (-0.94,3.72);
\draw (-0.92,3.1)-- (-0.4,2.74);
\draw (0.06,2.74)-- (0.56,3.1);
\draw (0.56,3.1)-- (0.54,3.72);
\draw (0.56,3.1)-- (1.08,2.74);
\draw (3.64,2.7)-- (4.14,3.06);
\draw (4.14,3.06)-- (4.66,2.7);
\draw (4.14,3.06)-- (4.1,3.68);
\draw (5.62,2.7)-- (5.58,3.32);
\draw (5.58,3.32)-- (6.06,3.68);
\draw (5.58,3.32)-- (5.08,3.68);
\draw (7.06,3.08)-- (6.58,2.7);
\draw (7.06,3.08)-- (7.04,3.68);
\draw (7.06,3.08)-- (7.6,2.7);
\draw [shift={(5.62,3.9793750000000006)}] plot[domain=4.068653423330618:5.356124537438761,variable=\t]({1.0*1.5995000439590492*cos(\t r)+-0.0*1.5995000439590492*sin(\t r)},{0.0*1.5995000439590492*cos(\t r)+1.0*1.5995000439590492*sin(\t r)});
\begin{scriptsize}
\draw [fill=black] (-2.4,3.72) circle (2.5pt);
\draw [fill=black] (-1.86,2.74) circle (2.5pt);
\draw [fill=black] (-2.96,2.74) circle (2.5pt);
\draw [fill=black] (-0.94,3.72) circle (2.5pt);
\draw [fill=black] (-0.4,2.74) circle (2.5pt);
\draw [fill=black] (0.54,3.72) circle (2.5pt);
\draw [fill=black] (1.08,2.74) circle (2.5pt);
\draw [fill=black] (0.06,2.74) circle (2.5pt);
\draw [fill=black] (4.1,3.68) circle (2.5pt);
\draw [fill=black] (4.66,2.7) circle (2.5pt);
\draw [fill=black] (3.64,2.7) circle (2.5pt);
\draw [fill=black] (5.08,3.68) circle (2.5pt);
\draw [fill=black] (5.62,2.7) circle (2.5pt);
\draw [fill=black] (7.04,3.68) circle (2.5pt);
\draw [fill=black] (7.6,2.7) circle (2.5pt);
\draw [fill=black] (6.58,2.7) circle (2.5pt);
\draw [fill=black] (6.06,3.68) circle (2.5pt);
\draw [fill=black] (-1.38,2.74) circle (2.5pt);
\draw [fill=black] (-2.4,3.08) circle (2.5pt);
\draw [fill=black] (-0.92,3.1) circle (2.5pt);
\draw [fill=black] (0.56,3.1) circle (2.5pt);
\draw [fill=black] (4.14,3.06) circle (2.5pt);
\draw [fill=black] (5.58,3.32) circle (2.5pt);
\draw [fill=black] (7.06,3.08) circle (2.5pt);
\end{scriptsize}
\end{tikzpicture}

%% file: graph71.tex
\begin{tikzpicture}[line cap=round,line join=round,>=triangle 45,x=0.8316498316498296cm,y=1.0cm]
\clip(3.0199999999999694,0.4999999999999986) rectangle (20.839999999999772,6.279999999999996);
\draw (5.859223164628248,5.799618603959381)-- (5.04578153541165,5.1838973067335585);
\draw (5.38,4.22)-- (5.04578153541165,5.1838973067335585);
\draw (5.38,4.22)-- (6.4,4.24);
\draw (6.4,4.24)-- (6.6961762039365444,5.216257986508556);
\draw (6.6961762039365444,5.216257986508556)-- (5.859223164628248,5.799618603959381);
\draw (5.859223164628248,5.799618603959381)-- (5.38,4.22);
\draw (5.38,4.22)-- (6.6961762039365444,5.216257986508556);
\draw (6.6961762039365444,5.216257986508556)-- (5.04578153541165,5.1838973067335585);
\draw (5.04578153541165,5.1838973067335585)-- (6.4,4.24);
\draw (6.4,4.24)-- (5.859223164628248,5.799618603959381);
\draw (4.180776835371754,3.501172274702884)-- (3.3114631162884596,2.8943002471603623);
\draw (3.62,1.88)-- (3.3114631162884596,2.8943002471603623);
\draw (3.62,1.88)-- (4.68,1.86);
\draw (4.68,1.86)-- (5.026579144363348,2.8619395673853627);
\draw (5.026579144363348,2.8619395673853627)-- (4.180776835371754,3.501172274702884);
\draw (4.180776835371754,3.501172274702884)-- (3.62,1.88);
\draw (3.62,1.88)-- (5.026579144363348,2.8619395673853627);
\draw (5.026579144363348,2.8619395673853627)-- (3.3114631162884596,2.8943002471603623);
\draw (3.3114631162884596,2.8943002471603623)-- (4.68,1.86);
\draw (4.68,1.86)-- (4.180776835371754,3.501172274702884);
\draw (7.3607768353717535,3.4596186039593784)-- (6.523823796063457,2.8762579865085556);
\draw (6.82,1.9)-- (6.523823796063457,2.8762579865085556);
\draw (6.82,1.9)-- (7.84,1.88);
\draw (7.84,1.88)-- (8.17421846458835,2.8438973067335565);
\draw (8.17421846458835,2.8438973067335565)-- (7.3607768353717535,3.4596186039593784);
\draw (7.3607768353717535,3.4596186039593784)-- (6.82,1.9);
\draw (6.82,1.9)-- (8.17421846458835,2.8438973067335565);
\draw (8.17421846458835,2.8438973067335565)-- (6.523823796063457,2.8762579865085556);
\draw (6.523823796063457,2.8762579865085556)-- (7.84,1.88);
\draw (7.84,1.88)-- (7.3607768353717535,3.4596186039593784);
\draw (5.38,4.22)-- (5.026579144363348,2.8619395673853627);
\draw (5.38,4.22)-- (6.523823796063457,2.8762579865085556);
\draw (11.839223164628123,5.761172274702883)-- (10.99342085563653,5.121939567385363);
\draw (11.339999999999877,4.12)-- (10.99342085563653,5.121939567385363);
\draw (11.339999999999877,4.12)-- (12.399999999999878,4.14);
\draw (12.399999999999878,4.14)-- (12.708536883711417,5.154300247160362);
\draw (12.708536883711417,5.154300247160362)-- (11.839223164628123,5.761172274702883);
\draw (11.839223164628123,5.761172274702883)-- (11.339999999999877,4.12);
\draw (11.339999999999877,4.12)-- (12.708536883711417,5.154300247160362);
\draw (12.708536883711417,5.154300247160362)-- (10.99342085563653,5.121939567385363);
\draw (10.99342085563653,5.121939567385363)-- (12.399999999999878,4.14);
\draw (12.399999999999878,4.14)-- (11.839223164628123,5.761172274702883);
\draw (10.119223164628137,3.4811722747028577)-- (9.273420855636555,2.841939567385348);
\draw (9.619999999999896,1.84)-- (9.273420855636555,2.841939567385348);
\draw (9.619999999999896,1.84)-- (10.679999999999879,1.86);
\draw (10.679999999999879,1.86)-- (10.988536883711415,2.8743002471603445);
\draw (10.988536883711415,2.8743002471603445)-- (10.119223164628137,3.4811722747028577);
\draw (10.119223164628137,3.4811722747028577)-- (9.619999999999896,1.84);
\draw (9.619999999999896,1.84)-- (10.988536883711415,2.8743002471603445);
\draw (10.988536883711415,2.8743002471603445)-- (9.273420855636555,2.841939567385348);
\draw (9.273420855636555,2.841939567385348)-- (10.679999999999879,1.86);
\draw (10.679999999999879,1.86)-- (10.119223164628137,3.4811722747028577);
\draw (13.609223164628114,3.4703954393311625)-- (12.779601195524002,2.84291843705948);
\draw (13.119999999999857,1.86)-- (12.779601195524002,2.84291843705948);
\draw (13.119999999999857,1.86)-- (14.159999999999878,1.88);
\draw (14.159999999999878,1.88)-- (14.462356543823926,2.875279116834478);
\draw (14.462356543823926,2.875279116834478)-- (13.609223164628114,3.4703954393311625);
\draw (13.609223164628114,3.4703954393311625)-- (13.119999999999857,1.86);
\draw (13.119999999999857,1.86)-- (14.462356543823926,2.875279116834478);
\draw (14.462356543823926,2.875279116834478)-- (12.779601195524002,2.84291843705948);
\draw (12.779601195524002,2.84291843705948)-- (14.159999999999878,1.88);
\draw (14.159999999999878,1.88)-- (13.609223164628114,3.4703954393311625);
\draw (11.339999999999877,4.12)-- (10.988536883711415,2.8743002471603445);
\draw (18.009223164628057,5.751949110074644)-- (17.147240515748962,5.100960697711272);
\draw (17.499999999999808,4.08)-- (17.147240515748962,5.100960697711272);
\draw (17.499999999999808,4.08)-- (18.579999999999814,4.1);
\draw (18.579999999999814,4.1)-- (18.894717223598853,5.133321377486269);
\draw (18.894717223598853,5.133321377486269)-- (18.009223164628057,5.751949110074644);
\draw (18.009223164628057,5.751949110074644)-- (17.499999999999808,4.08);
\draw (17.499999999999808,4.08)-- (18.894717223598853,5.133321377486269);
\draw (18.894717223598853,5.133321377486269)-- (17.147240515748962,5.100960697711272);
\draw (17.147240515748962,5.100960697711272)-- (18.579999999999814,4.1);
\draw (18.579999999999814,4.1)-- (18.009223164628057,5.751949110074644);
\draw (16.249223164628074,3.5319491100746188)-- (15.387240515748987,2.8809606977112567);
\draw (15.739999999999828,1.86)-- (15.387240515748987,2.8809606977112567);
\draw (15.739999999999828,1.86)-- (16.819999999999816,1.88);
\draw (16.819999999999816,1.88)-- (17.134717223598855,2.913321377486252);
\draw (17.134717223598855,2.913321377486252)-- (16.249223164628074,3.5319491100746188);
\draw (16.249223164628074,3.5319491100746188)-- (15.739999999999828,1.86);
\draw (15.739999999999828,1.86)-- (17.134717223598855,2.913321377486252);
\draw (17.134717223598855,2.913321377486252)-- (15.387240515748987,2.8809606977112567);
\draw (15.387240515748987,2.8809606977112567)-- (16.819999999999816,1.88);
\draw (16.819999999999816,1.88)-- (16.249223164628074,3.5319491100746188);
\draw (19.309999999999807,3.5096186039594164)-- (18.48480266573734,2.91007764662108);
\draw (18.799999999999795,1.94)-- (18.48480266573734,2.91007764662108);
\draw (18.799999999999795,1.94)-- (19.81999999999982,1.94);
\draw (19.81999999999982,1.94)-- (20.135197334262273,2.910077646621079);
\draw (20.135197334262273,2.910077646621079)-- (19.309999999999807,3.5096186039594164);
\draw (19.309999999999807,3.5096186039594164)-- (18.799999999999795,1.94);
\draw (18.799999999999795,1.94)-- (20.135197334262273,2.910077646621079);
\draw (20.135197334262273,2.910077646621079)-- (18.48480266573734,2.91007764662108);
\draw (18.48480266573734,2.91007764662108)-- (19.81999999999982,1.94);
\draw (19.81999999999982,1.94)-- (19.309999999999807,3.5096186039594164);
\draw (17.499999999999808,4.08)-- (17.134717223598855,2.913321377486252);
\draw (17.499999999999808,4.08)-- (18.48480266573734,2.91007764662108);
\draw (12.399999999999878,4.14)-- (12.779601195524002,2.84291843705948);
\draw (17.134717223598855,2.913321377486252)-- (18.48480266573734,2.91007764662108);
\draw (5.259999999999944,1.539999999999998) node[anchor=north west] {$R_{15,1}$};
\draw (11.359999999999877,1.579999999999998) node[anchor=north west] {$R_{15,2}$};
\draw (17.41999999999981,1.559999999999998) node[anchor=north west] {$R_{15,3}$};
\begin{scriptsize}
\draw [fill=black] (5.38,4.22) circle (1.5pt);
\draw [fill=black] (6.4,4.24) circle (1.5pt);
\draw [fill=black] (6.6961762039365444,5.216257986508556) circle (1.5pt);
\draw [fill=black] (5.859223164628248,5.799618603959381) circle (1.5pt);
\draw [fill=black] (5.04578153541165,5.1838973067335585) circle (1.5pt);
\draw [fill=black] (3.62,1.88) circle (1.5pt);
\draw [fill=black] (4.68,1.86) circle (1.5pt);
\draw [fill=black] (5.026579144363348,2.8619395673853627) circle (1.5pt);
\draw [fill=black] (4.180776835371754,3.501172274702884) circle (1.5pt);
\draw [fill=black] (3.3114631162884596,2.8943002471603623) circle (1.5pt);
\draw [fill=black] (6.82,1.9) circle (1.5pt);
\draw [fill=black] (7.84,1.88) circle (1.5pt);
\draw [fill=black] (8.17421846458835,2.8438973067335565) circle (1.5pt);
\draw [fill=black] (7.3607768353717535,3.4596186039593784) circle (1.5pt);
\draw [fill=black] (6.523823796063457,2.8762579865085556) circle (1.5pt);
\draw [fill=black] (11.339999999999877,4.12) circle (1.5pt);
\draw [fill=black] (12.399999999999878,4.14) circle (1.5pt);
\draw [fill=black] (12.708536883711417,5.154300247160362) circle (1.5pt);
\draw [fill=black] (11.839223164628123,5.761172274702883) circle (1.5pt);
\draw [fill=black] (10.99342085563653,5.121939567385363) circle (1.5pt);
\draw [fill=black] (9.619999999999896,1.84) circle (1.5pt);
\draw [fill=black] (10.679999999999879,1.86) circle (1.5pt);
\draw [fill=black] (10.988536883711415,2.8743002471603445) circle (1.5pt);
\draw [fill=black] (10.119223164628137,3.4811722747028577) circle (1.5pt);
\draw [fill=black] (9.273420855636555,2.841939567385348) circle (1.5pt);
\draw [fill=black] (13.119999999999857,1.86) circle (1.5pt);
\draw [fill=black] (14.159999999999878,1.88) circle (1.5pt);
\draw [fill=black] (14.462356543823926,2.875279116834478) circle (1.5pt);
\draw [fill=black] (13.609223164628114,3.4703954393311625) circle (1.5pt);
\draw [fill=black] (12.779601195524002,2.84291843705948) circle (1.5pt);
\draw [fill=black] (17.499999999999808,4.08) circle (1.5pt);
\draw [fill=black] (18.579999999999814,4.1) circle (1.5pt);
\draw [fill=black] (18.894717223598853,5.133321377486269) circle (1.5pt);
\draw [fill=black] (18.009223164628057,5.751949110074644) circle (1.5pt);
\draw [fill=black] (17.147240515748962,5.100960697711272) circle (1.5pt);
\draw [fill=black] (15.739999999999828,1.86) circle (1.5pt);
\draw [fill=black] (16.819999999999816,1.88) circle (1.5pt);
\draw [fill=black] (17.134717223598855,2.913321377486252) circle (1.5pt);
\draw [fill=black] (16.249223164628074,3.5319491100746188) circle (1.5pt);
\draw [fill=black] (15.387240515748987,2.8809606977112567) circle (1.5pt);
\draw [fill=black] (18.799999999999795,1.94) circle (1.5pt);
\draw [fill=black] (19.81999999999982,1.94) circle (1.5pt);
\draw [fill=black] (20.135197334262273,2.910077646621079) circle (1.5pt);
\draw [fill=black] (19.309999999999807,3.5096186039594164) circle (1.5pt);
\draw [fill=black] (18.48480266573734,2.91007764662108) circle (1.5pt);
\end{scriptsize}
\end{tikzpicture}

%% file: graphl31.tex
\definecolor{zzttqq}{rgb}{0.6,0.2,0.0}
\begin{tikzpicture}[line cap=round,line join=round,>=triangle 45,x=1.0cm,y=1.0cm]
\clip(-5.519690617822582,1.5149457986058634) rectangle (1.8162095498877742,7.368461994602707);
\fill[color=zzttqq,fill=zzttqq,fill opacity=0.1] (-0.19831225782543246,7.254432458317054) -- (-0.17930733511115693,6.722294622317341) -- (0.35283050088855555,6.741299545031617) -- (0.3338255781742802,7.2734373810313295) -- cycle;
\fill[color=zzttqq,fill=zzttqq,fill opacity=0.1] (-0.236322103253984,6.095132172746251) -- (-0.21731718053970855,5.562994336746538) -- (0.31482065546000393,5.581999259460813) -- (0.29581573274572864,6.114137095460526) -- cycle;
\fill[color=zzttqq,fill=zzttqq,fill opacity=0.1] (-0.2363221032539838,5.068866346175375) -- (-0.2173171805397078,4.536728510175663) -- (0.3148206554600047,4.555733432889938) -- (0.29581573274572914,5.087871268889651) -- cycle;
\fill[color=zzttqq,fill=zzttqq,fill opacity=0.1] (-0.10328764425405487,3.4154380700334097) -- (-0.08428272153977942,2.883300234033697) -- (0.4478551144599331,2.9023051567479725) -- (0.42885019174565775,3.434442992747685) -- cycle;
\fill[color=zzttqq,fill=zzttqq,fill opacity=0.1] (-5.386656158822651,5.068866346175375) -- (-5.367651236108376,4.536728510175663) -- (-4.835513400108663,4.555733432889938) -- (-4.854518322822939,5.087871268889651) -- cycle;
\draw (-3.258104814823792,4.859812196318345)-- (-3.0680555876810365,5.791053409317843);
\draw (-3.2771097375380673,6.608265086031688)-- (-3.0680555876810365,5.791053409317843);
\draw (-3.2771097375380673,6.608265086031688)-- (-3.258104814823792,4.859812196318345);
\draw (-3.2581048148237923,4.003694051175736)-- (-3.239099892109517,3.1864823744618898);
\draw [rotate around={-88.78112476486909:(-3.267607276180931,4.850309734961213)}] (-3.267607276180931,4.850309734961213) ellipse (2.333163428367182cm and 0.6743552024673358cm);
\draw (-5.101582318108518,4.802797428175519)-- (-3.2771097375380673,6.608265086031688);
\draw (-5.101582318108518,4.802797428175519)-- (-3.0680555876810365,5.791053409317843);
\draw (-5.101582318108518,4.802797428175519)-- (-3.258104814823792,4.859812196318345);
\draw (-5.101582318108518,4.802797428175519)-- (-3.2581048148237923,4.003694051175736);
\draw (-5.101582318108518,4.802797428175519)-- (-3.239099892109517,3.1864823744618898);
\draw (-3.2771097375380673,6.608265086031688)-- (0.048751737460149464,6.988363540317198);
\draw (-3.2771097375380673,6.608265086031688)-- (0.048751737460149464,6.456225704317484);
\draw (-3.0680555876810365,5.791053409317843)-- (0.04875173746014945,5.829063254746394);
\draw (-3.0680555876810365,5.791053409317843)-- (0.06775666017442498,5.315930341460956);
\draw (-3.258104814823792,4.859812196318345)-- (0.048751737460149464,4.802797428175519);
\draw (-3.258104814823792,4.859812196318345)-- (0.06775666017442492,4.289664514890081);
\draw (0.18178619646007813,3.1493691520335534)-- (-3.2581048148237923,4.003694051175736);
\draw (-3.239099892109517,3.1864823744618898)-- (0.18178619646007813,3.1493691520335534);
\draw (-0.008263030682677116,0.2226110540351318) node[anchor=north west] {$v_1$};
\draw (0.5428797280313133,3.491457760890512) node[anchor=north west] {$v$};
\draw (0.48586495988848677,7.235427535602779) node[anchor=north west] {$v_1$};
\draw (0.5428797280313133,6.076127250031975) node[anchor=north west] {$v_2$};
\draw (0.5428797280313133,5.0308565007468244) node[anchor=north west] {$v_3$};
\draw (-0.008263030682677116,0.2226110540351318) node[anchor=north west] {$v$};
\draw (-5.386656158822653,4.517723587461387) node[anchor=north west] {$w$};
\draw (-3.733227882680682,2.4081771661768103) node[anchor=north west] {$\Gamma(w)$};
\draw [color=zzttqq] (-0.19831225782543246,7.254432458317054)-- (-0.17930733511115693,6.722294622317341);
\draw [color=zzttqq] (-0.17930733511115693,6.722294622317341)-- (0.35283050088855555,6.741299545031617);
\draw [color=zzttqq] (0.35283050088855555,6.741299545031617)-- (0.3338255781742802,7.2734373810313295);
\draw [color=zzttqq] (0.3338255781742802,7.2734373810313295)-- (-0.19831225782543246,7.254432458317054);
\draw [color=zzttqq] (-0.236322103253984,6.095132172746251)-- (-0.21731718053970855,5.562994336746538);
\draw [color=zzttqq] (-0.21731718053970855,5.562994336746538)-- (0.31482065546000393,5.581999259460813);
\draw [color=zzttqq] (0.31482065546000393,5.581999259460813)-- (0.29581573274572864,6.114137095460526);
\draw [color=zzttqq] (0.29581573274572864,6.114137095460526)-- (-0.236322103253984,6.095132172746251);
\draw [color=zzttqq] (-0.2363221032539838,5.068866346175375)-- (-0.2173171805397078,4.536728510175663);
\draw [color=zzttqq] (-0.2173171805397078,4.536728510175663)-- (0.3148206554600047,4.555733432889938);
\draw [color=zzttqq] (0.3148206554600047,4.555733432889938)-- (0.29581573274572914,5.087871268889651);
\draw [color=zzttqq] (0.29581573274572914,5.087871268889651)-- (-0.2363221032539838,5.068866346175375);
\draw [color=zzttqq] (-0.10328764425405487,3.4154380700334097)-- (-0.08428272153977942,2.883300234033697);
\draw [color=zzttqq] (-0.08428272153977942,2.883300234033697)-- (0.4478551144599331,2.9023051567479725);
\draw [color=zzttqq] (0.4478551144599331,2.9023051567479725)-- (0.42885019174565775,3.434442992747685);
\draw [color=zzttqq] (0.42885019174565775,3.434442992747685)-- (-0.10328764425405487,3.4154380700334097);
\draw [color=zzttqq] (-5.386656158822651,5.068866346175375)-- (-5.367651236108376,4.536728510175663);
\draw [color=zzttqq] (-5.367651236108376,4.536728510175663)-- (-4.835513400108663,4.555733432889938);
\draw [color=zzttqq] (-4.835513400108663,4.555733432889938)-- (-4.854518322822939,5.087871268889651);
\draw [color=zzttqq] (-4.854518322822939,5.087871268889651)-- (-5.386656158822651,5.068866346175375);
\begin{scriptsize}
\draw [fill=black] (-3.2771097375380673,6.608265086031688) circle (2.5pt);
\draw [fill=black] (-3.0680555876810365,5.791053409317843) circle (2.5pt);
\draw [fill=black] (-3.258104814823792,4.859812196318345) circle (2.5pt);
\draw [fill=black] (-3.2581048148237923,4.003694051175736) circle (2.5pt);
\draw [fill=black] (-3.239099892109517,3.1864823744618898) circle (2.5pt);
\draw [fill=black] (-5.101582318108518,4.802797428175519) circle (2.5pt);
\draw [fill=black] (0.04875173746014945,5.829063254746394) circle (2.5pt);
\draw [fill=black] (0.06775666017442498,5.315930341460956) circle (2.5pt);
\draw [fill=black] (0.048751737460149464,4.802797428175519) circle (2.5pt);
\draw [fill=black] (0.06775666017442492,4.289664514890081) circle (2.5pt);
\draw [fill=black] (0.08676158288870052,6.988363540317198) circle (2.5pt);
\draw [fill=black] (0.10576650560297604,6.418215858888933) circle (2.5pt);
\draw [fill=black] (0.18178619646007813,3.1493691520335534) circle (2.5pt);
\end{scriptsize}
\end{tikzpicture}

%% file: graphl41.tex
\begin{tikzpicture}[line cap=round,line join=round,>=triangle 45,x=1.0cm,y=1.0cm]
\clip(-7.060000000000001,-0.3999999999999975) rectangle (7.740000000000006,6.159999999999999);
\draw (-2.03689265851299,4.921186317993395)-- (0.4402175437534477,3.022574971290883);
\draw (-0.6,0.08)-- (0.4402175437534477,3.022574971290883);
\draw (-0.6,0.08)-- (-3.72,0.16);
\draw (-3.72,0.16)-- (-4.608048501146223,3.1520176903908754);
\draw (-4.608048501146223,3.1520176903908754)-- (-2.03689265851299,4.921186317993395);
\draw (3.06,3.6)-- (-2.03689265851299,4.921186317993395);
\draw (3.06,3.6)-- (-0.6,0.08);
\draw (5.279999999999999,1.2600000000000005)-- (-4.608048501146223,3.1520176903908754);
\draw (2.880000000000004,4.54) node[anchor=north west] {$y_1$};
\draw (5.360000000000005,4.36) node[anchor=north west] {$y_2$};
\draw (2.880000000000004,0.9400000000000018) node[anchor=north west] {$y_3$};
\draw (5.260000000000005,1.0000000000000018) node[anchor=north west] {$y_4$};
\draw (0.30000000000000243,3.9200000000000004) node[anchor=north west] {$u_2$};
\draw (-3.1799999999999993,5.499999999999999) node[anchor=north west] {$u_3$};
\draw (-5.460000000000001,3.5400000000000005) node[anchor=north west] {$u_4$};
\draw (-4.54,0.9600000000000017) node[anchor=north west] {$u_5$};
\draw (-0.059999999999997757,0.520000000000002) node[anchor=north west] {$u_1$};
\draw(3.06,3.6) circle (0.3cm);
\draw [shift={(2.1107003261646624,10.66193514793171)}] plot[domain=4.497109640911972:5.138091705272323,variable=\t]({1.0*7.819868070164488*cos(\t r)+-0.0*7.819868070164488*sin(\t r)},{0.0*7.819868070164488*cos(\t r)+1.0*7.819868070164488*sin(\t r)});
\draw(3.0,1.3200000000000014) circle (0.3cm);
\draw(5.279999999999999,1.2600000000000005) circle (0.3cm);
\draw(0.44,3.02) circle (0.3cm);
\draw (-2.03689265851299,4.921186317993395)-- (-0.6,0.08);
\draw (5.34,3.54)-- (3.7,4.68);
\draw (5.34,3.54)-- (4.2,4.98);
\draw (5.34,3.54)-- (4.72,5.22);
\draw (5.34,3.54)-- (5.4,5.44);
\draw (2.840000000000004,5.279999999999999) node[anchor=north west] {$w_1$};
\draw (3.400000000000004,5.739999999999999) node[anchor=north west] {$w_2$};
\draw (3.9600000000000044,6.1) node[anchor=north west] {$w_3$};
\draw (4.960000000000005,6.299999999999999) node[anchor=north west] {$w_4$};
\draw [shift={(-8.292010050251228,46.69095477386913)}] plot[domain=4.810331977861224:4.9563150504617735,variable=\t]({1.0*46.7550321149227*cos(\t r)+-0.0*46.7550321149227*sin(\t r)},{0.0*46.7550321149227*cos(\t r)+1.0*46.7550321149227*sin(\t r)});
\begin{scriptsize}
\draw [fill=black] (-3.72,0.16) circle (1.5pt);
\draw [fill=black] (-0.6,0.08) circle (1.5pt);
\draw [fill=black] (0.4402175437534477,3.022574971290883) circle (1.5pt);
\draw [fill=black] (-2.03689265851299,4.921186317993395) circle (1.5pt);
\draw [fill=black] (-4.608048501146223,3.1520176903908754) circle (1.5pt);
\draw [fill=black] (5.34,3.54) circle (1.5pt);
\draw [fill=black] (3.06,3.6) circle (1.5pt);
\draw [fill=black] (3.0,1.3200000000000014) circle (1.5pt);
\draw [fill=black] (5.279999999999999,1.2600000000000005) circle (1.5pt);
\draw [fill=black] (0.44,3.02) circle (1.5pt);
\draw [fill=black] (3.7,4.68) circle (1.5pt);
\draw [fill=black] (4.2,4.98) circle (1.5pt);
\draw [fill=black] (4.72,5.22) circle (1.5pt);
\draw [fill=black] (5.4,5.44) circle (1.5pt);
\end{scriptsize}
\end{tikzpicture}

%% file: graphl42.tex
\begin{tikzpicture}[line cap=round,line join=round,>=triangle 45,x=1.0cm,y=1.0cm]
\clip(-6.940000000000001,-0.6799999999999974) rectangle (7.180000000000004,5.359999999999999);
\draw (-2.026892658512989,4.890409482621642)-- (0.43403720386594946,3.0035538409649796);
\draw (-0.6,0.08)-- (0.43403720386594946,3.0035538409649796);
\draw (-0.6,0.08)-- (-3.7,0.16);
\draw (-3.7,0.16)-- (-4.581868161258725,3.1329965600649725);
\draw (-4.581868161258725,3.1329965600649725)-- (-2.026892658512989,4.890409482621642);
\draw (3.06,3.6)-- (-2.026892658512989,4.890409482621642);
\draw (3.06,3.6)-- (-0.6,0.08);
\draw [shift={(-1.0637678207739312,4.862624482056323)}] plot[domain=4.201453578468767:5.56619410537268,variable=\t]({1.0*5.3911406141949705*cos(\t r)+-0.0*5.3911406141949705*sin(\t r)},{0.0*5.3911406141949705*cos(\t r)+1.0*5.3911406141949705*sin(\t r)});
\draw (5.279999999999999,1.2600000000000005)-- (-4.581868161258725,3.1329965600649725);
\draw (2.880000000000002,4.54) node[anchor=north west] {$y_1$};
\draw (5.160000000000003,4.48) node[anchor=north west] {$y_2$};
\draw (2.880000000000002,0.9400000000000018) node[anchor=north west] {$y_3$};
\draw (5.260000000000003,1.0000000000000018) node[anchor=north west] {$y_4$};
\draw (0.30000000000000127,3.9200000000000004) node[anchor=north west] {$u_2$};
\draw (-3.18,5.499999999999999) node[anchor=north west] {$u_3$};
\draw (-5.460000000000001,3.5400000000000005) node[anchor=north west] {$u_4$};
\draw (-4.540000000000001,0.9600000000000017) node[anchor=north west] {$u_5$};
\draw (-0.05999999999999889,0.520000000000002) node[anchor=north west] {$u_1$};
\draw (-6.580000000000001,0.9600000000000017) node[anchor=north west] {$w$};
\draw(3.06,3.6) circle (0.3cm);
\draw [shift={(2.0570844994575936,10.861775956567147)}] plot[domain=4.508711707917064:5.133892349189461,variable=\t]({1.0*8.02408482893011*cos(\t r)+-0.0*8.02408482893011*sin(\t r)},{0.0*8.02408482893011*cos(\t r)+1.0*8.02408482893011*sin(\t r)});
\draw(5.34,3.54) circle (0.3cm);
\draw(3.0,1.3200000000000014) circle (0.3cm);
\draw(5.279999999999999,1.2600000000000005) circle (0.3cm);
\draw(-5.82,0.16) circle (0.3cm);
\draw (-3.7,0.16)-- (-5.82,0.16);
\draw [shift={(1.9050333915468993,5.102327292509178)}] plot[domain=3.195437249689105:4.99418094679084,variable=\t]({1.0*3.9376327433227036*cos(\t r)+-0.0*3.9376327433227036*sin(\t r)},{0.0*3.9376327433227036*cos(\t r)+1.0*3.9376327433227036*sin(\t r)});
\begin{scriptsize}
\draw [fill=black] (-3.7,0.16) circle (1.5pt);
\draw [fill=black] (-0.6,0.08) circle (1.5pt);
\draw [fill=black] (0.43403720386594946,3.0035538409649796) circle (1.5pt);
\draw [fill=black] (-2.026892658512989,4.890409482621642) circle (1.5pt);
\draw [fill=black] (-4.581868161258725,3.1329965600649725) circle (1.5pt);
\draw [fill=black] (5.34,3.54) circle (1.5pt);
\draw [fill=black] (3.06,3.6) circle (1.5pt);
\draw [fill=black] (3.0,1.3200000000000014) circle (1.5pt);
\draw [fill=black] (5.279999999999999,1.2600000000000005) circle (1.5pt);
\draw [fill=black] (-5.82,0.16) circle (1.5pt);
\end{scriptsize}
\end{tikzpicture}

%% file: graph1.tex
\begin{tikzpicture}[line cap=round,line join=round,>=triangle 45,x=1.0cm,y=1.0cm]
\clip(-5.84,3.999999999999992) rectangle (5.460000000000004,10.619999999999987);
\draw (0.0800000000000023,6.49999999999999) node[anchor=north west] {$u_{1}$};
\draw (-4.02,10.439999999999987) node[anchor=north west] {$u_{4}$};
\draw (-4.819999999999999,1.8199999999999938) node[anchor=north west] {$v_{1,3}$};
\draw (0.28000000000000236,8.879999999999988) node[anchor=north west] {$u_{2}$};
\draw (-1.019999999999998,10.379999999999987) node[anchor=north west] {$u_{3}$};
\draw (-0.639999999999998,1.6399999999999941) node[anchor=north west] {$v_{ 2,3}$};
\draw (1.98,3.5)-- (3.38,3.5);
\draw (2.2737367544323214E-15,0.23999999999999524) node[anchor=north west] {$v_{2,1}$};
\draw (-5.34,9.139999999999988) node[anchor=north west] {$u_{5}$};
\draw (-5.38,6.6199999999999894) node[anchor=north west] {$u_{6}$};
\draw (-3.959999999999999,5.059999999999991) node[anchor=north west] {$u_{7}$};
\draw (-1.4599999999999982,5.019999999999992) node[anchor=north west] {$u_{8}$};
\draw (2.460000000000003,9.419999999999987) node[anchor=north west] {$y_{1}$};
\draw (5.000000000000004,9.359999999999987) node[anchor=north west] {$y_{2}$};
\draw (2.500000000000003,6.09999999999999) node[anchor=north west] {$y_{4}$};
\draw (4.780000000000004,6.15999999999999) node[anchor=north west] {$y_{3}$};
\draw (2.660000000000001,8.620000000000001)-- (-1.3399999999999999,9.863574311003802);
\draw (2.660000000000001,8.620000000000001)-- (0.03178715550190114,6.551787155501901);
\draw (0.03178715550190114,6.551787155501901)-- (-1.3399999999999999,9.863574311003802);
\draw [shift={(2.5288458004173706,6.819474528448952)}] plot[domain=0.6507309768653503:2.551484016792735,variable=\t]({1.0*3.0053171877104443*cos(\t r)+-0.0*3.0053171877104443*sin(\t r)},{0.0*3.0053171877104443*cos(\t r)+1.0*3.0053171877104443*sin(\t r)});
\draw (4.920000000000002,8.639999999999999)-- (-4.651787155501902,6.551787155501902);
\draw (2.68,6.36)-- (-1.34,5.18);
\draw (2.68,6.36)-- (-3.2799999999999994,9.863574311003802);
\draw (4.94,6.38)-- (-4.651787155501902,8.491787155501902);
\draw (-3.2799999999999994,9.863574311003802)-- (-1.34,5.18);
\draw (-4.651787155501902,8.491787155501902)-- (-3.28,5.18);
\draw (-3.2799999999999994,9.863574311003802)-- (-1.3399999999999999,9.863574311003802);
\draw (0.03178715550190114,8.4917871555019)-- (-1.3399999999999999,9.863574311003802);
\draw (0.03178715550190114,8.4917871555019)-- (0.03178715550190114,6.551787155501901);
\draw (0.03178715550190114,6.551787155501901)-- (-1.34,5.18);
\draw (-1.34,5.18)-- (-3.28,5.18);
\draw (-3.28,5.18)-- (-4.651787155501902,6.551787155501902);
\draw (-4.651787155501902,6.551787155501902)-- (-4.651787155501902,8.491787155501902);
\draw (-4.651787155501902,8.491787155501902)-- (-3.2799999999999994,9.863574311003802);
\draw [shift={(0.11706568034073109,10.663600089665996)}] plot[domain=4.157757739873641:5.556943013678742,variable=\t]({1.0*6.450575569663047*cos(\t r)+-0.0*6.450575569663047*sin(\t r)},{0.0*6.450575569663047*cos(\t r)+1.0*6.450575569663047*sin(\t r)});
\draw (0.03178715550190114,8.4917871555019)-- (-4.651787155501902,6.551787155501902);
\begin{scriptsize}
\draw [fill=black] (1.98,3.5) circle (1.5pt);
\draw [fill=black] (3.38,3.5) circle (1.5pt);
\draw [fill=black] (-3.28,5.18) circle (1.5pt);
\draw [fill=black] (-1.34,5.18) circle (1.5pt);
\draw [fill=black] (0.03178715550190114,6.551787155501901) circle (1.5pt);
\draw [fill=black] (0.03178715550190114,8.4917871555019) circle (1.5pt);
\draw [fill=black] (-1.3399999999999999,9.863574311003802) circle (1.5pt);
\draw [fill=black] (-3.2799999999999994,9.863574311003802) circle (1.5pt);
\draw [fill=black] (-4.651787155501902,8.491787155501902) circle (1.5pt);
\draw [fill=black] (-4.651787155501902,6.551787155501902) circle (1.5pt);
\draw [fill=black] (4.94,6.38) circle (1.5pt);
\draw [fill=black] (2.68,6.36) circle (1.5pt);
\draw [fill=black] (4.920000000000002,8.639999999999999) circle (1.5pt);
\draw [fill=black] (2.660000000000001,8.620000000000001) circle (1.5pt);
\end{scriptsize}
\end{tikzpicture}

%% file: graph2.tex
\begin{tikzpicture}[line cap=round,line join=round,>=triangle 45,x=1.0cm,y=1.0cm]
\clip(-5.84,3.999999999999992) rectangle (5.460000000000004,10.619999999999987);
\draw (0.0800000000000023,6.49999999999999) node[anchor=north west] {$u_{1}$};
\draw (-4.02,10.439999999999987) node[anchor=north west] {$u_{4}$};
\draw (-4.819999999999999,1.8199999999999938) node[anchor=north west] {$v_{1,3}$};
\draw (0.28000000000000236,8.879999999999988) node[anchor=north west] {$u_{2}$};
\draw (-1.019999999999998,10.379999999999987) node[anchor=north west] {$u_{3}$};
\draw (-0.639999999999998,1.6399999999999941) node[anchor=north west] {$v_{ 2,3}$};
\draw (1.98,3.5)-- (3.38,3.5);
\draw (2.2737367544323214E-15,0.23999999999999524) node[anchor=north west] {$v_{2,1}$};
\draw (-5.34,9.139999999999988) node[anchor=north west] {$u_{5}$};
\draw (-5.38,6.6199999999999894) node[anchor=north west] {$u_{6}$};
\draw (-3.959999999999999,5.059999999999991) node[anchor=north west] {$u_{7}$};
\draw (-1.4599999999999982,5.019999999999992) node[anchor=north west] {$u_{8}$};
\draw (2.460000000000003,9.419999999999987) node[anchor=north west] {$y_{1}$};
\draw (5.000000000000004,9.359999999999987) node[anchor=north west] {$y_{2}$};
\draw (2.500000000000003,6.09999999999999) node[anchor=north west] {$y_{4}$};
\draw (4.780000000000004,6.15999999999999) node[anchor=north west] {$y_{3}$};
\draw (2.660000000000001,8.620000000000001)-- (-1.3399999999999999,9.863574311003802);
\draw (2.660000000000001,8.620000000000001)-- (0.03178715550190114,6.551787155501901);
\draw (0.03178715550190114,6.551787155501901)-- (-1.3399999999999999,9.863574311003802);
\draw [shift={(2.5288458004173706,6.819474528448952)}] plot[domain=0.6507309768653503:2.551484016792735,variable=\t]({1.0*3.0053171877104443*cos(\t r)+-0.0*3.0053171877104443*sin(\t r)},{0.0*3.0053171877104443*cos(\t r)+1.0*3.0053171877104443*sin(\t r)});
\draw (-3.2799999999999994,9.863574311003802)-- (-1.3399999999999999,9.863574311003802);
\draw (0.03178715550190114,8.4917871555019)-- (-1.3399999999999999,9.863574311003802);
\draw (0.03178715550190114,8.4917871555019)-- (0.03178715550190114,6.551787155501901);
\draw (0.03178715550190114,6.551787155501901)-- (-1.34,5.18);
\draw (-1.34,5.18)-- (-3.28,5.18);
\draw (-3.28,5.18)-- (-4.651787155501902,6.551787155501902);
\draw (-4.651787155501902,6.551787155501902)-- (-4.651787155501902,8.491787155501902);
\draw (-4.651787155501902,8.491787155501902)-- (-3.2799999999999994,9.863574311003802);
\draw (-3.2799999999999994,9.863574311003802)-- (-3.28,5.18);
\draw (-4.651787155501902,8.491787155501902)-- (-1.34,5.18);
\draw (-4.651787155501902,6.551787155501902)-- (0.03178715550190114,8.4917871555019);
\draw (4.920000000000002,8.639999999999999)-- (-4.651787155501902,6.551787155501902);
\draw (2.68,6.36)-- (-3.2799999999999994,9.863574311003802);
\draw (2.68,6.36)-- (-3.28,5.18);
\draw (4.94,6.38)-- (-4.651787155501902,8.491787155501902);
\draw [shift={(9.615720524017439,-35.122270742357934)}] plot[domain=1.6829850641143353:1.8362212174520522,variable=\t]({1.0*41.76482777638018*cos(\t r)+-0.0*41.76482777638018*sin(\t r)},{0.0*41.76482777638018*cos(\t r)+1.0*41.76482777638018*sin(\t r)});
\begin{scriptsize}
\draw [fill=black] (1.98,3.5) circle (1.5pt);
\draw [fill=black] (3.38,3.5) circle (1.5pt);
\draw [fill=black] (-3.28,5.18) circle (1.5pt);
\draw [fill=black] (-1.34,5.18) circle (1.5pt);
\draw [fill=black] (0.03178715550190114,6.551787155501901) circle (1.5pt);
\draw [fill=black] (0.03178715550190114,8.4917871555019) circle (1.5pt);
\draw [fill=black] (-1.3399999999999999,9.863574311003802) circle (1.5pt);
\draw [fill=black] (-3.2799999999999994,9.863574311003802) circle (1.5pt);
\draw [fill=black] (-4.651787155501902,8.491787155501902) circle (1.5pt);
\draw [fill=black] (-4.651787155501902,6.551787155501902) circle (1.5pt);
\draw [fill=black] (4.94,6.38) circle (1.5pt);
\draw [fill=black] (2.68,6.36) circle (1.5pt);
\draw [fill=black] (4.920000000000002,8.639999999999999) circle (1.5pt);
\draw [fill=black] (2.660000000000001,8.620000000000001) circle (1.5pt);
\end{scriptsize}
\end{tikzpicture}

%% file: graph3.tex
\begin{tikzpicture}[line cap=round,line join=round,>=triangle 45,x=1.0cm,y=1.0cm]
\clip(-5.76,4.219999999999994) rectangle (5.5400000000000045,10.519999999999992);
\draw (0.0800000000000023,6.499999999999993) node[anchor=north west] {$u_{1}$};
\draw (-4.02,10.439999999999992) node[anchor=north west] {$u_{4}$};
\draw (-4.819999999999999,1.819999999999995) node[anchor=north west] {$v_{1,3}$};
\draw (0.28000000000000236,8.879999999999992) node[anchor=north west] {$u_{2}$};
\draw (-1.019999999999998,10.379999999999992) node[anchor=north west] {$u_{3}$};
\draw (-0.639999999999998,1.639999999999995) node[anchor=north west] {$v_{ 2,3}$};
\draw (1.98,3.5)-- (3.38,3.5);
\draw (2.2737367544323214E-15,0.23999999999999538) node[anchor=north west] {$v_{2,1}$};
\draw (-5.34,9.139999999999992) node[anchor=north west] {$u_{5}$};
\draw (-5.38,6.619999999999993) node[anchor=north west] {$u_{6}$};
\draw (-3.959999999999999,5.059999999999993) node[anchor=north west] {$u_{7}$};
\draw (-1.4599999999999982,5.019999999999994) node[anchor=north west] {$u_{8}$};
\draw (2.440000000000003,9.299999999999992) node[anchor=north west] {$y_{1}$};
\draw (5.000000000000004,9.259999999999993) node[anchor=north west] {$y_{2}$};
\draw (2.500000000000003,6.099999999999993) node[anchor=north west] {$y_{4}$};
\draw (4.780000000000004,6.159999999999993) node[anchor=north west] {$y_{3}$};
\draw (2.660000000000001,8.620000000000001)-- (0.03178715550190114,6.551787155501901);
\draw (-3.2799999999999994,9.863574311003802)-- (-1.3399999999999999,9.863574311003802);
\draw (0.03178715550190114,8.4917871555019)-- (-1.3399999999999999,9.863574311003802);
\draw (0.03178715550190114,8.4917871555019)-- (0.03178715550190114,6.551787155501901);
\draw (0.03178715550190114,6.551787155501901)-- (-1.34,5.18);
\draw (-1.34,5.18)-- (-3.28,5.18);
\draw (-3.28,5.18)-- (-4.651787155501902,6.551787155501902);
\draw (-4.651787155501902,6.551787155501902)-- (-4.651787155501902,8.491787155501902);
\draw (-4.651787155501902,8.491787155501902)-- (-3.2799999999999994,9.863574311003802);
\draw (2.660000000000001,8.620000000000001)-- (-3.2799999999999994,9.863574311003802);
\draw (-3.2799999999999994,9.863574311003802)-- (0.03178715550190114,6.551787155501901);
\draw (-1.3399999999999999,9.863574311003802)-- (-4.651787155501902,6.551787155501902);
\draw [shift={(-2.578432999054916,-13.097808065771739)}] plot[domain=1.238628440765869:1.5169130636384296,variable=\t]({1.0*22.99475585749165*cos(\t r)+-0.0*22.99475585749165*sin(\t r)},{0.0*22.99475585749165*cos(\t r)+1.0*22.99475585749165*sin(\t r)});
\draw (4.920000000000002,8.639999999999999)-- (-4.651787155501902,6.551787155501902);
\draw (4.94,6.38)-- (0.03178715550190114,8.4917871555019);
\draw (4.94,6.38)-- (-3.28,5.18);
\draw (0.03178715550190114,8.4917871555019)-- (-3.28,5.18);
\begin{scriptsize}
\draw [fill=black] (1.98,3.5) circle (1.5pt);
\draw [fill=black] (3.38,3.5) circle (1.5pt);
\draw [fill=black] (-3.28,5.18) circle (1.5pt);
\draw [fill=black] (-1.34,5.18) circle (1.5pt);
\draw [fill=black] (0.03178715550190114,6.551787155501901) circle (1.5pt);
\draw [fill=black] (0.03178715550190114,8.4917871555019) circle (1.5pt);
\draw [fill=black] (-1.3399999999999999,9.863574311003802) circle (1.5pt);
\draw [fill=black] (-3.2799999999999994,9.863574311003802) circle (1.5pt);
\draw [fill=black] (-4.651787155501902,8.491787155501902) circle (1.5pt);
\draw [fill=black] (-4.651787155501902,6.551787155501902) circle (1.5pt);
\draw [fill=black] (4.94,6.38) circle (1.5pt);
\draw [fill=black] (2.68,6.36) circle (1.5pt);
\draw [fill=black] (4.920000000000002,8.639999999999999) circle (1.5pt);
\draw [fill=black] (2.660000000000001,8.620000000000001) circle (1.5pt);
\end{scriptsize}
\end{tikzpicture}

%% file: graph4.tex
\definecolor{uuuuuu}{rgb}{0.26666666666666666,0.26666666666666666,0.26666666666666666}
\begin{tikzpicture}[line cap=round,line join=round,>=triangle 45,x=1.0cm,y=1.0cm]
\clip(-4.779999999999998,4.879999999999996) rectangle (5.760000000000007,10.27999999999999);
\draw (-4.819999999999998,1.82) node[anchor=north west] {$v_{1,3}$};
\draw (-0.6399999999999957,1.6400000000000003) node[anchor=north west] {$v_{ 2,3}$};
\draw (4.547473508864643E-15,0.240000000000002) node[anchor=north west] {$v_{2,1}$};
\draw (2.5600000000000054,9.339999999999991) node[anchor=north west] {$y_{1}$};
\draw (4.7800000000000065,9.379999999999992) node[anchor=north west] {$y_{2}$};
\draw (2.4800000000000053,6.119999999999995) node[anchor=north west] {$y_{4}$};
\draw (4.7800000000000065,6.079999999999995) node[anchor=north west] {$y_{3}$};
\draw (-2.3199999999999963,10.33999999999999) node[anchor=north west] {$u_{4}$};
\draw (-0.25999999999999557,9.45999999999999) node[anchor=north west] {$u_{3}$};
\draw (-4.359999999999997,9.319999999999991) node[anchor=north west] {$u_{5}$};
\draw (0.1600000000000046,7.779999999999993) node[anchor=north west] {$u_{2}$};
\draw (4.547473508864643E-15,0.240000000000002) node[anchor=north west] {$u_{4}$};
\draw (-4.719999999999997,7.299999999999994) node[anchor=north west] {$u_{6}$};
\draw (-1.0999999999999959,5.519999999999996) node[anchor=north west] {$u_{1}$};
\draw (-2.86,5.76)-- (-1.1,5.76);
\draw (-1.1,5.76)-- (-0.0026579487286291226,7.136023409143732);
\draw (-0.0026579487286291226,7.136023409143732)-- (-0.3942947924917424,8.85189653458374);
\draw (-0.3942947924917424,8.85189653458374)-- (-1.9799999999999998,9.615531915430644);
\draw (-1.9799999999999998,9.615531915430644)-- (-3.565705207508257,8.851896534583743);
\draw (-3.565705207508257,8.851896534583743)-- (-3.957342051271371,7.136023409143733);
\draw (-3.957342051271371,7.136023409143733)-- (-2.86,5.76);
\draw (-3.2599999999999967,5.579999999999996) node[anchor=north west] {$u_{7}$};
\draw (-1.1,5.76)-- (2.660000000000001,8.620000000000001);
\draw [shift={(-0.40445867448399253,5.64797438592881)}] plot[domain=0.7700869602593796:1.9488054819323313,variable=\t]({1.0*4.268939355192978*cos(\t r)+-0.0*4.268939355192978*sin(\t r)},{0.0*4.268939355192978*cos(\t r)+1.0*4.268939355192978*sin(\t r)});
\draw [shift={(0.3041915496703065,14.939827259556957)}] plot[domain=4.673088722381368:5.344718096586081,variable=\t]({1.0*7.809834258829822*cos(\t r)+-0.0*7.809834258829822*sin(\t r)},{0.0*7.809834258829822*cos(\t r)+1.0*7.809834258829822*sin(\t r)});
\draw [shift={(58.245260581215845,-333.0682120943088)}] plot[domain=1.7256026690398594:1.7516378019493288,variable=\t]({1.0*345.8440192178313*cos(\t r)+-0.0*345.8440192178313*sin(\t r)},{0.0*345.8440192178313*cos(\t r)+1.0*345.8440192178313*sin(\t r)});
\draw (-1.1,5.76)-- (-1.9799999999999998,9.615531915430644);
\draw (-0.0026579487286291226,7.136023409143732)-- (-3.957342051271371,7.136023409143733);
\begin{scriptsize}
\draw [fill=black] (4.94,6.38) circle (1.5pt);
\draw [fill=black] (2.68,6.36) circle (1.5pt);
\draw [fill=black] (4.920000000000002,8.639999999999999) circle (1.5pt);
\draw [fill=black] (2.660000000000001,8.620000000000001) circle (1.5pt);
\draw [fill=black] (-2.86,5.76) circle (1.5pt);
\draw [fill=black] (-1.1,5.76) circle (1.5pt);
\draw [fill=uuuuuu] (-0.0026579487286291226,7.136023409143732) circle (1.5pt);
\draw [fill=uuuuuu] (-0.3942947924917424,8.85189653458374) circle (1.5pt);
\draw [fill=uuuuuu] (-1.9799999999999998,9.615531915430644) circle (1.5pt);
\draw [fill=uuuuuu] (-3.565705207508257,8.851896534583743) circle (1.5pt);
\draw [fill=uuuuuu] (-3.957342051271371,7.136023409143733) circle (1.5pt);
\end{scriptsize}
\end{tikzpicture}

%% file: graph5.tex
\begin{tikzpicture}[line cap=round,line join=round,>=triangle 45,x=1.0cm,y=1.0cm]
\clip(-2.459999999999998,5.439999999999992) rectangle (10.180000000000016,9.379999999999988);
\draw (-4.82,1.8199999999999958) node[anchor=north west] {$v_{1,3}$};
\draw (-0.6399999999999961,1.6399999999999961) node[anchor=north west] {$v_{ 2,3}$};
\draw (4.547473508864646E-15,0.2399999999999975) node[anchor=north west] {$v_{2,1}$};
\draw (7.280000000000013,6.099999999999992) node[anchor=north west] {$y_{1}$};
\draw (0.3200000000000049,9.479999999999988) node[anchor=north west] {$y_{2}$};
\draw (9.140000000000015,9.459999999999988) node[anchor=north west] {$y_{4}$};
\draw (-2.039999999999998,6.139999999999992) node[anchor=north west] {$y_{3}$};
\draw (2.660000000000001,8.620000000000001)-- (4.920000000000002,8.639999999999999);
\draw (4.920000000000002,8.639999999999999)-- (4.94,6.38);
\draw (4.94,6.38)-- (2.68,6.36);
\draw (2.68,6.36)-- (2.660000000000001,8.620000000000001);
\draw (2.380000000000007,9.499999999999988) node[anchor=north west] {$u_{4}$};
\draw (4.94000000000001,9.479999999999988) node[anchor=north west] {$u_{3}$};
\draw (0.40000000000000496,6.139999999999992) node[anchor=north west] {$u_{5}$};
\draw (6.940000000000012,9.479999999999988) node[anchor=north west] {$u_{2}$};
\draw (4.547473508864646E-15,0.2399999999999975) node[anchor=north west] {$u_{4}$};
\draw (2.380000000000007,6.099999999999992) node[anchor=north west] {$u_{6}$};
\draw (4.92000000000001,6.039999999999991) node[anchor=north west] {$u_{1}$};
\draw (2.660000000000001,8.620000000000001)-- (0.62,6.36);
\draw (0.62,6.36)-- (2.68,6.36);
\draw (2.660000000000001,8.620000000000001)-- (0.52,8.62);
\draw (0.52,8.62)-- (2.68,6.36);
\draw (4.920000000000002,8.639999999999999)-- (7.04,8.68);
\draw (7.04,8.68)-- (4.94,6.38);
\draw (4.920000000000002,8.639999999999999)-- (7.22,6.38);
\draw (7.22,6.38)-- (4.94,6.38);
\draw (7.04,8.68)-- (9.26,8.66);
\draw (0.62,6.36)-- (-1.56,6.34);
\draw [dash pattern=on 1pt off 3pt on 5pt off 4pt] (2.68,6.36)-- (4.920000000000002,8.639999999999999);
\begin{scriptsize}
\draw [fill=black] (4.94,6.38) circle (1.5pt);
\draw [fill=black] (2.68,6.36) circle (1.5pt);
\draw [fill=black] (4.920000000000002,8.639999999999999) circle (1.5pt);
\draw [fill=black] (2.660000000000001,8.620000000000001) circle (1.5pt);
\draw [fill=black] (0.62,6.36) circle (1.5pt);
\draw [fill=black] (0.52,8.62) circle (1.5pt);
\draw [fill=black] (7.04,8.68) circle (1.5pt);
\draw [fill=black] (7.22,6.38) circle (1.5pt);
\draw [fill=black] (9.26,8.66) circle (1.5pt);
\draw [fill=black] (-1.56,6.34) circle (1.5pt);
\end{scriptsize}
\end{tikzpicture}

%% file: graph6.tex
\definecolor{uuuuuu}{rgb}{0.26666666666666666,0.26666666666666666,0.26666666666666666}
\begin{tikzpicture}[line cap=round,line join=round,>=triangle 45,x=1.0cm,y=1.0cm]
\clip(-4.479999999999998,5.159999999999992) rectangle (5.980000000000007,10.239999999999988);
\draw (-4.819999999999998,1.819999999999995) node[anchor=north west] {$v_{1,3}$};
\draw (-0.6399999999999962,1.6399999999999952) node[anchor=north west] {$v_{ 2,3}$};
\draw (3.979039320256562E-15,0.23999999999999638) node[anchor=north west] {$v_{2,1}$};
\draw (2.440000000000005,9.299999999999988) node[anchor=north west] {$y_{1}$};
\draw (5.000000000000006,9.25999999999999) node[anchor=north west] {$y_{2}$};
\draw (2.480000000000005,6.119999999999991) node[anchor=north west] {$y_{4}$};
\draw (4.780000000000006,6.159999999999991) node[anchor=north west] {$y_{3}$};
\draw (-2.8599999999999994,9.068409421321224)-- (-1.0999999999999999,9.068409421321224);
\draw (-3.74,7.544204710660613)-- (-2.8599999999999994,9.068409421321224);
\draw (-1.0999999999999999,9.068409421321224)-- (-0.21999999999999997,7.544204710660612);
\draw (-0.21999999999999997,7.544204710660612)-- (-1.1,6.02);
\draw (-1.1,6.02)-- (-2.86,6.02);
\draw (-2.86,6.02)-- (-3.74,7.544204710660613);
\draw (-1.1,6.02)-- (2.660000000000001,8.620000000000001);
\draw (4.920000000000002,8.639999999999999)-- (-0.21999999999999997,7.544204710660612);
\draw (4.920000000000002,8.639999999999999)-- (-3.74,7.544204710660613);
\draw (-3.4799999999999973,9.679999999999989) node[anchor=north west] {$u_{4}$};
\draw (-0.9599999999999964,9.679999999999989) node[anchor=north west] {$u_{3}$};
\draw (-4.359999999999998,7.85999999999999) node[anchor=north west] {$u_{5}$};
\draw (-0.09999999999999606,7.63999999999999) node[anchor=north west] {$u_{2}$};
\draw (3.979039320256562E-15,0.23999999999999638) node[anchor=north west] {$u_{4}$};
\draw (-3.079999999999997,5.959999999999992) node[anchor=north west] {$u_{6}$};
\draw (-1.2599999999999965,5.919999999999992) node[anchor=north west] {$u_{1}$};
\draw [shift={(-0.32375965080700764,6.089683476884233)}] plot[domain=0.7033490365581737:2.276132395114745,variable=\t]({1.0*3.9122018558526133*cos(\t r)+-0.0*3.9122018558526133*sin(\t r)},{0.0*3.9122018558526133*cos(\t r)+1.0*3.9122018558526133*sin(\t r)});
\draw (-2.8599999999999994,9.068409421321224)-- (-1.1,6.02);
\begin{scriptsize}
\draw [fill=black] (4.94,6.38) circle (1.5pt);
\draw [fill=black] (2.68,6.36) circle (1.5pt);
\draw [fill=black] (4.920000000000002,8.639999999999999) circle (1.5pt);
\draw [fill=black] (2.660000000000001,8.620000000000001) circle (1.5pt);
\draw [fill=black] (-2.86,6.02) circle (1.5pt);
\draw [fill=black] (-1.1,6.02) circle (1.5pt);
\draw [fill=uuuuuu] (-0.21999999999999997,7.544204710660612) circle (1.5pt);
\draw [fill=uuuuuu] (-1.0999999999999999,9.068409421321224) circle (1.5pt);
\draw [fill=uuuuuu] (-2.8599999999999994,9.068409421321224) circle (1.5pt);
\draw [fill=uuuuuu] (-3.74,7.544204710660613) circle (1.5pt);
\end{scriptsize}
\end{tikzpicture}

%% file: grapht71.tex
\definecolor{qqqqff}{rgb}{0.0,0.0,1.0}
\begin{tikzpicture}[line cap=round,line join=round,>=triangle 45,x=1.0cm,y=1.0cm]
\clip(8.51999999999986,-0.019999999999995428) rectangle (17.059999999999725,7.239999999999996);
\draw (13.739999999999851,5.78)-- (10.999999999999881,5.78);
\draw (10.999999999999881,5.78)-- (9.629999999999896,3.4070903936306642);
\draw (9.629999999999896,3.4070903936306642)-- (10.999999999999881,1.0341807872613273);
\draw (10.999999999999881,1.0341807872613273)-- (13.739999999999851,1.034180787261327);
\draw (13.739999999999851,1.034180787261327)-- (15.109999999999836,3.407090393630662);
\draw (15.109999999999836,3.407090393630662)-- (13.739999999999851,5.78);
\draw (13.739999999999851,5.78)-- (12.419999999999865,4.96);
\draw (12.419999999999865,4.96)-- (10.999999999999881,5.78);
\draw (9.629999999999896,3.4070903936306642)-- (10.979999999999881,2.56);
\draw (10.979999999999881,2.56)-- (10.999999999999881,1.0341807872613273);
\draw (13.739999999999851,1.034180787261327)-- (13.559999999999851,2.74);
\draw (13.559999999999851,2.74)-- (15.109999999999836,3.407090393630662);
\draw (12.419999999999865,3.58)-- (12.439999999999872,4.12);
\draw (12.419999999999865,3.58)-- (11.999999999999876,3.22);
\draw (12.419999999999865,3.58)-- (12.859999999999866,3.16);
\draw (12.41999999999987,4.14)-- (12.39999999999987,4.98);
\draw [shift={(15.461162790697863,2.6844258720927128)}] plot[domain=2.0782402072298716:2.989798983171639,variable=\t]({1.0*3.541889401640361*cos(\t r)+-0.0*3.541889401640361*sin(\t r)},{0.0*3.541889401640361*cos(\t r)+1.0*3.541889401640361*sin(\t r)});
\draw [shift={(13.296725521669158,6.8030818619583835)}] plot[domain=4.3553150704901205:5.197492845947945,variable=\t]({1.0*3.824305263675662*cos(\t r)+-0.0*3.824305263675662*sin(\t r)},{0.0*3.824305263675662*cos(\t r)+1.0*3.824305263675662*sin(\t r)});
\draw (12.839999999999867,3.16)-- (13.559999999999858,2.74);
\draw (11.959999999999798,3.22)-- (10.959999999999887,2.58);
\draw [shift={(14.814385406366979,1.1914937279122126)}] plot[domain=2.252852163613375:3.186758623824585,variable=\t]({1.0*3.7982589051781086*cos(\t r)+-0.0*3.7982589051781086*sin(\t r)},{0.0*3.7982589051781086*cos(\t r)+1.0*3.7982589051781086*sin(\t r)});
\draw [shift={(7.678884147393201,0.3636921132024486)}] plot[domain=0.11513289656071614:0.6726043249119994,variable=\t]({1.0*6.0612441613688235*cos(\t r)+-0.0*6.0612441613688235*sin(\t r)},{0.0*6.0612441613688235*cos(\t r)+1.0*6.0612441613688235*sin(\t r)});
\draw [shift={(9.672175439901563,2.869231849776593)}] plot[domain=0.09153148197352397:1.1402084219055215,variable=\t]({1.0*3.1811410784098135*cos(\t r)+-0.0*3.1811410784098135*sin(\t r)},{0.0*3.1811410784098135*cos(\t r)+1.0*3.1811410784098135*sin(\t r)});
\draw [shift={(11.473431546959734,6.163003654889765)}] plot[domain=4.123185344936516:5.139448790742529,variable=\t]({1.0*3.2993242168852985*cos(\t r)+-0.0*3.2993242168852985*sin(\t r)},{0.0*3.2993242168852985*cos(\t r)+1.0*3.2993242168852985*sin(\t r)});
\draw (12.719999999999793,4.159999999999999) node[anchor=north west] {$x$};
\draw (11.999999999999805,5.679999999999998) node[anchor=north west] {$y_1$};
\draw (10.27999999999983,2.760000000000001) node[anchor=north west] {$y_2$};
\draw (13.739999999999776,2.820000000000001) node[anchor=north west] {$y_3$};
\begin{scriptsize}
\draw [fill=black] (13.739999999999851,5.78) circle (1.5pt);
\draw [fill=black] (10.999999999999881,5.78) circle (1.5pt);
\draw [fill=black] (9.629999999999896,3.4070903936306642) circle (1.5pt);
\draw [fill=black] (10.999999999999881,1.0341807872613273) circle (1.5pt);
\draw [fill=black] (13.739999999999851,1.034180787261327) circle (1.5pt);
\draw [fill=black] (15.109999999999836,3.407090393630662) circle (1.5pt);
\draw [fill=qqqqff] (12.419999999999865,4.96) circle (1.5pt);
\draw [fill=qqqqff] (10.979999999999881,2.56) circle (1.5pt);
\draw [fill=qqqqff] (13.559999999999851,2.74) circle (1.5pt);
\draw [fill=qqqqff] (12.419999999999865,3.58) circle (1.5pt);
\draw [fill=qqqqff] (12.439999999999872,4.12) circle (1.5pt);
\draw [fill=qqqqff] (11.999999999999876,3.22) circle (1.5pt);
\draw [fill=qqqqff] (12.859999999999866,3.16) circle (1.5pt);
\draw [fill=black] (13.739999999999856,5.78) circle (2.5pt);
\draw [fill=black] (10.999999999999881,5.76) circle (2.5pt);
\draw [fill=black] (12.41999999999987,3.56) circle (2.5pt);
\draw [fill=black] (12.39999999999987,4.98) circle (2.5pt);
\draw [fill=black] (12.41999999999987,4.14) circle (2.5pt);
\draw [fill=black] (11.959999999999798,3.22) circle (2.5pt);
\draw [fill=black] (12.839999999999867,3.16) circle (2.5pt);
\draw [fill=black] (9.64,3.42) circle (2.5pt);
\draw [fill=black] (15.079999999999844,3.42) circle (2.5pt);
\draw [fill=black] (13.559999999999858,2.74) circle (2.5pt);
\draw [fill=black] (11.019999999999886,1.02) circle (2.5pt);
\draw [fill=black] (13.699999999999857,1.06) circle (2.5pt);
\draw [fill=black] (10.959999999999887,2.58) circle (2.5pt);
\end{scriptsize}
\end{tikzpicture}

%% file: grapht72.tex
\definecolor{qqqqff}{rgb}{0.0,0.0,1.0}
\begin{tikzpicture}[line cap=round,line join=round,>=triangle 45,x=1.0cm,y=1.0cm]
\clip(8.53999999999987,0.5800000000000038) rectangle (17.039999999999747,6.359999999999997);
\draw (13.739999999999851,5.78)-- (10.999999999999881,5.78);
\draw (10.999999999999881,5.78)-- (9.629999999999896,3.4070903936306642);
\draw (9.629999999999896,3.4070903936306642)-- (10.999999999999881,1.0341807872613273);
\draw (10.999999999999881,1.0341807872613273)-- (13.739999999999851,1.034180787261327);
\draw (13.739999999999851,1.034180787261327)-- (15.109999999999836,3.407090393630662);
\draw (15.109999999999836,3.407090393630662)-- (13.739999999999851,5.78);
\draw (13.739999999999851,5.78)-- (12.419999999999865,4.96);
\draw (12.419999999999865,4.96)-- (10.999999999999881,5.78);
\draw (9.629999999999896,3.4070903936306642)-- (10.979999999999881,2.56);
\draw (10.979999999999881,2.56)-- (10.999999999999881,1.0341807872613273);
\draw (13.739999999999851,1.034180787261327)-- (13.559999999999851,2.74);
\draw (13.559999999999851,2.74)-- (15.109999999999836,3.407090393630662);
\draw (12.419999999999865,3.58)-- (12.439999999999872,4.12);
\draw (12.419999999999865,3.58)-- (11.999999999999876,3.22);
\draw (12.419999999999865,3.58)-- (12.859999999999866,3.16);
\draw (12.41999999999987,4.14)-- (12.39999999999987,4.98);
\draw [shift={(15.461162790697863,2.6844258720927128)}] plot[domain=2.0782402072298716:2.989798983171639,variable=\t]({1.0*3.541889401640361*cos(\t r)+-0.0*3.541889401640361*sin(\t r)},{0.0*3.541889401640361*cos(\t r)+1.0*3.541889401640361*sin(\t r)});
\draw [shift={(13.296725521669158,6.8030818619583835)}] plot[domain=4.3553150704901205:5.197492845947945,variable=\t]({1.0*3.824305263675662*cos(\t r)+-0.0*3.824305263675662*sin(\t r)},{0.0*3.824305263675662*cos(\t r)+1.0*3.824305263675662*sin(\t r)});
\draw (12.839999999999867,3.16)-- (13.559999999999858,2.74);
\draw (11.959999999999798,3.22)-- (10.959999999999887,2.58);
\draw [shift={(14.814385406366979,1.1914937279122126)}] plot[domain=2.252852163613375:3.186758623824585,variable=\t]({1.0*3.7982589051781086*cos(\t r)+-0.0*3.7982589051781086*sin(\t r)},{0.0*3.7982589051781086*cos(\t r)+1.0*3.7982589051781086*sin(\t r)});
\draw [shift={(7.678884147393201,0.3636921132024486)}] plot[domain=0.11513289656071614:0.6726043249119994,variable=\t]({1.0*6.0612441613688235*cos(\t r)+-0.0*6.0612441613688235*sin(\t r)},{0.0*6.0612441613688235*cos(\t r)+1.0*6.0612441613688235*sin(\t r)});
\draw [shift={(9.672175439901563,2.869231849776593)}] plot[domain=0.09153148197352397:1.1402084219055215,variable=\t]({1.0*3.1811410784098135*cos(\t r)+-0.0*3.1811410784098135*sin(\t r)},{0.0*3.1811410784098135*cos(\t r)+1.0*3.1811410784098135*sin(\t r)});
\draw [shift={(11.473431546959734,6.163003654889765)}] plot[domain=4.123185344936516:5.139448790742529,variable=\t]({1.0*3.2993242168852985*cos(\t r)+-0.0*3.2993242168852985*sin(\t r)},{0.0*3.2993242168852985*cos(\t r)+1.0*3.2993242168852985*sin(\t r)});
\draw (12.719999999999809,4.159999999999999) node[anchor=north west] {$x$};
\draw (15.839999999999764,5.619999999999997) node[anchor=north west] {$y$};
\draw (11.999999999999819,5.679999999999998) node[anchor=north west] {$y_1$};
\draw (10.279999999999843,2.760000000000001) node[anchor=north west] {$y_2$};
\draw (13.739999999999794,2.820000000000001) node[anchor=north west] {$y_3$};
\draw (15.619999999999838,5.18)-- (12.39999999999987,4.98);
\draw (15.619999999999838,5.18)-- (13.559999999999858,2.74);
\draw [shift={(10.637875150059784,8.633423769507935)}] plot[domain=4.7655524998660574:5.67706280057867,variable=\t]({1.0*6.0619884322054425*cos(\t r)+-0.0*6.0619884322054425*sin(\t r)},{0.0*6.0619884322054425*cos(\t r)+1.0*6.0619884322054425*sin(\t r)});
\begin{scriptsize}
\draw [fill=black] (13.739999999999851,5.78) circle (1.5pt);
\draw [fill=black] (10.999999999999881,5.78) circle (1.5pt);
\draw [fill=black] (9.629999999999896,3.4070903936306642) circle (1.5pt);
\draw [fill=black] (10.999999999999881,1.0341807872613273) circle (1.5pt);
\draw [fill=black] (13.739999999999851,1.034180787261327) circle (1.5pt);
\draw [fill=black] (15.109999999999836,3.407090393630662) circle (1.5pt);
\draw [fill=qqqqff] (12.419999999999865,4.96) circle (1.5pt);
\draw [fill=qqqqff] (10.979999999999881,2.56) circle (1.5pt);
\draw [fill=qqqqff] (13.559999999999851,2.74) circle (1.5pt);
\draw [fill=qqqqff] (12.419999999999865,3.58) circle (1.5pt);
\draw [fill=qqqqff] (12.439999999999872,4.12) circle (1.5pt);
\draw [fill=qqqqff] (11.999999999999876,3.22) circle (1.5pt);
\draw [fill=qqqqff] (12.859999999999866,3.16) circle (1.5pt);
\draw [fill=black] (13.739999999999856,5.78) circle (2.5pt);
\draw [fill=black] (10.999999999999881,5.76) circle (2.5pt);
\draw [fill=black] (12.41999999999987,3.56) circle (2.5pt);
\draw [fill=black] (12.39999999999987,4.98) circle (2.5pt);
\draw [fill=black] (12.41999999999987,4.14) circle (2.5pt);
\draw [fill=black] (11.959999999999798,3.22) circle (2.5pt);
\draw [fill=black] (12.839999999999867,3.16) circle (2.5pt);
\draw [fill=black] (9.64,3.42) circle (2.5pt);
\draw [fill=black] (15.079999999999844,3.42) circle (2.5pt);
\draw [fill=black] (13.559999999999858,2.74) circle (2.5pt);
\draw [fill=black] (11.019999999999886,1.02) circle (2.5pt);
\draw [fill=black] (13.699999999999857,1.06) circle (2.5pt);
\draw [fill=black] (15.619999999999838,5.18) circle (2.5pt);
\draw [fill=black] (10.959999999999887,2.58) circle (2.5pt);
\end{scriptsize}
\end{tikzpicture}

%% file: graph7.tex
\begin{tikzpicture}[line cap=round,line join=round,>=triangle 45,x=0.8787878787878788cm,y=1.0cm]
\clip(-5.260000000000001,-0.5599999999999986) rectangle (11.240000000000004,6.519999999999998);
\draw [rotate around={0.0:(-3.129999999999998,2.3)}] (-3.129999999999998,2.3) ellipse (1.364478596925364cm and 0.6918981674693683cm);
\draw [rotate around={0.0:(0.8500000000000001,2.24)}] (0.8500000000000001,2.24) ellipse (1.3644785969253699cm and 0.6918981674693704cm);
\draw [rotate around={0.0:(9.089999999999998,2.24)}] (9.089999999999998,2.24) ellipse (1.3644785969253448cm and 0.6918981674693577cm);
\draw [line width=2.0pt,dash pattern=on 5pt off 5pt] (1.38,5.5)-- (-1.74,3.42);
\draw (1.98,5.58)-- (7.6,3.1);
\draw (-4.98,1.5800000000000003) node[anchor=north west] {$x_1, x_2, x_3,$ $  ...  $ $ $ $     ,x_{n-1}$};
\draw (0.5000000000000013,1.5400000000000003) node[anchor=north west] {$x_{n}, ...$};
\draw (4.320000000000003,1.5000000000000004) node[anchor=north west] {$x_{2n-1}, ... $};
\draw (2.740000000000002,5.739999999999998) node[anchor=north west] {red degree 1};
\draw (2.720000000000002,6.299999999999998) node[anchor=north west] {blue degree};
\draw (5.260000000000003,6.339999999999998) node[anchor=north west] {$3(n-1)$};
\draw (-3.3600000000000003,3.7999999999999994) node[anchor=north west] {$K_{n-1}$};
\draw (1.2800000000000016,3.7599999999999993) node[anchor=north west] {$K_{n-1}$};
\draw (8.420000000000003,3.7399999999999993) node[anchor=north west] {$K_{n-1}$};
\draw (-4.12,3.7999999999999994) node[anchor=north west] {red};
\draw (0.5400000000000014,3.7999999999999994) node[anchor=north west] {red};
\draw (7.680000000000004,3.7399999999999993) node[anchor=north west] {red};
\draw [line width=2.4000000000000004pt,dash pattern=on 5pt off 5pt] (-0.86,2.26)-- (-1.44,2.26);
\draw [line width=2.0pt,dash pattern=on 5pt off 5pt] (3.24,2.26)-- (2.58,2.26);
\draw [line width=2.0pt,dash pattern=on 5pt off 5pt] (1.86,5.44)-- (3.48,3.3);
\draw [rotate around={0.0:(4.929999999999998,2.2600000000000002)}] (4.929999999999998,2.2600000000000002) ellipse (1.364478596925369cm and 0.6918981674693704cm);
\draw (4.900000000000003,3.7399999999999993) node[anchor=north west] {$K_{n-1}$};
\draw (4.160000000000003,3.7599999999999993) node[anchor=north west] {red};
\draw [line width=2.0pt,dash pattern=on 5pt off 5pt] (7.3,2.24)-- (6.64,2.24);
\draw (7.840000000000004,1.5200000000000005) node[anchor=north west] {$x_{3n-2}, ... ,x_{4(n-1)}$};
\draw [line width=2.0pt,dash pattern=on 5pt off 5pt] (1.58,5.32)-- (1.12,3.86);
\draw [shift={(0.2820000000000002,9.632)},line width=1.6pt,dash pattern=on 5pt off 5pt]  plot[domain=4.32310744680871:5.112689686208233,variable=\t]({1.0*9.438143249601588*cos(\t r)+-0.0*9.438143249601588*sin(\t r)},{0.0*9.438143249601588*cos(\t r)+1.0*9.438143249601588*sin(\t r)});
\draw [shift={(4.753419020801163,8.634857843582521)},line width=2.0pt,dash pattern=on 5pt off 5pt]  plot[domain=4.305783209729312:5.112915742404995,variable=\t]({1.0*8.37792234628314*cos(\t r)+-0.0*8.37792234628314*sin(\t r)},{0.0*8.37792234628314*cos(\t r)+1.0*8.37792234628314*sin(\t r)});
\draw [shift={(2.3699999999999997,15.532456140350872)},line width=2.0pt,dash pattern=on 5pt off 5pt]  plot[domain=4.329002014748123:5.095775946021256,variable=\t]({1.0*15.799448437996368*cos(\t r)+-0.0*15.799448437996368*sin(\t r)},{0.0*15.799448437996368*cos(\t r)+1.0*15.799448437996368*sin(\t r)});
\begin{scriptsize}
\draw [fill=black] (-4.44,2.28) circle (2.0pt);
\draw [fill=black] (-4.0,2.28) circle (2.0pt);
\draw [fill=black] (-3.58,2.26) circle (2.0pt);
\draw [fill=black] (-1.96,2.28) circle (2.0pt);
\draw [fill=black] (-3.32,2.46) circle (0.5pt);
\draw [fill=black] (-3.1,2.46) circle (0.5pt);
\draw [fill=black] (-2.82,2.46) circle (0.5pt);
\draw [fill=black] (-0.1,2.22) circle (2.0pt);
\draw [fill=black] (0.34,2.22) circle (2.0pt);
\draw [fill=black] (0.76,2.2) circle (2.0pt);
\draw [fill=black] (1.02,2.2) circle (0.5pt);
\draw [fill=black] (1.24,2.2) circle (0.5pt);
\draw [fill=black] (1.52,2.2) circle (0.5pt);
\draw [fill=black] (8.26,2.24) circle (2.0pt);
\draw [fill=black] (9.3,2.24) circle (0.5pt);
\draw [fill=black] (9.62,2.28) circle (0.5pt);
\draw [fill=black] (9.02,2.24) circle (0.5pt);
\draw [fill=black] (3.98,2.24) circle (2.0pt);
\draw [fill=black] (4.42,2.24) circle (2.0pt);
\draw [fill=black] (4.84,2.22) circle (2.0pt);
\draw [fill=black] (5.1,2.22) circle (0.5pt);
\draw [fill=black] (5.32,2.22) circle (0.5pt);
\draw [fill=black] (5.6,2.22) circle (0.5pt);
\draw [fill=black] (8.56,2.24) circle (2.0pt);
\draw [fill=black] (7.92,2.24) circle (2.0pt);
\end{scriptsize}
\end{tikzpicture}

%% file: graph81.tex
\begin{tikzpicture}[line cap=round,line join=round,>=triangle 45,x=1.0cm,y=1.0cm]
\clip(-4.860000000000002,1.5400000000000027) rectangle (9.940000000000005,6.24);
\draw [rotate around={1.145762838175131:(-3.8999999999999995,5.429999999999998)}] (-3.8999999999999995,5.429999999999998) ellipse (0.7473599204960948cm and 0.5553799157009054cm);
\draw [rotate around={1.1457628381750795:(-1.48,5.450000000000001)}] (-1.48,5.450000000000001) ellipse (0.7473599204961263cm and 0.5553799157009286cm);
\draw [rotate around={1.145762838175105:(-3.9400000000000004,3.289999999999999)}] (-3.9400000000000004,3.289999999999999) ellipse (0.7473599204961175cm and 0.5553799157009224cm);
\draw [rotate around={1.145762838175105:(-1.4000000000000001,3.290000000000001)}] (-1.4000000000000001,3.290000000000001) ellipse (0.7473599204961234cm and 0.5553799157009267cm);
\draw (-1.98,5.44)-- (-3.4,5.44);
\draw (-3.4,5.44)-- (-3.44,3.3);
\draw [rotate around={1.1457628381751304:(5.96,5.470000000000001)}] (5.96,5.470000000000001) ellipse (0.7473599204961332cm and 0.5553799157009339cm);
\draw [rotate around={1.1457628381750788:(8.380000000000003,5.49)}] (8.380000000000003,5.49) ellipse (0.7473599204961118cm and 0.5553799157009173cm);
\draw [rotate around={1.1457628381751053:(5.92,3.3299999999999996)}] (5.92,3.3299999999999996) ellipse (0.7473599204961027cm and 0.5553799157009112cm);
\draw [rotate around={1.1457628381751028:(8.460000000000003,3.33)}] (8.460000000000003,3.33) ellipse (0.7473599204961426cm and 0.55537991570094cm);
\draw (7.88,5.48)-- (6.46,5.48);
\draw (6.46,5.48)-- (6.42,3.34);
\draw (7.96,3.32)-- (6.42,3.34);
\draw [rotate around={1.1457628381751308:(0.9600000000000003,5.429999999999999)}] (0.9600000000000003,5.429999999999999) ellipse (0.747359920496095cm and 0.5553799157009055cm);
\draw [rotate around={1.1457628381750795:(3.3799999999999994,5.450000000000001)}] (3.3799999999999994,5.450000000000001) ellipse (0.7473599204961185cm and 0.5553799157009228cm);
\draw [rotate around={1.145762838175105:(0.9199999999999995,3.29)}] (0.9199999999999995,3.29) ellipse (0.7473599204961217cm and 0.5553799157009254cm);
\draw [rotate around={1.145762838175105:(3.459999999999999,3.289999999999999)}] (3.459999999999999,3.289999999999999) ellipse (0.7473599204961178cm and 0.5553799157009225cm);
\draw (2.88,5.44)-- (1.46,5.44);
\draw (1.16,3.3)-- (1.2,5.44);
\draw (-3.4400000000000013,2.780000000000002) node[anchor=north west] {$R_{4n-4,3}$};
\draw (1.4200000000000013,2.800000000000002) node[anchor=north west] {$R_{4n-4,4}$};
\draw (6.440000000000003,2.800000000000002) node[anchor=north west] {$R_{4n-4,5}$};
\begin{scriptsize}
\draw [fill=black] (-4.18,5.42) circle (0.5pt);
\draw [fill=black] (-4.02,5.42) circle (0.5pt);
\draw [fill=black] (-3.9,5.42) circle (0.5pt);
\draw [fill=black] (-4.4,5.42) circle (1.5pt);
\draw [fill=black] (-3.4,5.44) circle (1.5pt);
\draw [fill=black] (-3.66,5.44) circle (1.5pt);
\draw [fill=black] (-1.76,5.44) circle (0.5pt);
\draw [fill=black] (-1.6,5.44) circle (0.5pt);
\draw [fill=black] (-1.48,5.44) circle (0.5pt);
\draw [fill=black] (-1.98,5.44) circle (1.5pt);
\draw [fill=black] (-0.98,5.46) circle (1.5pt);
\draw [fill=black] (-1.24,5.46) circle (1.5pt);
\draw [fill=black] (-4.22,3.28) circle (0.5pt);
\draw [fill=black] (-4.06,3.28) circle (0.5pt);
\draw [fill=black] (-3.94,3.28) circle (0.5pt);
\draw [fill=black] (-4.44,3.28) circle (1.5pt);
\draw [fill=black] (-3.44,3.3) circle (1.5pt);
\draw [fill=black] (-3.7,3.3) circle (1.5pt);
\draw [fill=black] (-1.68,3.28) circle (0.5pt);
\draw [fill=black] (-1.52,3.28) circle (0.5pt);
\draw [fill=black] (-1.4,3.28) circle (0.5pt);
\draw [fill=black] (-1.9,3.28) circle (1.5pt);
\draw [fill=black] (-0.9,3.3) circle (1.5pt);
\draw [fill=black] (-1.16,3.3) circle (1.5pt);
\draw [fill=black] (5.68,5.46) circle (0.5pt);
\draw [fill=black] (5.84,5.46) circle (0.5pt);
\draw [fill=black] (5.96,5.46) circle (0.5pt);
\draw [fill=black] (5.46,5.46) circle (1.5pt);
\draw [fill=black] (6.46,5.48) circle (1.5pt);
\draw [fill=black] (6.2,5.48) circle (1.5pt);
\draw [fill=black] (8.1,5.48) circle (0.5pt);
\draw [fill=black] (8.26,5.48) circle (0.5pt);
\draw [fill=black] (8.38,5.48) circle (0.5pt);
\draw [fill=black] (7.88,5.48) circle (1.5pt);
\draw [fill=black] (8.88,5.5) circle (1.5pt);
\draw [fill=black] (8.62,5.5) circle (1.5pt);
\draw [fill=black] (5.64,3.32) circle (0.5pt);
\draw [fill=black] (5.8,3.32) circle (0.5pt);
\draw [fill=black] (5.92,3.32) circle (0.5pt);
\draw [fill=black] (5.42,3.32) circle (1.5pt);
\draw [fill=black] (6.42,3.34) circle (1.5pt);
\draw [fill=black] (6.16,3.34) circle (1.5pt);
\draw [fill=black] (8.18,3.32) circle (0.5pt);
\draw [fill=black] (8.34,3.32) circle (0.5pt);
\draw [fill=black] (8.46,3.32) circle (0.5pt);
\draw [fill=black] (7.96,3.32) circle (1.5pt);
\draw [fill=black] (8.96,3.34) circle (1.5pt);
\draw [fill=black] (8.7,3.34) circle (1.5pt);
\draw [fill=black] (0.68,5.42) circle (0.5pt);
\draw [fill=black] (0.84,5.42) circle (0.5pt);
\draw [fill=black] (0.96,5.42) circle (0.5pt);
\draw [fill=black] (0.46,5.42) circle (1.5pt);
\draw [fill=black] (1.46,5.44) circle (1.5pt);
\draw [fill=black] (1.2,5.44) circle (1.5pt);
\draw [fill=black] (3.1,5.44) circle (0.5pt);
\draw [fill=black] (3.26,5.44) circle (0.5pt);
\draw [fill=black] (3.38,5.44) circle (0.5pt);
\draw [fill=black] (2.88,5.44) circle (1.5pt);
\draw [fill=black] (3.88,5.46) circle (1.5pt);
\draw [fill=black] (3.62,5.46) circle (1.5pt);
\draw [fill=black] (0.64,3.28) circle (0.5pt);
\draw [fill=black] (0.8,3.28) circle (0.5pt);
\draw [fill=black] (0.92,3.28) circle (0.5pt);
\draw [fill=black] (0.42,3.28) circle (1.5pt);
\draw [fill=black] (1.42,3.3) circle (1.5pt);
\draw [fill=black] (1.16,3.3) circle (1.5pt);
\draw [fill=black] (3.18,3.28) circle (0.5pt);
\draw [fill=black] (3.34,3.28) circle (0.5pt);
\draw [fill=black] (3.46,3.28) circle (0.5pt);
\draw [fill=black] (2.96,3.28) circle (1.5pt);
\draw [fill=black] (3.96,3.3) circle (1.5pt);
\draw [fill=black] (3.7,3.3) circle (1.5pt);
\end{scriptsize}
\end{tikzpicture}

%% file: graph82.tex
\begin{tikzpicture}[line cap=round,line join=round,>=triangle 45,x=1.0cm,y=1.0cm]
\clip(-4.860000000000002,1.5199999999999994) rectangle (9.500000000000004,6.240000000000001);
\draw [rotate around={1.145762838175131:(-3.8999999999999995,5.429999999999998)}] (-3.8999999999999995,5.429999999999998) ellipse (0.7473599204960948cm and 0.5553799157009054cm);
\draw [rotate around={1.1457628381750795:(-1.48,5.450000000000001)}] (-1.48,5.450000000000001) ellipse (0.7473599204961263cm and 0.5553799157009286cm);
\draw [rotate around={1.145762838175105:(-3.9400000000000004,3.289999999999999)}] (-3.9400000000000004,3.289999999999999) ellipse (0.7473599204961175cm and 0.5553799157009224cm);
\draw [rotate around={1.145762838175105:(-1.4000000000000001,3.290000000000001)}] (-1.4000000000000001,3.290000000000001) ellipse (0.7473599204961234cm and 0.5553799157009267cm);
\draw (-1.98,5.44)-- (-3.4,5.44);
\draw [rotate around={1.1457628381751304:(5.96,5.470000000000001)}] (5.96,5.470000000000001) ellipse (0.7473599204961332cm and 0.5553799157009339cm);
\draw [rotate around={1.1457628381750788:(8.380000000000003,5.49)}] (8.380000000000003,5.49) ellipse (0.7473599204961118cm and 0.5553799157009173cm);
\draw [rotate around={1.1457628381751053:(5.92,3.3299999999999996)}] (5.92,3.3299999999999996) ellipse (0.7473599204961027cm and 0.5553799157009112cm);
\draw [rotate around={1.1457628381751028:(8.460000000000003,3.33)}] (8.460000000000003,3.33) ellipse (0.7473599204961426cm and 0.55537991570094cm);
\draw (7.88,5.48)-- (6.46,5.48);
\draw (7.96,3.32)-- (6.42,3.34);
\draw [rotate around={1.1457628381751308:(0.9600000000000003,5.429999999999999)}] (0.9600000000000003,5.429999999999999) ellipse (0.747359920496095cm and 0.5553799157009055cm);
\draw [rotate around={1.1457628381750795:(3.3799999999999994,5.450000000000001)}] (3.3799999999999994,5.450000000000001) ellipse (0.7473599204961185cm and 0.5553799157009228cm);
\draw [rotate around={1.145762838175105:(0.9199999999999995,3.29)}] (0.9199999999999995,3.29) ellipse (0.7473599204961217cm and 0.5553799157009254cm);
\draw [rotate around={1.145762838175105:(3.459999999999999,3.289999999999999)}] (3.459999999999999,3.289999999999999) ellipse (0.7473599204961178cm and 0.5553799157009225cm);
\draw (2.88,5.44)-- (1.46,5.44);
\draw (-3.7,3.3)-- (-3.66,5.44);
\draw (-1.9,3.28)-- (-3.44,3.3);
\draw (1.46,5.44)-- (1.16,3.3);
\draw (2.96,3.28)-- (1.42,3.3);
\draw (-3.4000000000000012,2.8) node[anchor=north west] {$R_{4n-2, 6}$};
\draw (1.5200000000000005,2.82) node[anchor=north west] {$R_{4n-2, 7}$};
\draw (6.560000000000002,2.84) node[anchor=north west] {$R_{4n-2, 8}$};
\begin{scriptsize}
\draw [fill=black] (-4.18,5.42) circle (0.5pt);
\draw [fill=black] (-4.02,5.42) circle (0.5pt);
\draw [fill=black] (-3.9,5.42) circle (0.5pt);
\draw [fill=black] (-4.4,5.42) circle (1.5pt);
\draw [fill=black] (-3.4,5.44) circle (1.5pt);
\draw [fill=black] (-3.66,5.44) circle (1.5pt);
\draw [fill=black] (-1.76,5.44) circle (0.5pt);
\draw [fill=black] (-1.6,5.44) circle (0.5pt);
\draw [fill=black] (-1.48,5.44) circle (0.5pt);
\draw [fill=black] (-1.98,5.44) circle (1.5pt);
\draw [fill=black] (-0.98,5.46) circle (1.5pt);
\draw [fill=black] (-1.24,5.46) circle (1.5pt);
\draw [fill=black] (-4.22,3.28) circle (0.5pt);
\draw [fill=black] (-4.06,3.28) circle (0.5pt);
\draw [fill=black] (-3.94,3.28) circle (0.5pt);
\draw [fill=black] (-4.44,3.28) circle (1.5pt);
\draw [fill=black] (-3.44,3.3) circle (1.5pt);
\draw [fill=black] (-3.7,3.3) circle (1.5pt);
\draw [fill=black] (-1.68,3.28) circle (0.5pt);
\draw [fill=black] (-1.52,3.28) circle (0.5pt);
\draw [fill=black] (-1.4,3.28) circle (0.5pt);
\draw [fill=black] (-1.9,3.28) circle (1.5pt);
\draw [fill=black] (-0.9,3.3) circle (1.5pt);
\draw [fill=black] (-1.16,3.3) circle (1.5pt);
\draw [fill=black] (5.68,5.46) circle (0.5pt);
\draw [fill=black] (5.84,5.46) circle (0.5pt);
\draw [fill=black] (5.96,5.46) circle (0.5pt);
\draw [fill=black] (5.46,5.46) circle (1.5pt);
\draw [fill=black] (6.46,5.48) circle (1.5pt);
\draw [fill=black] (6.2,5.48) circle (1.5pt);
\draw [fill=black] (8.1,5.48) circle (0.5pt);
\draw [fill=black] (8.26,5.48) circle (0.5pt);
\draw [fill=black] (8.38,5.48) circle (0.5pt);
\draw [fill=black] (7.88,5.48) circle (1.5pt);
\draw [fill=black] (8.88,5.5) circle (1.5pt);
\draw [fill=black] (8.62,5.5) circle (1.5pt);
\draw [fill=black] (5.64,3.32) circle (0.5pt);
\draw [fill=black] (5.8,3.32) circle (0.5pt);
\draw [fill=black] (5.92,3.32) circle (0.5pt);
\draw [fill=black] (5.42,3.32) circle (1.5pt);
\draw [fill=black] (6.42,3.34) circle (1.5pt);
\draw [fill=black] (6.16,3.34) circle (1.5pt);
\draw [fill=black] (8.18,3.32) circle (0.5pt);
\draw [fill=black] (8.34,3.32) circle (0.5pt);
\draw [fill=black] (8.46,3.32) circle (0.5pt);
\draw [fill=black] (7.96,3.32) circle (1.5pt);
\draw [fill=black] (8.96,3.34) circle (1.5pt);
\draw [fill=black] (8.7,3.34) circle (1.5pt);
\draw [fill=black] (0.68,5.42) circle (0.5pt);
\draw [fill=black] (0.84,5.42) circle (0.5pt);
\draw [fill=black] (0.96,5.42) circle (0.5pt);
\draw [fill=black] (0.46,5.42) circle (1.5pt);
\draw [fill=black] (1.46,5.44) circle (1.5pt);
\draw [fill=black] (1.2,5.44) circle (1.5pt);
\draw [fill=black] (3.1,5.44) circle (0.5pt);
\draw [fill=black] (3.26,5.44) circle (0.5pt);
\draw [fill=black] (3.38,5.44) circle (0.5pt);
\draw [fill=black] (2.88,5.44) circle (1.5pt);
\draw [fill=black] (3.88,5.46) circle (1.5pt);
\draw [fill=black] (3.62,5.46) circle (1.5pt);
\draw [fill=black] (0.64,3.28) circle (0.5pt);
\draw [fill=black] (0.8,3.28) circle (0.5pt);
\draw [fill=black] (0.92,3.28) circle (0.5pt);
\draw [fill=black] (0.42,3.28) circle (1.5pt);
\draw [fill=black] (1.42,3.3) circle (1.5pt);
\draw [fill=black] (1.16,3.3) circle (1.5pt);
\draw [fill=black] (3.18,3.28) circle (0.5pt);
\draw [fill=black] (3.34,3.28) circle (0.5pt);
\draw [fill=black] (3.46,3.28) circle (0.5pt);
\draw [fill=black] (2.96,3.28) circle (1.5pt);
\draw [fill=black] (3.96,3.3) circle (1.5pt);
\draw [fill=black] (3.7,3.3) circle (1.5pt);
\end{scriptsize}
\end{tikzpicture}

%% file: graph83.tex
\begin{tikzpicture}[line cap=round,line join=round,>=triangle 45,x=1.0cm,y=1.0cm]
\clip(-4.860000000000002,1.5000000000000027) rectangle (9.600000000000003,6.24);
\draw [rotate around={1.145762838175131:(-3.8999999999999995,5.429999999999998)}] (-3.8999999999999995,5.429999999999998) ellipse (0.7473599204960948cm and 0.5553799157009054cm);
\draw [rotate around={1.1457628381750795:(-1.48,5.450000000000001)}] (-1.48,5.450000000000001) ellipse (0.7473599204961263cm and 0.5553799157009286cm);
\draw [rotate around={1.145762838175105:(-3.9400000000000004,3.289999999999999)}] (-3.9400000000000004,3.289999999999999) ellipse (0.7473599204961175cm and 0.5553799157009224cm);
\draw [rotate around={1.145762838175105:(-1.4000000000000001,3.290000000000001)}] (-1.4000000000000001,3.290000000000001) ellipse (0.7473599204961234cm and 0.5553799157009267cm);
\draw (-1.98,5.44)-- (-3.4,5.44);
\draw (-3.4,5.44)-- (-3.44,3.3);
\draw [rotate around={1.1457628381751304:(5.96,5.470000000000001)}] (5.96,5.470000000000001) ellipse (0.7473599204961332cm and 0.5553799157009339cm);
\draw [rotate around={1.1457628381750788:(8.380000000000003,5.49)}] (8.380000000000003,5.49) ellipse (0.7473599204961118cm and 0.5553799157009173cm);
\draw [rotate around={1.1457628381751053:(5.92,3.3299999999999996)}] (5.92,3.3299999999999996) ellipse (0.7473599204961027cm and 0.5553799157009112cm);
\draw [rotate around={1.1457628381751028:(8.460000000000003,3.33)}] (8.460000000000003,3.33) ellipse (0.7473599204961426cm and 0.55537991570094cm);
\draw (7.88,5.48)-- (6.46,5.48);
\draw (6.46,5.48)-- (6.42,3.34);
\draw (7.96,3.32)-- (6.42,3.34);
\draw [rotate around={1.1457628381751308:(0.9600000000000003,5.429999999999999)}] (0.9600000000000003,5.429999999999999) ellipse (0.747359920496095cm and 0.5553799157009055cm);
\draw [rotate around={1.1457628381750795:(3.3799999999999994,5.450000000000001)}] (3.3799999999999994,5.450000000000001) ellipse (0.7473599204961185cm and 0.5553799157009228cm);
\draw [rotate around={1.145762838175105:(0.9199999999999995,3.29)}] (0.9199999999999995,3.29) ellipse (0.7473599204961217cm and 0.5553799157009254cm);
\draw [rotate around={1.145762838175105:(3.459999999999999,3.289999999999999)}] (3.459999999999999,3.289999999999999) ellipse (0.7473599204961178cm and 0.5553799157009225cm);
\draw (2.88,5.44)-- (1.46,5.44);
\draw (-1.98,5.44)-- (-3.44,3.3);
\draw (1.46,5.44)-- (1.42,3.3);
\draw (1.42,3.3)-- (2.96,3.28);
\draw (2.96,3.28)-- (2.88,5.44);
\draw (7.88,5.48)-- (6.42,3.34);
\draw (-3.4000000000000017,2.860000000000002) node[anchor=north west] {$R_{4n-4,9}$};
\draw (1.46,2.920000000000002) node[anchor=north west] {$R_{4n-4,10}$};
\draw (6.480000000000002,2.9600000000000017) node[anchor=north west] {$R_{4n-4,11}$};
\begin{scriptsize}
\draw [fill=black] (-4.18,5.42) circle (0.5pt);
\draw [fill=black] (-4.02,5.42) circle (0.5pt);
\draw [fill=black] (-3.9,5.42) circle (0.5pt);
\draw [fill=black] (-4.4,5.42) circle (1.5pt);
\draw [fill=black] (-3.4,5.44) circle (1.5pt);
\draw [fill=black] (-3.66,5.44) circle (1.5pt);
\draw [fill=black] (-1.76,5.44) circle (0.5pt);
\draw [fill=black] (-1.6,5.44) circle (0.5pt);
\draw [fill=black] (-1.48,5.44) circle (0.5pt);
\draw [fill=black] (-1.98,5.44) circle (1.5pt);
\draw [fill=black] (-0.98,5.46) circle (1.5pt);
\draw [fill=black] (-1.24,5.46) circle (1.5pt);
\draw [fill=black] (-4.22,3.28) circle (0.5pt);
\draw [fill=black] (-4.06,3.28) circle (0.5pt);
\draw [fill=black] (-3.94,3.28) circle (0.5pt);
\draw [fill=black] (-4.44,3.28) circle (1.5pt);
\draw [fill=black] (-3.44,3.3) circle (1.5pt);
\draw [fill=black] (-3.7,3.3) circle (1.5pt);
\draw [fill=black] (-1.68,3.28) circle (0.5pt);
\draw [fill=black] (-1.52,3.28) circle (0.5pt);
\draw [fill=black] (-1.4,3.28) circle (0.5pt);
\draw [fill=black] (-1.9,3.28) circle (1.5pt);
\draw [fill=black] (-0.9,3.3) circle (1.5pt);
\draw [fill=black] (-1.16,3.3) circle (1.5pt);
\draw [fill=black] (5.68,5.46) circle (0.5pt);
\draw [fill=black] (5.84,5.46) circle (0.5pt);
\draw [fill=black] (5.96,5.46) circle (0.5pt);
\draw [fill=black] (5.46,5.46) circle (1.5pt);
\draw [fill=black] (6.46,5.48) circle (1.5pt);
\draw [fill=black] (6.2,5.48) circle (1.5pt);
\draw [fill=black] (8.1,5.48) circle (0.5pt);
\draw [fill=black] (8.26,5.48) circle (0.5pt);
\draw [fill=black] (8.38,5.48) circle (0.5pt);
\draw [fill=black] (7.88,5.48) circle (1.5pt);
\draw [fill=black] (8.88,5.5) circle (1.5pt);
\draw [fill=black] (8.62,5.5) circle (1.5pt);
\draw [fill=black] (5.64,3.32) circle (0.5pt);
\draw [fill=black] (5.8,3.32) circle (0.5pt);
\draw [fill=black] (5.92,3.32) circle (0.5pt);
\draw [fill=black] (5.42,3.32) circle (1.5pt);
\draw [fill=black] (6.42,3.34) circle (1.5pt);
\draw [fill=black] (6.16,3.34) circle (1.5pt);
\draw [fill=black] (8.18,3.32) circle (0.5pt);
\draw [fill=black] (8.34,3.32) circle (0.5pt);
\draw [fill=black] (8.46,3.32) circle (0.5pt);
\draw [fill=black] (7.96,3.32) circle (1.5pt);
\draw [fill=black] (8.96,3.34) circle (1.5pt);
\draw [fill=black] (8.7,3.34) circle (1.5pt);
\draw [fill=black] (0.68,5.42) circle (0.5pt);
\draw [fill=black] (0.84,5.42) circle (0.5pt);
\draw [fill=black] (0.96,5.42) circle (0.5pt);
\draw [fill=black] (0.46,5.42) circle (1.5pt);
\draw [fill=black] (1.46,5.44) circle (1.5pt);
\draw [fill=black] (1.2,5.44) circle (1.5pt);
\draw [fill=black] (3.1,5.44) circle (0.5pt);
\draw [fill=black] (3.26,5.44) circle (0.5pt);
\draw [fill=black] (3.38,5.44) circle (0.5pt);
\draw [fill=black] (2.88,5.44) circle (1.5pt);
\draw [fill=black] (3.88,5.46) circle (1.5pt);
\draw [fill=black] (3.62,5.46) circle (1.5pt);
\draw [fill=black] (0.64,3.28) circle (0.5pt);
\draw [fill=black] (0.8,3.28) circle (0.5pt);
\draw [fill=black] (0.92,3.28) circle (0.5pt);
\draw [fill=black] (0.42,3.28) circle (1.5pt);
\draw [fill=black] (1.42,3.3) circle (1.5pt);
\draw [fill=black] (1.16,3.3) circle (1.5pt);
\draw [fill=black] (3.18,3.28) circle (0.5pt);
\draw [fill=black] (3.34,3.28) circle (0.5pt);
\draw [fill=black] (3.46,3.28) circle (0.5pt);
\draw [fill=black] (2.96,3.28) circle (1.5pt);
\draw [fill=black] (3.96,3.3) circle (1.5pt);
\draw [fill=black] (3.7,3.3) circle (1.5pt);
\end{scriptsize}
\end{tikzpicture}

%% file: graph84.tex
\begin{tikzpicture}[line cap=round,line join=round,>=triangle 45,x=1.0cm,y=1.0cm]
\clip(-4.860000000000002,1.4800000000000026) rectangle (9.700000000000003,6.24);
\draw [rotate around={1.145762838175131:(-3.8999999999999995,5.429999999999998)}] (-3.8999999999999995,5.429999999999998) ellipse (0.7473599204960948cm and 0.5553799157009054cm);
\draw [rotate around={1.1457628381750795:(-1.48,5.450000000000001)}] (-1.48,5.450000000000001) ellipse (0.7473599204961263cm and 0.5553799157009286cm);
\draw [rotate around={1.145762838175105:(-3.9400000000000004,3.289999999999999)}] (-3.9400000000000004,3.289999999999999) ellipse (0.7473599204961175cm and 0.5553799157009224cm);
\draw [rotate around={1.145762838175105:(-1.4000000000000001,3.290000000000001)}] (-1.4000000000000001,3.290000000000001) ellipse (0.7473599204961234cm and 0.5553799157009267cm);
\draw (-1.98,5.44)-- (-3.4,5.44);
\draw (-3.4,5.44)-- (-3.44,3.3);
\draw [rotate around={1.1457628381750788:(8.380000000000003,5.49)}] (8.380000000000003,5.49) ellipse (0.7473599204961118cm and 0.5553799157009173cm);
\draw [rotate around={2.202598161765806:(8.440000000000003,3.3200000000000003)}] (8.440000000000003,3.3200000000000003) ellipse (0.7697019896443038cm and 0.5671341577284595cm);
\draw (7.88,5.48)-- (6.46,5.48);
\draw [rotate around={0.0:(0.9600000000000002,5.419999999999999)}] (0.9600000000000002,5.419999999999999) ellipse (0.7573491155352758cm and 0.5688388900225353cm);
\draw [rotate around={1.1691393279073934:(3.39,5.450000000000001)}] (3.39,5.450000000000001) ellipse (0.7382531117140578cm and 0.5521029405423266cm);
\draw [rotate around={0.0:(0.919999999999999,3.2799999999999994)}] (0.919999999999999,3.2799999999999994) ellipse (0.757349115535285cm and 0.5688388900225422cm);
\draw [rotate around={1.101706115206375:(3.440000000000001,3.2900000000000005)}] (3.440000000000001,3.2900000000000005) ellipse (0.7656657234962463cm and 0.5619110251072007cm);
\draw (2.9,5.44)-- (1.46,5.42);
\draw (-1.9,3.28)-- (-3.4,5.44);
\draw (1.2,5.44)-- (2.92,3.28);
\draw [rotate around={0.8814039965821584:(5.810000000000006,5.47)}] (5.810000000000006,5.47) ellipse (0.8612988864660768cm and 0.5650095325104786cm);
\draw [rotate around={0.881403996582139:(5.9300000000000015,3.3099999999999996)}] (5.9300000000000015,3.3099999999999996) ellipse (0.8771331588907284cm and 0.5888655011338544cm);
\draw (5.9,5.48)-- (6.02,3.32);
\draw (6.16,5.48)-- (7.92,3.3);
\draw (1.4800000000000006,2.9800000000000018) node[anchor=north west] {$R_{4n-4,13}$};
\draw (-3.4200000000000013,3.020000000000002) node[anchor=north west] {$R_{4n-4,12}$};
\draw (6.620000000000003,2.9800000000000018) node[anchor=north west] {$R_{4n-4,14}$};
\draw (1.16,3.3)-- (1.2,5.44);
\begin{scriptsize}
\draw [fill=black] (-4.18,5.42) circle (0.5pt);
\draw [fill=black] (-4.02,5.42) circle (0.5pt);
\draw [fill=black] (-3.9,5.42) circle (0.5pt);
\draw [fill=black] (-4.4,5.42) circle (1.5pt);
\draw [fill=black] (-3.4,5.44) circle (1.5pt);
\draw [fill=black] (-3.66,5.44) circle (1.5pt);
\draw [fill=black] (-1.76,5.44) circle (0.5pt);
\draw [fill=black] (-1.6,5.44) circle (0.5pt);
\draw [fill=black] (-1.48,5.44) circle (0.5pt);
\draw [fill=black] (-1.98,5.44) circle (1.5pt);
\draw [fill=black] (-0.98,5.46) circle (1.5pt);
\draw [fill=black] (-1.24,5.46) circle (1.5pt);
\draw [fill=black] (-4.22,3.28) circle (0.5pt);
\draw [fill=black] (-4.06,3.28) circle (0.5pt);
\draw [fill=black] (-3.94,3.28) circle (0.5pt);
\draw [fill=black] (-4.44,3.28) circle (1.5pt);
\draw [fill=black] (-3.44,3.3) circle (1.5pt);
\draw [fill=black] (-3.7,3.3) circle (1.5pt);
\draw [fill=black] (-1.68,3.28) circle (0.5pt);
\draw [fill=black] (-1.52,3.28) circle (0.5pt);
\draw [fill=black] (-1.4,3.28) circle (0.5pt);
\draw [fill=black] (-1.9,3.28) circle (1.5pt);
\draw [fill=black] (-0.9,3.3) circle (1.5pt);
\draw [fill=black] (-1.16,3.3) circle (1.5pt);
\draw [fill=black] (5.38,5.46) circle (0.5pt);
\draw [fill=black] (5.54,5.46) circle (0.5pt);
\draw [fill=black] (5.66,5.46) circle (0.5pt);
\draw [fill=black] (5.16,5.46) circle (1.5pt);
\draw [fill=black] (6.46,5.48) circle (1.5pt);
\draw [fill=black] (5.9,5.48) circle (1.5pt);
\draw [fill=black] (8.1,5.48) circle (0.5pt);
\draw [fill=black] (8.26,5.48) circle (0.5pt);
\draw [fill=black] (8.38,5.48) circle (0.5pt);
\draw [fill=black] (7.88,5.48) circle (1.5pt);
\draw [fill=black] (8.88,5.5) circle (1.5pt);
\draw [fill=black] (8.62,5.5) circle (1.5pt);
\draw [fill=black] (8.18,3.32) circle (0.5pt);
\draw [fill=black] (8.34,3.32) circle (0.5pt);
\draw [fill=black] (8.46,3.32) circle (0.5pt);
\draw [fill=black] (7.92,3.3) circle (1.5pt);
\draw [fill=black] (8.96,3.34) circle (1.5pt);
\draw [fill=black] (8.7,3.34) circle (1.5pt);
\draw [fill=black] (0.68,5.42) circle (0.5pt);
\draw [fill=black] (0.84,5.42) circle (0.5pt);
\draw [fill=black] (0.96,5.42) circle (0.5pt);
\draw [fill=black] (0.46,5.42) circle (1.5pt);
\draw [fill=black] (1.46,5.42) circle (1.5pt);
\draw [fill=black] (1.2,5.44) circle (1.5pt);
\draw [fill=black] (3.1,5.44) circle (0.5pt);
\draw [fill=black] (3.26,5.44) circle (0.5pt);
\draw [fill=black] (3.38,5.44) circle (0.5pt);
\draw [fill=black] (2.9,5.44) circle (1.5pt);
\draw [fill=black] (3.88,5.46) circle (1.5pt);
\draw [fill=black] (3.62,5.46) circle (1.5pt);
\draw [fill=black] (0.64,3.28) circle (0.5pt);
\draw [fill=black] (0.8,3.28) circle (0.5pt);
\draw [fill=black] (0.92,3.28) circle (0.5pt);
\draw [fill=black] (0.42,3.28) circle (1.5pt);
\draw [fill=black] (1.42,3.28) circle (1.5pt);
\draw [fill=black] (1.16,3.3) circle (1.5pt);
\draw [fill=black] (3.18,3.28) circle (0.5pt);
\draw [fill=black] (3.34,3.28) circle (0.5pt);
\draw [fill=black] (3.46,3.28) circle (0.5pt);
\draw [fill=black] (2.92,3.28) circle (1.5pt);
\draw [fill=black] (3.96,3.3) circle (1.5pt);
\draw [fill=black] (3.7,3.3) circle (1.5pt);
\draw [fill=black] (6.16,5.48) circle (1.5pt);
\draw [fill=black] (5.5,3.3) circle (0.5pt);
\draw [fill=black] (5.66,3.3) circle (0.5pt);
\draw [fill=black] (5.78,3.3) circle (0.5pt);
\draw [fill=black] (5.28,3.3) circle (1.5pt);
\draw [fill=black] (6.58,3.32) circle (1.5pt);
\draw [fill=black] (6.02,3.32) circle (1.5pt);
\draw [fill=black] (6.28,3.32) circle (1.5pt);
\end{scriptsize}
\end{tikzpicture}

%% file: graph85.tex
\begin{tikzpicture}[line cap=round,line join=round,>=triangle 45,x=1.0cm,y=1.0cm]
\clip(-4.860000000000002,1.3400000000000052) rectangle (9.520000000000003,6.239999999999999);
\draw [rotate around={1.1017061152063998:(-3.879999999999999,5.430000000000001)}] (-3.879999999999999,5.430000000000001) ellipse (0.7492539944357897cm and 0.539334356571126cm);
\draw [rotate around={1.1457628381750795:(-1.48,5.450000000000001)}] (-1.48,5.450000000000001) ellipse (0.7473599204961263cm and 0.5553799157009286cm);
\draw [rotate around={1.145762838175105:(-3.9400000000000004,3.289999999999999)}] (-3.9400000000000004,3.289999999999999) ellipse (0.7473599204961175cm and 0.5553799157009224cm);
\draw [rotate around={1.145762838175105:(-1.4000000000000001,3.290000000000001)}] (-1.4000000000000001,3.290000000000001) ellipse (0.7473599204961234cm and 0.5553799157009267cm);
\draw (-1.98,5.44)-- (-3.36,5.44);
\draw [rotate around={1.1457628381751304:(5.96,5.470000000000001)}] (5.96,5.470000000000001) ellipse (0.7473599204961332cm and 0.5553799157009339cm);
\draw [rotate around={1.1457628381750788:(8.380000000000003,5.49)}] (8.380000000000003,5.49) ellipse (0.7473599204961118cm and 0.5553799157009173cm);
\draw [rotate around={1.1457628381751053:(5.92,3.3299999999999996)}] (5.92,3.3299999999999996) ellipse (0.7473599204961027cm and 0.5553799157009112cm);
\draw [rotate around={2.202598161765806:(8.440000000000003,3.3200000000000003)}] (8.440000000000003,3.3200000000000003) ellipse (0.7697019896443038cm and 0.5671341577284595cm);
\draw (7.88,5.48)-- (6.46,5.48);
\draw (6.46,5.48)-- (6.42,3.34);
\draw (7.92,3.3)-- (6.42,3.34);
\draw [rotate around={1.1457628381751308:(0.9600000000000003,5.429999999999999)}] (0.9600000000000003,5.429999999999999) ellipse (0.747359920496095cm and 0.5553799157009055cm);
\draw [rotate around={1.1457628381750795:(3.3799999999999994,5.450000000000001)}] (3.3799999999999994,5.450000000000001) ellipse (0.7473599204961185cm and 0.5553799157009228cm);
\draw [rotate around={1.145762838175105:(0.9199999999999995,3.29)}] (0.9199999999999995,3.29) ellipse (0.7473599204961217cm and 0.5553799157009254cm);
\draw [rotate around={1.0809241866606896:(3.4299999999999997,3.2900000000000005)}] (3.4299999999999997,3.2900000000000005) ellipse (0.7748617289560169cm and 0.5651643115065772cm);
\draw (2.88,5.44)-- (1.46,5.44);
\draw (1.46,5.44)-- (1.42,3.3);
\draw (1.42,3.3)-- (2.9,3.28);
\draw (7.88,5.48)-- (6.42,3.34);
\draw (-1.9,3.28)-- (-3.44,3.3);
\draw (2.88,5.44)-- (1.42,3.3);
\draw (7.88,5.48)-- (7.92,3.3);
\draw (7.92,3.3)-- (6.46,5.48);
\draw (2.88,5.44)-- (2.9,3.28);
\draw (-3.4200000000000013,2.9400000000000035) node[anchor=north west] {$R_{4n-4,15}$};
\draw (1.4600000000000006,2.8800000000000034) node[anchor=north west] {$R_{4n-4,16}$};
\draw (6.480000000000002,2.980000000000003) node[anchor=north west] {$R_{4n-4,17}$};
\draw (-3.66,3.3)-- (-3.36,5.44);
\draw (-3.66,3.3)-- (-1.98,5.44);
\begin{scriptsize}
\draw [fill=black] (-4.18,5.42) circle (0.5pt);
\draw [fill=black] (-4.02,5.42) circle (0.5pt);
\draw [fill=black] (-3.9,5.42) circle (0.5pt);
\draw [fill=black] (-4.4,5.42) circle (1.5pt);
\draw [fill=black] (-3.36,5.44) circle (1.5pt);
\draw [fill=black] (-3.6,5.46) circle (1.5pt);
\draw [fill=black] (-1.76,5.44) circle (0.5pt);
\draw [fill=black] (-1.6,5.44) circle (0.5pt);
\draw [fill=black] (-1.48,5.44) circle (0.5pt);
\draw [fill=black] (-1.98,5.44) circle (1.5pt);
\draw [fill=black] (-0.98,5.46) circle (1.5pt);
\draw [fill=black] (-1.24,5.46) circle (1.5pt);
\draw [fill=black] (-4.22,3.28) circle (0.5pt);
\draw [fill=black] (-4.06,3.28) circle (0.5pt);
\draw [fill=black] (-3.94,3.28) circle (0.5pt);
\draw [fill=black] (-4.44,3.28) circle (1.5pt);
\draw [fill=black] (-3.44,3.3) circle (1.5pt);
\draw [fill=black] (-3.66,3.3) circle (1.5pt);
\draw [fill=black] (-1.68,3.28) circle (0.5pt);
\draw [fill=black] (-1.52,3.28) circle (0.5pt);
\draw [fill=black] (-1.4,3.28) circle (0.5pt);
\draw [fill=black] (-1.9,3.28) circle (1.5pt);
\draw [fill=black] (-0.9,3.3) circle (1.5pt);
\draw [fill=black] (-1.16,3.3) circle (1.5pt);
\draw [fill=black] (5.68,5.46) circle (0.5pt);
\draw [fill=black] (5.84,5.46) circle (0.5pt);
\draw [fill=black] (5.96,5.46) circle (0.5pt);
\draw [fill=black] (5.46,5.46) circle (1.5pt);
\draw [fill=black] (6.46,5.48) circle (1.5pt);
\draw [fill=black] (6.2,5.48) circle (1.5pt);
\draw [fill=black] (8.1,5.48) circle (0.5pt);
\draw [fill=black] (8.26,5.48) circle (0.5pt);
\draw [fill=black] (8.38,5.48) circle (0.5pt);
\draw [fill=black] (7.88,5.48) circle (1.5pt);
\draw [fill=black] (8.88,5.5) circle (1.5pt);
\draw [fill=black] (8.62,5.5) circle (1.5pt);
\draw [fill=black] (5.64,3.32) circle (0.5pt);
\draw [fill=black] (5.8,3.32) circle (0.5pt);
\draw [fill=black] (5.92,3.32) circle (0.5pt);
\draw [fill=black] (5.42,3.32) circle (1.5pt);
\draw [fill=black] (6.42,3.34) circle (1.5pt);
\draw [fill=black] (6.16,3.34) circle (1.5pt);
\draw [fill=black] (8.18,3.32) circle (0.5pt);
\draw [fill=black] (8.34,3.32) circle (0.5pt);
\draw [fill=black] (8.46,3.32) circle (0.5pt);
\draw [fill=black] (7.92,3.3) circle (1.5pt);
\draw [fill=black] (8.96,3.34) circle (1.5pt);
\draw [fill=black] (8.7,3.34) circle (1.5pt);
\draw [fill=black] (0.68,5.42) circle (0.5pt);
\draw [fill=black] (0.84,5.42) circle (0.5pt);
\draw [fill=black] (0.96,5.42) circle (0.5pt);
\draw [fill=black] (0.46,5.42) circle (1.5pt);
\draw [fill=black] (1.46,5.44) circle (1.5pt);
\draw [fill=black] (1.2,5.44) circle (1.5pt);
\draw [fill=black] (3.1,5.44) circle (0.5pt);
\draw [fill=black] (3.26,5.44) circle (0.5pt);
\draw [fill=black] (3.38,5.44) circle (0.5pt);
\draw [fill=black] (2.88,5.44) circle (1.5pt);
\draw [fill=black] (3.88,5.46) circle (1.5pt);
\draw [fill=black] (3.62,5.46) circle (1.5pt);
\draw [fill=black] (0.64,3.28) circle (0.5pt);
\draw [fill=black] (0.8,3.28) circle (0.5pt);
\draw [fill=black] (0.92,3.28) circle (0.5pt);
\draw [fill=black] (0.42,3.28) circle (1.5pt);
\draw [fill=black] (1.42,3.3) circle (1.5pt);
\draw [fill=black] (1.16,3.3) circle (1.5pt);
\draw [fill=black] (3.18,3.28) circle (0.5pt);
\draw [fill=black] (3.34,3.28) circle (0.5pt);
\draw [fill=black] (3.46,3.28) circle (0.5pt);
\draw [fill=black] (2.9,3.28) circle (1.5pt);
\draw [fill=black] (3.96,3.3) circle (1.5pt);
\draw [fill=black] (3.7,3.3) circle (1.5pt);
\end{scriptsize}
\end{tikzpicture}

%% file: graph86.tex
\begin{tikzpicture}[line cap=round,line join=round,>=triangle 45,x=1.0cm,y=1.0cm]
\clip(-4.860000000000002,1.4600000000000026) rectangle (9.500000000000004,6.24);
\draw [rotate around={1.145762838175131:(-3.8999999999999995,5.429999999999998)}] (-3.8999999999999995,5.429999999999998) ellipse (0.7473599204960948cm and 0.5553799157009054cm);
\draw [rotate around={1.1457628381750795:(-1.4199999999999973,5.450000000000001)}] (-1.4199999999999973,5.450000000000001) ellipse (0.7473599204961262cm and 0.5553799157009286cm);
\draw [rotate around={1.145762838175105:(-3.9400000000000004,3.289999999999999)}] (-3.9400000000000004,3.289999999999999) ellipse (0.7473599204961175cm and 0.5553799157009224cm);
\draw [rotate around={1.145762838175105:(-1.4000000000000001,3.290000000000001)}] (-1.4000000000000001,3.290000000000001) ellipse (0.7473599204961234cm and 0.5553799157009267cm);
\draw [rotate around={1.1233027140754563:(5.9500000000000055,5.4700000000000015)}] (5.9500000000000055,5.4700000000000015) ellipse (0.7564979717662241cm and 0.5586494261040549cm);
\draw [rotate around={3.3664606634297805:(8.370000000000001,5.470000000000001)}] (8.370000000000001,5.470000000000001) ellipse (0.7648549453707757cm and 0.5692126908793644cm);
\draw [rotate around={-1.1457628381751046:(5.92,3.31)}] (5.92,3.31) ellipse (0.7673392487442487cm and 0.5819875622926957cm);
\draw [rotate around={0.0:(8.409999999999997,3.3)}] (8.409999999999997,3.3) ellipse (0.7988120476976098cm and 0.6148176051048363cm);
\draw (7.9,3.3)-- (6.42,3.3);
\draw (-1.9,3.28)-- (-3.44,3.3);
\draw (-1.92,5.44)-- (-3.4,5.44);
\draw (-1.66,5.42)-- (-1.64,3.28);
\draw (-3.66,5.44)-- (-3.7,3.3);
\draw [rotate around={1.1457628381751304:(0.899999999999999,5.429999999999999)}] (0.899999999999999,5.429999999999999) ellipse (0.7473599204961104cm and 0.5553799157009169cm);
\draw [rotate around={1.1457628381750795:(3.3799999999999994,5.450000000000001)}] (3.3799999999999994,5.450000000000001) ellipse (0.7473599204961185cm and 0.5553799157009228cm);
\draw [rotate around={1.145762838175105:(0.8599999999999995,3.29)}] (0.8599999999999995,3.29) ellipse (0.7473599204961197cm and 0.5553799157009239cm);
\draw [rotate around={1.145762838175105:(3.3999999999999995,3.29)}] (3.3999999999999995,3.29) ellipse (0.7473599204961178cm and 0.5553799157009225cm);
\draw (2.9,3.28)-- (1.36,3.3);
\draw (2.88,5.44)-- (1.4,5.44);
\draw (1.14,5.44)-- (1.1,3.3);
\draw (2.88,5.44)-- (3.16,3.28);
\draw (8.88,5.5)-- (8.92,3.3);
\draw (5.44,5.46)-- (5.42,3.32);
\draw (8.12,5.46)-- (6.2,3.32);
\draw (6.2,5.48)-- (8.12,3.3);
\draw (7.86,5.44)-- (6.46,5.48);
\draw (-3.4400000000000013,2.900000000000002) node[anchor=north west] {$R_{4n-4,18}$};
\draw (1.3800000000000006,2.880000000000002) node[anchor=north west] {$R_{4n-4,19}$};
\draw (6.540000000000003,2.920000000000002) node[anchor=north west] {$R_{4n-4,20}$};
\begin{scriptsize}
\draw [fill=black] (-4.18,5.42) circle (0.5pt);
\draw [fill=black] (-4.02,5.42) circle (0.5pt);
\draw [fill=black] (-3.9,5.42) circle (0.5pt);
\draw [fill=black] (-4.4,5.42) circle (1.5pt);
\draw [fill=black] (-3.4,5.44) circle (1.5pt);
\draw [fill=black] (-3.66,5.44) circle (1.5pt);
\draw [fill=black] (-1.48,5.42) circle (0.5pt);
\draw [fill=black] (-1.32,5.42) circle (0.5pt);
\draw [fill=black] (-1.2,5.42) circle (0.5pt);
\draw [fill=black] (-1.92,5.44) circle (1.5pt);
\draw [fill=black] (-0.92,5.46) circle (1.5pt);
\draw [fill=black] (-1.66,5.42) circle (1.5pt);
\draw [fill=black] (-4.22,3.28) circle (0.5pt);
\draw [fill=black] (-4.06,3.28) circle (0.5pt);
\draw [fill=black] (-3.94,3.28) circle (0.5pt);
\draw [fill=black] (-4.44,3.28) circle (1.5pt);
\draw [fill=black] (-3.44,3.3) circle (1.5pt);
\draw [fill=black] (-3.7,3.3) circle (1.5pt);
\draw [fill=black] (-1.42,3.66) circle (0.5pt);
\draw [fill=black] (-1.26,3.66) circle (0.5pt);
\draw [fill=black] (-1.14,3.66) circle (0.5pt);
\draw [fill=black] (-1.9,3.28) circle (1.5pt);
\draw [fill=black] (-0.9,3.3) circle (1.5pt);
\draw [fill=black] (-1.64,3.28) circle (1.5pt);
\draw [fill=black] (5.68,5.46) circle (0.5pt);
\draw [fill=black] (5.84,5.46) circle (0.5pt);
\draw [fill=black] (5.96,5.46) circle (0.5pt);
\draw [fill=black] (5.44,5.46) circle (1.5pt);
\draw [fill=black] (6.46,5.48) circle (1.5pt);
\draw [fill=black] (6.2,5.48) circle (1.5pt);
\draw [fill=black] (8.7,5.48) circle (0.5pt);
\draw [fill=black] (8.5,5.46) circle (0.5pt);
\draw [fill=black] (7.86,5.44) circle (1.5pt);
\draw [fill=black] (8.88,5.5) circle (1.5pt);
\draw [fill=black] (8.12,5.46) circle (1.5pt);
\draw [fill=black] (5.64,3.32) circle (0.5pt);
\draw [fill=black] (5.8,3.32) circle (0.5pt);
\draw [fill=black] (5.92,3.32) circle (0.5pt);
\draw [fill=black] (5.42,3.32) circle (1.5pt);
\draw [fill=black] (6.42,3.3) circle (1.5pt);
\draw [fill=black] (6.2,3.32) circle (1.5pt);
\draw [fill=black] (8.88,3.32) circle (0.5pt);
\draw [fill=black] (8.56,3.32) circle (0.5pt);
\draw [fill=black] (8.72,3.32) circle (0.5pt);
\draw [fill=black] (7.9,3.3) circle (1.5pt);
\draw [fill=black] (8.92,3.3) circle (1.5pt);
\draw [fill=black] (8.12,3.3) circle (1.5pt);
\draw [fill=black] (0.62,5.42) circle (0.5pt);
\draw [fill=black] (0.78,5.42) circle (0.5pt);
\draw [fill=black] (0.9,5.42) circle (0.5pt);
\draw [fill=black] (0.4,5.42) circle (1.5pt);
\draw [fill=black] (1.4,5.44) circle (1.5pt);
\draw [fill=black] (1.14,5.44) circle (1.5pt);
\draw [fill=black] (3.32,5.42) circle (0.5pt);
\draw [fill=black] (3.48,5.42) circle (0.5pt);
\draw [fill=black] (3.6,5.42) circle (0.5pt);
\draw [fill=black] (2.88,5.44) circle (1.5pt);
\draw [fill=black] (3.88,5.46) circle (1.5pt);
\draw [fill=black] (3.14,5.42) circle (1.5pt);
\draw [fill=black] (0.58,3.28) circle (0.5pt);
\draw [fill=black] (0.74,3.28) circle (0.5pt);
\draw [fill=black] (0.86,3.28) circle (0.5pt);
\draw [fill=black] (0.36,3.28) circle (1.5pt);
\draw [fill=black] (1.36,3.3) circle (1.5pt);
\draw [fill=black] (1.1,3.3) circle (1.5pt);
\draw [fill=black] (3.36,3.3) circle (0.5pt);
\draw [fill=black] (3.52,3.3) circle (0.5pt);
\draw [fill=black] (3.64,3.3) circle (0.5pt);
\draw [fill=black] (2.9,3.28) circle (1.5pt);
\draw [fill=black] (3.9,3.3) circle (1.5pt);
\draw [fill=black] (3.16,3.28) circle (1.5pt);
\draw [fill=black] (8.3,5.46) circle (0.5pt);
\draw [fill=black] (8.34,3.34) circle (0.5pt);
\end{scriptsize}
\end{tikzpicture}

%% file: graph87.tex
\begin{tikzpicture}[line cap=round,line join=round,>=triangle 45,x=1.0cm,y=1.0cm]
\clip(-4.860000000000002,1.4600000000000026) rectangle (9.500000000000004,6.24);
\draw [rotate around={1.145762838175131:(-3.8999999999999995,5.429999999999998)}] (-3.8999999999999995,5.429999999999998) ellipse (0.7473599204960948cm and 0.5553799157009054cm);
\draw [rotate around={1.1457628381750795:(-1.4199999999999973,5.450000000000001)}] (-1.4199999999999973,5.450000000000001) ellipse (0.7473599204961262cm and 0.5553799157009286cm);
\draw [rotate around={1.145762838175105:(-3.9400000000000004,3.289999999999999)}] (-3.9400000000000004,3.289999999999999) ellipse (0.7473599204961175cm and 0.5553799157009224cm);
\draw [rotate around={1.145762838175105:(-1.4000000000000001,3.290000000000001)}] (-1.4000000000000001,3.290000000000001) ellipse (0.7473599204961234cm and 0.5553799157009267cm);
\draw (-1.92,5.44)-- (-3.4,5.44);
\draw (-1.66,5.42)-- (-1.64,3.28);
\draw [rotate around={0.0:(0.8999999999999992,5.419999999999999)}] (0.8999999999999992,5.419999999999999) ellipse (0.757349115535283cm and 0.5688388900225407cm);
\draw [rotate around={1.1457628381750795:(3.3799999999999994,5.450000000000001)}] (3.3799999999999994,5.450000000000001) ellipse (0.7473599204961185cm and 0.5553799157009228cm);
\draw [rotate around={1.145762838175105:(0.8599999999999995,3.29)}] (0.8599999999999995,3.29) ellipse (0.7473599204961197cm and 0.5553799157009239cm);
\draw [rotate around={0.0:(3.410000000000001,3.28)}] (3.410000000000001,3.28) ellipse (0.7580291804235784cm and 0.5608103408226713cm);
\draw (2.88,5.44)-- (1.4,5.42);
\draw (-3.4600000000000013,2.900000000000002) node[anchor=north west] {$R_{4n-4,21}$};
\draw (1.3800000000000006,2.880000000000002) node[anchor=north west] {$R_{4n-4,22}$};
\draw (-3.44,3.3)-- (-3.66,5.44);
\draw (-3.66,5.44)-- (-1.9,3.28);
\draw (-1.9,3.28)-- (-3.44,3.3);
\draw (1.4,5.42)-- (1.12,3.28);
\draw (1.12,3.28)-- (2.88,5.44);
\draw [rotate around={0.0:(5.88,5.380000000000001)}] (5.88,5.380000000000001) ellipse (0.7573491155352899cm and 0.5688388900225461cm);
\draw [rotate around={1.145762838175081:(8.359999999999998,5.41)}] (8.359999999999998,5.41) ellipse (0.7473599204960798cm and 0.5553799157008944cm);
\draw [rotate around={1.1457628381750795:(5.839999999999996,3.25)}] (5.839999999999996,3.25) ellipse (0.7473599204961178cm and 0.5553799157009225cm);
\draw [rotate around={0.0:(8.389999999999997,3.24)}] (8.389999999999997,3.24) ellipse (0.758029180423587cm and 0.5608103408226773cm);
\draw (7.86,5.4)-- (6.38,5.38);
\draw (6.360000000000002,2.840000000000002) node[anchor=north west] {$R_{4n-4,23}$};
\draw (6.38,5.38)-- (6.1,3.24);
\draw (7.88,3.24)-- (6.34,3.26);
\draw (6.1,5.34)-- (8.14,3.24);
\draw (3.16,3.28)-- (3.12,5.42);
\draw (2.9,3.28)-- (1.36,3.3);
\draw (3.92,3.28)-- (1.12,5.38);
\begin{scriptsize}
\draw [fill=black] (-4.18,5.42) circle (0.5pt);
\draw [fill=black] (-4.02,5.42) circle (0.5pt);
\draw [fill=black] (-3.9,5.42) circle (0.5pt);
\draw [fill=black] (-4.4,5.42) circle (1.5pt);
\draw [fill=black] (-3.4,5.44) circle (1.5pt);
\draw [fill=black] (-3.66,5.44) circle (1.5pt);
\draw [fill=black] (-1.48,5.42) circle (0.5pt);
\draw [fill=black] (-1.32,5.42) circle (0.5pt);
\draw [fill=black] (-1.2,5.42) circle (0.5pt);
\draw [fill=black] (-1.92,5.44) circle (1.5pt);
\draw [fill=black] (-0.92,5.46) circle (1.5pt);
\draw [fill=black] (-1.66,5.42) circle (1.5pt);
\draw [fill=black] (-4.22,3.28) circle (0.5pt);
\draw [fill=black] (-4.06,3.28) circle (0.5pt);
\draw [fill=black] (-3.94,3.28) circle (0.5pt);
\draw [fill=black] (-4.44,3.28) circle (1.5pt);
\draw [fill=black] (-3.44,3.3) circle (1.5pt);
\draw [fill=black] (-3.7,3.3) circle (1.5pt);
\draw [fill=black] (-1.4,3.3) circle (0.5pt);
\draw [fill=black] (-1.24,3.3) circle (0.5pt);
\draw [fill=black] (-1.12,3.3) circle (0.5pt);
\draw [fill=black] (-1.9,3.28) circle (1.5pt);
\draw [fill=black] (-0.9,3.3) circle (1.5pt);
\draw [fill=black] (-1.64,3.28) circle (1.5pt);
\draw [fill=black] (0.62,5.42) circle (0.5pt);
\draw [fill=black] (0.78,5.42) circle (0.5pt);
\draw [fill=black] (0.9,5.42) circle (0.5pt);
\draw [fill=black] (0.4,5.42) circle (1.5pt);
\draw [fill=black] (1.4,5.42) circle (1.5pt);
\draw [fill=black] (1.12,5.38) circle (1.5pt);
\draw [fill=black] (3.32,5.42) circle (0.5pt);
\draw [fill=black] (3.48,5.42) circle (0.5pt);
\draw [fill=black] (3.6,5.42) circle (0.5pt);
\draw [fill=black] (2.88,5.44) circle (1.5pt);
\draw [fill=black] (3.88,5.46) circle (1.5pt);
\draw [fill=black] (3.12,5.42) circle (1.5pt);
\draw [fill=black] (0.58,3.28) circle (0.5pt);
\draw [fill=black] (0.74,3.28) circle (0.5pt);
\draw [fill=black] (0.86,3.28) circle (0.5pt);
\draw [fill=black] (0.36,3.28) circle (1.5pt);
\draw [fill=black] (1.36,3.3) circle (1.5pt);
\draw [fill=black] (1.12,3.28) circle (1.5pt);
\draw [fill=black] (3.36,3.3) circle (0.5pt);
\draw [fill=black] (3.52,3.3) circle (0.5pt);
\draw [fill=black] (3.64,3.3) circle (0.5pt);
\draw [fill=black] (2.9,3.28) circle (1.5pt);
\draw [fill=black] (3.92,3.28) circle (1.5pt);
\draw [fill=black] (3.16,3.28) circle (1.5pt);
\draw [fill=black] (5.6,5.38) circle (0.5pt);
\draw [fill=black] (5.76,5.38) circle (0.5pt);
\draw [fill=black] (5.88,5.38) circle (0.5pt);
\draw [fill=black] (5.38,5.38) circle (1.5pt);
\draw [fill=black] (6.38,5.38) circle (1.5pt);
\draw [fill=black] (6.1,5.34) circle (1.5pt);
\draw [fill=black] (8.3,5.38) circle (0.5pt);
\draw [fill=black] (8.46,5.38) circle (0.5pt);
\draw [fill=black] (8.58,5.38) circle (0.5pt);
\draw [fill=black] (7.86,5.4) circle (1.5pt);
\draw [fill=black] (8.86,5.42) circle (1.5pt);
\draw [fill=black] (8.1,5.38) circle (1.5pt);
\draw [fill=black] (5.56,3.24) circle (0.5pt);
\draw [fill=black] (5.72,3.24) circle (0.5pt);
\draw [fill=black] (5.84,3.24) circle (0.5pt);
\draw [fill=black] (5.34,3.24) circle (1.5pt);
\draw [fill=black] (6.34,3.26) circle (1.5pt);
\draw [fill=black] (6.1,3.24) circle (1.5pt);
\draw [fill=black] (8.34,3.26) circle (0.5pt);
\draw [fill=black] (8.5,3.26) circle (0.5pt);
\draw [fill=black] (8.62,3.26) circle (0.5pt);
\draw [fill=black] (7.88,3.24) circle (1.5pt);
\draw [fill=black] (8.9,3.24) circle (1.5pt);
\draw [fill=black] (8.14,3.24) circle (1.5pt);
\end{scriptsize}
\end{tikzpicture}

%% file: graph88.tex
\begin{tikzpicture}[line cap=round,line join=round,>=triangle 45,x=1.0cm,y=1.0cm]
\clip(-4.860000000000002,1.4600000000000026) rectangle (9.500000000000004,6.24);
\draw [rotate around={1.145762838175131:(-3.8999999999999995,5.429999999999998)}] (-3.8999999999999995,5.429999999999998) ellipse (0.7473599204960948cm and 0.5553799157009054cm);
\draw [rotate around={1.1457628381750795:(-1.4199999999999973,5.450000000000001)}] (-1.4199999999999973,5.450000000000001) ellipse (0.7473599204961262cm and 0.5553799157009286cm);
\draw [rotate around={1.145762838175105:(-3.9400000000000004,3.289999999999999)}] (-3.9400000000000004,3.289999999999999) ellipse (0.7473599204961175cm and 0.5553799157009224cm);
\draw [rotate around={1.145762838175105:(-1.4000000000000001,3.290000000000001)}] (-1.4000000000000001,3.290000000000001) ellipse (0.7473599204961234cm and 0.5553799157009267cm);
\draw (-1.92,5.44)-- (-3.4,5.44);
\draw (-1.66,5.42)-- (-1.64,3.28);
\draw (-3.4600000000000013,2.900000000000002) node[anchor=north west] {$R_{4n-4,24}$};
\draw (1.3600000000000005,2.880000000000002) node[anchor=north west] {$R_{4n-4,25}$};
\draw (-3.66,5.44)-- (-1.9,3.28);
\draw (6.360000000000002,2.840000000000002) node[anchor=north west] {$R_{4n-4,26}$};
\draw (-3.44,3.3)-- (-0.92,5.46);
\draw (-3.7,3.3)-- (-4.4,5.42);
\draw [rotate around={1.1457628381750795:(1.04,5.45)}] (1.04,5.45) ellipse (0.7473599204961185cm and 0.5553799157009228cm);
\draw [rotate around={1.1457628381751312:(3.5200000000000027,5.47)}] (3.5200000000000027,5.47) ellipse (0.7473599204961179cm and 0.5553799157009227cm);
\draw [rotate around={1.145762838175105:(0.9999999999999996,3.309999999999999)}] (0.9999999999999996,3.309999999999999) ellipse (0.747359920496112cm and 0.5553799157009182cm);
\draw [rotate around={1.1457628381751053:(3.5399999999999983,3.3099999999999987)}] (3.5399999999999983,3.3099999999999987) ellipse (0.7473599204961104cm and 0.5553799157009169cm);
\draw (3.02,5.46)-- (1.54,5.46);
\draw (3.28,5.44)-- (3.3,3.3);
\draw (1.28,5.46)-- (3.04,3.3);
\draw (1.24,3.32)-- (0.54,5.44);
\draw [rotate around={1.1457628381751304:(5.939999999999999,5.429999999999999)}] (5.939999999999999,5.429999999999999) ellipse (0.7473599204961031cm and 0.5553799157009114cm);
\draw [rotate around={1.1457628381750797:(8.42,5.449999999999999)}] (8.42,5.449999999999999) ellipse (0.7473599204961104cm and 0.5553799157009169cm);
\draw [rotate around={1.1457628381751053:(5.899999999999998,3.2899999999999996)}] (5.899999999999998,3.2899999999999996) ellipse (0.7473599204961027cm and 0.5553799157009112cm);
\draw [rotate around={1.1691393279074207:(8.45,3.29)}] (8.45,3.29) ellipse (0.7382531117140095cm and 0.5521029405422909cm);
\draw (7.92,5.44)-- (6.44,5.44);
\draw (7.96,3.28)-- (6.4,3.3);
\draw (6.4,3.3)-- (6.2,5.44);
\draw (8.2,3.28)-- (7.92,5.44);
\draw (6.14,3.3)-- (8.18,5.42);
\begin{scriptsize}
\draw [fill=black] (-4.18,5.42) circle (0.5pt);
\draw [fill=black] (-4.02,5.42) circle (0.5pt);
\draw [fill=black] (-3.9,5.42) circle (0.5pt);
\draw [fill=black] (-4.4,5.42) circle (1.5pt);
\draw [fill=black] (-3.4,5.44) circle (1.5pt);
\draw [fill=black] (-3.66,5.44) circle (1.5pt);
\draw [fill=black] (-1.48,5.42) circle (0.5pt);
\draw [fill=black] (-1.32,5.42) circle (0.5pt);
\draw [fill=black] (-1.2,5.42) circle (0.5pt);
\draw [fill=black] (-1.92,5.44) circle (1.5pt);
\draw [fill=black] (-0.92,5.46) circle (1.5pt);
\draw [fill=black] (-1.66,5.42) circle (1.5pt);
\draw [fill=black] (-4.22,3.28) circle (0.5pt);
\draw [fill=black] (-4.06,3.28) circle (0.5pt);
\draw [fill=black] (-3.94,3.28) circle (0.5pt);
\draw [fill=black] (-4.44,3.28) circle (1.5pt);
\draw [fill=black] (-3.44,3.3) circle (1.5pt);
\draw [fill=black] (-3.7,3.3) circle (1.5pt);
\draw [fill=black] (-1.44,3.3) circle (0.5pt);
\draw [fill=black] (-1.28,3.3) circle (0.5pt);
\draw [fill=black] (-1.16,3.3) circle (0.5pt);
\draw [fill=black] (-1.9,3.28) circle (1.5pt);
\draw [fill=black] (-0.9,3.3) circle (1.5pt);
\draw [fill=black] (-1.64,3.28) circle (1.5pt);
\draw [fill=black] (0.76,5.44) circle (0.5pt);
\draw [fill=black] (0.92,5.44) circle (0.5pt);
\draw [fill=black] (1.04,5.44) circle (0.5pt);
\draw [fill=black] (0.54,5.44) circle (1.5pt);
\draw [fill=black] (1.54,5.46) circle (1.5pt);
\draw [fill=black] (1.28,5.46) circle (1.5pt);
\draw [fill=black] (3.46,5.44) circle (0.5pt);
\draw [fill=black] (3.62,5.44) circle (0.5pt);
\draw [fill=black] (3.74,5.44) circle (0.5pt);
\draw [fill=black] (3.02,5.46) circle (1.5pt);
\draw [fill=black] (4.02,5.48) circle (1.5pt);
\draw [fill=black] (3.28,5.44) circle (1.5pt);
\draw [fill=black] (0.72,3.3) circle (0.5pt);
\draw [fill=black] (0.88,3.3) circle (0.5pt);
\draw [fill=black] (1.0,3.3) circle (0.5pt);
\draw [fill=black] (0.5,3.3) circle (1.5pt);
\draw [fill=black] (1.5,3.32) circle (1.5pt);
\draw [fill=black] (1.24,3.32) circle (1.5pt);
\draw [fill=black] (3.52,3.32) circle (0.5pt);
\draw [fill=black] (3.68,3.32) circle (0.5pt);
\draw [fill=black] (3.8,3.32) circle (0.5pt);
\draw [fill=black] (3.04,3.3) circle (1.5pt);
\draw [fill=black] (4.04,3.32) circle (1.5pt);
\draw [fill=black] (3.3,3.3) circle (1.5pt);
\draw [fill=black] (5.66,5.42) circle (0.5pt);
\draw [fill=black] (5.82,5.42) circle (0.5pt);
\draw [fill=black] (5.94,5.42) circle (0.5pt);
\draw [fill=black] (5.44,5.42) circle (1.5pt);
\draw [fill=black] (6.44,5.44) circle (1.5pt);
\draw [fill=black] (6.2,5.44) circle (1.5pt);
\draw [fill=black] (8.36,5.42) circle (0.5pt);
\draw [fill=black] (8.52,5.42) circle (0.5pt);
\draw [fill=black] (8.64,5.42) circle (0.5pt);
\draw [fill=black] (7.92,5.44) circle (1.5pt);
\draw [fill=black] (8.92,5.46) circle (1.5pt);
\draw [fill=black] (8.18,5.42) circle (1.5pt);
\draw [fill=black] (5.62,3.28) circle (0.5pt);
\draw [fill=black] (5.78,3.28) circle (0.5pt);
\draw [fill=black] (5.9,3.28) circle (0.5pt);
\draw [fill=black] (5.4,3.28) circle (1.5pt);
\draw [fill=black] (6.4,3.3) circle (1.5pt);
\draw [fill=black] (6.14,3.3) circle (1.5pt);
\draw [fill=black] (8.42,3.32) circle (0.5pt);
\draw [fill=black] (8.58,3.32) circle (0.5pt);
\draw [fill=black] (8.7,3.32) circle (0.5pt);
\draw [fill=black] (7.96,3.28) circle (1.5pt);
\draw [fill=black] (8.94,3.3) circle (1.5pt);
\draw [fill=black] (8.2,3.28) circle (1.5pt);
\end{scriptsize}
\end{tikzpicture}

%% file: graph89.tex
\begin{tikzpicture}[line cap=round,line join=round,>=triangle 45,x=1.0cm,y=1.0cm]
\clip(-4.860000000000002,1.4600000000000026) rectangle (9.500000000000004,6.24);
\draw [rotate around={1.145762838175131:(-3.8999999999999995,5.429999999999998)}] (-3.8999999999999995,5.429999999999998) ellipse (0.7473599204960948cm and 0.5553799157009054cm);
\draw [rotate around={1.1457628381750795:(-1.4199999999999973,5.450000000000001)}] (-1.4199999999999973,5.450000000000001) ellipse (0.7473599204961262cm and 0.5553799157009286cm);
\draw [rotate around={1.145762838175105:(-3.9400000000000004,3.289999999999999)}] (-3.9400000000000004,3.289999999999999) ellipse (0.7473599204961175cm and 0.5553799157009224cm);
\draw [rotate around={1.145762838175105:(-1.4000000000000001,3.290000000000001)}] (-1.4000000000000001,3.290000000000001) ellipse (0.7473599204961234cm and 0.5553799157009267cm);
\draw (-1.92,5.44)-- (-3.4,5.44);
\draw (-3.4600000000000013,2.880000000000002) node[anchor=north west] {$R_{4n-4,27}$};
\draw (1.3600000000000005,2.880000000000002) node[anchor=north west] {$R_{4n-4,28}$};
\draw (6.360000000000002,2.840000000000002) node[anchor=north west] {$R_{4n-4,29}$};
\draw [rotate around={1.1457628381750795:(1.04,5.45)}] (1.04,5.45) ellipse (0.7473599204961185cm and 0.5553799157009228cm);
\draw [rotate around={1.1457628381751312:(3.5200000000000027,5.47)}] (3.5200000000000027,5.47) ellipse (0.7473599204961179cm and 0.5553799157009227cm);
\draw [rotate around={1.145762838175105:(0.9999999999999996,3.309999999999999)}] (0.9999999999999996,3.309999999999999) ellipse (0.747359920496112cm and 0.5553799157009182cm);
\draw [rotate around={1.1457628381751053:(3.5399999999999983,3.3099999999999987)}] (3.5399999999999983,3.3099999999999987) ellipse (0.7473599204961104cm and 0.5553799157009169cm);
\draw (3.02,5.46)-- (1.54,5.46);
\draw (3.28,5.44)-- (3.3,3.3);
\draw (1.28,5.46)-- (3.04,3.3);
\draw [rotate around={1.1457628381751304:(5.939999999999999,5.429999999999999)}] (5.939999999999999,5.429999999999999) ellipse (0.7473599204961031cm and 0.5553799157009114cm);
\draw [rotate around={1.1457628381750797:(8.42,5.449999999999999)}] (8.42,5.449999999999999) ellipse (0.7473599204961104cm and 0.5553799157009169cm);
\draw [rotate around={1.1457628381751053:(5.899999999999998,3.2899999999999996)}] (5.899999999999998,3.2899999999999996) ellipse (0.7473599204961027cm and 0.5553799157009112cm);
\draw [rotate around={1.1457628381751064:(8.439999999999996,3.289999999999999)}] (8.439999999999996,3.289999999999999) ellipse (0.747359920496079cm and 0.5553799157008941cm);
\draw (7.92,5.44)-- (6.44,5.44);
\draw (8.18,5.42)-- (8.2,3.28);
\draw (6.18,5.44)-- (7.94,3.28);
\draw (-3.66,5.44)-- (-3.44,3.3);
\draw (-3.44,3.3)-- (-1.9,3.28);
\draw (-1.64,3.28)-- (-1.92,5.44);
\draw (5.44,5.42)-- (5.4,3.28);
\begin{scriptsize}
\draw [fill=black] (-4.18,5.42) circle (0.5pt);
\draw [fill=black] (-4.02,5.42) circle (0.5pt);
\draw [fill=black] (-3.9,5.42) circle (0.5pt);
\draw [fill=black] (-4.4,5.42) circle (1.5pt);
\draw [fill=black] (-3.4,5.44) circle (1.5pt);
\draw [fill=black] (-3.66,5.44) circle (1.5pt);
\draw [fill=black] (-1.48,5.42) circle (0.5pt);
\draw [fill=black] (-1.32,5.42) circle (0.5pt);
\draw [fill=black] (-1.2,5.42) circle (0.5pt);
\draw [fill=black] (-1.92,5.44) circle (1.5pt);
\draw [fill=black] (-0.92,5.46) circle (1.5pt);
\draw [fill=black] (-1.66,5.42) circle (1.5pt);
\draw [fill=black] (-4.22,3.28) circle (0.5pt);
\draw [fill=black] (-4.06,3.28) circle (0.5pt);
\draw [fill=black] (-3.94,3.28) circle (0.5pt);
\draw [fill=black] (-4.44,3.28) circle (1.5pt);
\draw [fill=black] (-3.44,3.3) circle (1.5pt);
\draw [fill=black] (-3.7,3.3) circle (1.5pt);
\draw [fill=black] (-1.44,3.28) circle (0.5pt);
\draw [fill=black] (-1.28,3.28) circle (0.5pt);
\draw [fill=black] (-1.16,3.28) circle (0.5pt);
\draw [fill=black] (-1.9,3.28) circle (1.5pt);
\draw [fill=black] (-0.9,3.3) circle (1.5pt);
\draw [fill=black] (-1.64,3.28) circle (1.5pt);
\draw [fill=black] (0.76,5.44) circle (0.5pt);
\draw [fill=black] (0.92,5.44) circle (0.5pt);
\draw [fill=black] (1.04,5.44) circle (0.5pt);
\draw [fill=black] (0.54,5.44) circle (1.5pt);
\draw [fill=black] (1.54,5.46) circle (1.5pt);
\draw [fill=black] (1.28,5.46) circle (1.5pt);
\draw [fill=black] (3.46,5.44) circle (0.5pt);
\draw [fill=black] (3.62,5.44) circle (0.5pt);
\draw [fill=black] (3.74,5.44) circle (0.5pt);
\draw [fill=black] (3.02,5.46) circle (1.5pt);
\draw [fill=black] (4.02,5.48) circle (1.5pt);
\draw [fill=black] (3.28,5.44) circle (1.5pt);
\draw [fill=black] (0.72,3.3) circle (0.5pt);
\draw [fill=black] (0.88,3.3) circle (0.5pt);
\draw [fill=black] (1.0,3.3) circle (0.5pt);
\draw [fill=black] (0.5,3.3) circle (1.5pt);
\draw [fill=black] (1.5,3.32) circle (1.5pt);
\draw [fill=black] (1.24,3.32) circle (1.5pt);
\draw [fill=black] (3.5,3.3) circle (0.5pt);
\draw [fill=black] (3.66,3.3) circle (0.5pt);
\draw [fill=black] (3.78,3.3) circle (0.5pt);
\draw [fill=black] (3.04,3.3) circle (1.5pt);
\draw [fill=black] (4.04,3.32) circle (1.5pt);
\draw [fill=black] (3.3,3.3) circle (1.5pt);
\draw [fill=black] (5.66,5.42) circle (0.5pt);
\draw [fill=black] (5.82,5.42) circle (0.5pt);
\draw [fill=black] (5.94,5.42) circle (0.5pt);
\draw [fill=black] (5.44,5.42) circle (1.5pt);
\draw [fill=black] (6.44,5.44) circle (1.5pt);
\draw [fill=black] (6.18,5.44) circle (1.5pt);
\draw [fill=black] (8.36,5.42) circle (0.5pt);
\draw [fill=black] (8.52,5.42) circle (0.5pt);
\draw [fill=black] (8.64,5.42) circle (0.5pt);
\draw [fill=black] (7.92,5.44) circle (1.5pt);
\draw [fill=black] (8.92,5.46) circle (1.5pt);
\draw [fill=black] (8.18,5.42) circle (1.5pt);
\draw [fill=black] (5.62,3.28) circle (0.5pt);
\draw [fill=black] (5.78,3.28) circle (0.5pt);
\draw [fill=black] (5.9,3.28) circle (0.5pt);
\draw [fill=black] (5.4,3.28) circle (1.5pt);
\draw [fill=black] (6.4,3.3) circle (1.5pt);
\draw [fill=black] (6.14,3.3) circle (1.5pt);
\draw [fill=black] (8.42,3.28) circle (0.5pt);
\draw [fill=black] (8.58,3.28) circle (0.5pt);
\draw [fill=black] (8.7,3.28) circle (0.5pt);
\draw [fill=black] (7.94,3.28) circle (1.5pt);
\draw [fill=black] (8.94,3.3) circle (1.5pt);
\draw [fill=black] (8.2,3.28) circle (1.5pt);
\end{scriptsize}
\end{tikzpicture}

%% file: graph810.tex
\begin{tikzpicture}[line cap=round,line join=round,>=triangle 45,x=1.0cm,y=1.0cm]
\clip(-4.860000000000002,1.4600000000000026) rectangle (9.500000000000004,6.24);
\draw (1.3600000000000005,2.880000000000002) node[anchor=north west] {$R_{4n-4,30}$};
\draw [rotate around={1.1457628381750795:(1.04,5.45)}] (1.04,5.45) ellipse (0.7473599204961185cm and 0.5553799157009228cm);
\draw [rotate around={1.1457628381751312:(3.5200000000000027,5.47)}] (3.5200000000000027,5.47) ellipse (0.7473599204961179cm and 0.5553799157009227cm);
\draw [rotate around={1.145762838175105:(0.9999999999999996,3.309999999999999)}] (0.9999999999999996,3.309999999999999) ellipse (0.747359920496112cm and 0.5553799157009182cm);
\draw [rotate around={1.1457628381751053:(3.5399999999999983,3.3099999999999987)}] (3.5399999999999983,3.3099999999999987) ellipse (0.7473599204961104cm and 0.5553799157009169cm);
\draw (3.02,5.46)-- (1.54,5.46);
\draw (3.28,5.44)-- (3.3,3.3);
\draw (1.54,5.46)-- (1.5,3.32);
\draw (3.04,3.3)-- (1.5,3.32);
\begin{scriptsize}
\draw [fill=black] (0.76,5.44) circle (0.5pt);
\draw [fill=black] (0.92,5.44) circle (0.5pt);
\draw [fill=black] (1.04,5.44) circle (0.5pt);
\draw [fill=black] (0.54,5.44) circle (1.5pt);
\draw [fill=black] (1.54,5.46) circle (1.5pt);
\draw [fill=black] (1.28,5.46) circle (1.5pt);
\draw [fill=black] (3.46,5.44) circle (0.5pt);
\draw [fill=black] (3.62,5.44) circle (0.5pt);
\draw [fill=black] (3.74,5.44) circle (0.5pt);
\draw [fill=black] (3.02,5.46) circle (1.5pt);
\draw [fill=black] (4.02,5.48) circle (1.5pt);
\draw [fill=black] (3.28,5.44) circle (1.5pt);
\draw [fill=black] (0.72,3.3) circle (0.5pt);
\draw [fill=black] (0.88,3.3) circle (0.5pt);
\draw [fill=black] (1.0,3.3) circle (0.5pt);
\draw [fill=black] (0.5,3.3) circle (1.5pt);
\draw [fill=black] (1.5,3.32) circle (1.5pt);
\draw [fill=black] (1.24,3.32) circle (1.5pt);
\draw [fill=black] (3.54,3.3) circle (0.5pt);
\draw [fill=black] (3.7,3.3) circle (0.5pt);
\draw [fill=black] (3.82,3.3) circle (0.5pt);
\draw [fill=black] (3.04,3.3) circle (1.5pt);
\draw [fill=black] (4.04,3.32) circle (1.5pt);
\draw [fill=black] (3.3,3.3) circle (1.5pt);
\end{scriptsize}
\end{tikzpicture}

%% file: New_Star_critical_Ramsey_numbers_cycles_vs_K5_new.bbl
\begin{thebibliography}{99}

\bibitem{BoEr} J. A. Bondy and P. Erd{\"o}s, Ramsey number for cycles in graphs, {\it Journal of Combinatorial Theory Series B}, {\bf 14}, (1973), 46-54.
	
\bibitem{BaSuBr}	E.T. Baskoro,  H. Surahmat and H.J. Broersma,  The Ramsey numbers of fans versus $K_4$, {\it Bulletin of the Institute of Combinatorial Applications}, {\bf 43} (2005), 96-102.

\bibitem{BoJaSh}  B. Bollab{\'a}s, C. J. Jayawardene, Yang Jian Sheng, Huang Yi Ru, C. C. Rousseau and Zhang Ke Min, On a conjecture involving cycle-complete graph Ramsey numbers,  {\it The Australasian Journal},  {\bf 22}, (2000), 63-71.

\bibitem{HaMaSe} S. Haghi, H.R. Maimani and A. Seify, Star-critical Ramsey number of $F_n$ vs $K_4$, {\it Discrete Applied Mathematics}, {\bf 217: P2}, (2017), 203-209. 

\bibitem{Ho}	J. Hook, The Classification of Critical Graphs and Star-Critical Ramsey
Numbers, {\it Ph.D. Thesis, Lehigh University}, (2010).

\bibitem{HoIs}	J. Hook and G. Isaak, Star-critical Ramsey numbers, {\it Discrete Applied Mathematics}, {\bf 159}, (2011), 328-334.

\bibitem{JaNaRa} C.J. Jayawardene, D. Narvaez and S.P. Radziszowski, The star-critical Ramsey Number for any Cycle  vs. a $K_4$, {\it Discussiones Mathematicae Graph theory} {\it to appear}.

\bibitem{JaRo}  C. J. Jayawardene and C. C. Rousseau , The Ramsey number for quadrilateral vs. a complete graph on six vertices,  {\it Congressus Numerantium},  {\bf 123}, (1997), 97-108.

\bibitem{JaRo1} C.J. Jayawardene and C.C. Rousseau, The Ramsey Number for a Cycle of Length Five vs. a Complete Graph of Order Six, {\it Journal of Graph Theory}, {\bf 35} (2000) 99--108.

\bibitem{Ra} S.P. Radziszowski, Small Ramsey numbers, {\it Electronic Journal of Combinatorics} {\bf 14}, (2014), DS1.


\bibitem{YaYoRa} Wu Yali, Sun Yongqi and S.P. Radziszowski, Wheel and star-critical Ramsey numbers for quadrilaterals {\it Discrete Applied Mathematics}, {\bf 185}, (2015) 260-271.

\end{thebibliography}
